\documentclass[11pt,reqno]{amsart}

\usepackage{amsfonts,latexsym,amsthm,amssymb,amsmath,amscd,euscript}
\usepackage{framed}
\usepackage{fullpage}
\usepackage[margin=0.75in]{geometry}
\usepackage{hyperref}
    \hypersetup{colorlinks=true,citecolor=blue,urlcolor =gray,linkcolor=gray,linkbordercolor={1 0 0}}
\usepackage{tikz-cd}
\usepackage{comment}
\usepackage{mathtools}
\usepackage{mathdots}
\usepackage{enumerate}
\usepackage[new]{old-arrows}
\usepackage{mathrsfs}
\usepackage{pgfplots}
\usepackage{cancel}
\usepackage{rotating}
\usepackage{bookmark}
\usepackage{schemata}
\usepackage{environ}
\usepackage{graphics}
\usepackage{graphicx}
\usepackage{xfrac}
\usepackage{nicefrac}
\usepackage{tensor}
\usepackage{expdlist}
\usepackage{stmaryrd}
\usetikzlibrary{positioning}
\DeclareFontFamily{U}{mathx}{\hyphenchar\font45}
\DeclareFontShape{U}{mathx}{m}{n}{
      <5> <6> <7> <8> <9> <10>
      <10.95> <12> <14.4> <17.28> <20.74> <24.88>
      mathx10
      }{}
\DeclareSymbolFont{mathx}{U}{mathx}{m}{n}
\DeclareFontSubstitution{U}{mathx}{m}{n}
\DeclareMathAccent{\widecheck}{0}{mathx}{"71}
\DeclareMathAccent{\wideparen}{0}{mathx}{"75}

\makeatletter
\DeclareFontFamily{U}{tipa}{}
\DeclareFontShape{U}{tipa}{m}{n}{<->tipa10}{}
\newcommand{\arc@char}{{\usefont{U}{tipa}{m}{n}\symbol{62}}}%

\newcommand{\arc}[1]{\mathpalette\arc@arc{#1}}

\newcommand{\arc@arc}[2]{%
  \sbox0{$\m@th#1#2$}%
  \vbox{
    \hbox{\resizebox{\wd0}{\height}{\arc@char}}
    \nointerlineskip
    \box0
  }%
}
\makeatother

\makeatletter
\def\widebreve{\mathpalette\wide@breve}
\def\wide@breve#1#2{\sbox\z@{$#1#2$}%
     \mathop{\vbox{\m@th\ialign{##\crcr
\kern0.08em\brevefill#1{0.8\wd\z@}\crcr\noalign{\nointerlineskip}%
                    $\hss#1#2\hss$\crcr}}}\limits}
\def\brevefill#1#2{$\m@th\sbox\tw@{$#1($}%
  \hss\resizebox{#2}{\wd\tw@}{\rotatebox[origin=c]{90}{\upshape(}}\hss$}
\makeatletter

\usepackage{color}
\definecolor{blaze}{rgb}{0.00,0.420,0.00}

\DeclareRobustCommand\longtwoheadrightarrow{\relbar\joinrel\twoheadrightarrow}

\allowdisplaybreaks[1]

\NewEnviron{Thm*}[1][]{
	\vspace{10pt}\newline\noindent\ignorespaces\hspace*{-12.5pt}\DoBrackets\schema{}{\vspace{-7.5pt}\begin{theorem*}[#1]\BODY\newline\noindent$\square$\end{theorem*}\vspace{-7.5pt} }\ignorespacesafterend\leavevmode\vspace{12.5pt}\newline
}
\NewEnviron{Thm}[1][]{
	\vspace{10pt}\newline\noindent\ignorespaces\hspace*{-12.5pt}\DoBrackets\schema{}{\vspace{-7.5pt}\begin{theorem}[#1]\BODY\newline\noindent$\square$\end{theorem}\vspace{-7.5pt} }\ignorespacesafterend\leavevmode\vspace{12.5pt}\newline
}
\NewEnviron{Prop*}[1][]{
	\vspace{10pt}\newline\noindent\ignorespaces\hspace*{-12.5pt}\DoBrackets\schema{}{\vspace{-7.5pt}\begin{proposition*}[#1]\BODY\newline\noindent$\square$\end{proposition*}\vspace{-7.5pt} }\ignorespacesafterend\leavevmode\vspace{12.5pt}\newline
}
\NewEnviron{Prop}[1][]{
	\vspace{10pt}\newline\noindent\ignorespaces\hspace*{-12.5pt}\DoBrackets\schema{}{\vspace{-7.5pt}\begin{proposition}[#1]\BODY\newline\noindent$\square$\end{proposition}\vspace{-7.5pt} }\ignorespacesafterend\leavevmode\vspace{12.5pt}\newline
}
\NewEnviron{Def*}[1][]{
	\vspace{10pt}\newline\noindent\ignorespaces\hspace*{-7.5pt}\DoParens\schema{}{\vspace{-7.5pt}\begin{definition*}[#1]\BODY\end{definition*}\vspace{-7.5pt} }\ignorespacesafterend\leavevmode\vspace{12.5pt}\newline
}
\NewEnviron{Def}[1][]{
	\vspace{5pt}\newline\noindent\ignorespaces\hspace*{-7.5pt}\DoParens\schema{}{\vspace{-7.5pt}\begin{definition}[#1]\BODY\end{definition}\vspace{-7.5pt} }\ignorespacesafterend\leavevmode\vspace{7.5pt}\newline
}
\NewEnviron{Lem*}[1][]{
	\vspace{10pt}\newline\noindent\ignorespaces\hspace*{-12.5pt}\DoBrackets\schema{}{\vspace{-7.5pt}\begin{lemma*}[#1]\BODY\newline\noindent$\square$\end{lemma*}\vspace{-7.5pt} }\ignorespacesafterend\leavevmode\vspace{12.5pt}\newline
}
\NewEnviron{Lem}[1][]{
	\vspace{10pt}\newline\noindent\ignorespaces\hspace*{-12.5pt}\DoBrackets\schema{}{\vspace{-7.5pt}\begin{lemma}[#1]\BODY\newline\noindent$\square$\end{lemma}\vspace{-7.5pt} }\ignorespacesafterend\leavevmode\vspace{12.5pt}\newline
}
\NewEnviron{Fact}[1][]{
	\vspace{10pt}\newline\noindent\ignorespaces\hspace*{-12.5pt}\DoBrackets\schema{}{\vspace{-7.5pt}\begin{fact}[#1]\BODY\newline\noindent$\square$\end{fact}\vspace{-7.5pt} }\ignorespacesafterend\leavevmode\vspace{12.5pt}\newline
}
\NewEnviron{Fact*}[1][]{
	\vspace{10pt}\newline\noindent\ignorespaces\hspace*{-12.5pt}\DoBrackets\schema{}{\vspace{-7.5pt}\begin{fact*}[#1]\BODY\newline\noindent$\square$\end{fact*}\vspace{-7.5pt} }\ignorespacesafterend\leavevmode\vspace{12.5pt}\newline
}
\NewEnviron{THM*}[1][]{
	\vspace{10pt}\newline\noindent\ignorespaces\hspace*{-12.5pt}\DoBrackets\schema{}{\vspace{-7.5pt}\begin{theorem*}[#1]\BODY\end{theorem*}\vspace{-6pt} }\ignorespacesafterend\leavevmode\vspace{12.5pt}\newline
}
\NewEnviron{FACT*}[1][]{
	\vspace{10pt}\newline\noindent\ignorespaces\hspace*{-12.5pt}\DoBrackets\schema{}{\vspace{-7.5pt}\begin{fact*}[#1]\BODY\end{fact*}\vspace{-6pt} }\ignorespacesafterend\leavevmode\vspace{12.5pt}\newline
}
\NewEnviron{LEM*}[1][]{
	\vspace{10pt}\newline\noindent\ignorespaces\hspace*{-12.5pt}\DoBrackets\schema{}{\vspace{-7.5pt}\begin{lemma*}[#1]\BODY\end{lemma*}\vspace{-6pt} }\ignorespacesafterend\leavevmode\vspace{12.5pt}\newline
}
\NewEnviron{PROP*}[1][]{
	\vspace{10pt}\newline\noindent\ignorespaces\hspace*{-12.5pt}\DoBrackets\schema{}{\vspace{-7.5pt}\begin{proposition*}[#1]\BODY\end{proposition*}\vspace{-6pt} }\ignorespacesafterend\leavevmode\vspace{12.5pt}\newline
}
\NewEnviron{COR*}[1][]{
	\vspace{10pt}\newline\noindent\ignorespaces\hspace*{-12.5pt}\DoBrackets\schema{}{\vspace{-7.5pt}\begin{corollary*}[#1]\BODY\end{corollary*}\vspace{-6pt} }\ignorespacesafterend\leavevmode\vspace{12.5pt}\newline
}
\NewEnviron{THM}[1][]{
	\vspace{10pt}\newline\noindent\ignorespaces\hspace*{-12.5pt}\DoBrackets\addtocounter{theorem}{-1}\schema{}{\vspace{-7.5pt}\begin{theorem}[#1]\BODY\end{theorem}\vspace{-6pt} }\ignorespacesafterend\leavevmode\vspace{12.5pt}\newline
}
\NewEnviron{FACT}[1][]{
	\vspace{10pt}\newline\noindent\ignorespaces\hspace*{-12.5pt}\DoBrackets\addtocounter{theorem}{-1}\schema{}{\vspace{-7.5pt}\begin{fact}[#1]\BODY\end{fact}\vspace{-6pt} }\ignorespacesafterend\leavevmode\vspace{12.5pt}\newline
}
\NewEnviron{LEM}[1][]{
	\vspace{10pt}\newline\noindent\ignorespaces\hspace*{-12.5pt}\DoBrackets\addtocounter{theorem}{-1}\schema{}{\vspace{-7.5pt}\begin{lemma}[#1]\BODY\end{lemma}\vspace{-6pt} }\ignorespacesafterend\leavevmode\vspace{12.5pt}\newline
}
\NewEnviron{PROP}[1][]{
	\vspace{10pt}\newline\noindent\ignorespaces\hspace*{-12.5pt}\DoBrackets\addtocounter{theorem}{-1}\schema{}{\vspace{-7.5pt}\begin{proposition}[#1]\BODY\end{proposition}\vspace{-6pt} }\ignorespacesafterend\leavevmode\vspace{12.5pt}\newline
}
\NewEnviron{COR}[1][]{
	\vspace{10pt}\newline\noindent\ignorespaces\hspace*{-12.5pt}\DoBrackets\addtocounter{theorem}{-1}\schema{}{\vspace{-7.5pt}\begin{corollary}[#1]\BODY\end{corollary}\vspace{-6pt} }\ignorespacesafterend\leavevmode\vspace{12.5pt}\newline
}
\NewEnviron{EX*}[1][]{
    \vspace{-1em}\begin{adjustwidth}{2em}{0pt}
    \begin{example*}
    \BODY
    \end{example*}
    \end{adjustwidth}
}
\NewEnviron{PRF}[1][]{
    \begin{adjustwidth}{2em}{0pt}
    \begin{proof}
    \BODY
    \end{proof}
    \end{adjustwidth}\vspace{0.5em}
}
\NewEnviron{DEF}[1][]{
	\vspace{0.5em}\newline\noindent\ignorespaces\hspace*{-7.5pt}\DoParens\addtocounter{theorem}{-1}\schema{}{\vspace{-7.5pt}\begin{definition}[#1]\BODY\end{definition}\vspace{-7.5pt} }\ignorespacesafterend\leavevmode\vspace{0.6em}\newline
}
\NewEnviron{PRFNAME}[1][]{
    \vspace{0.5em}\begin{adjustwidth}{2em}{0pt}
    \begin{proof}[#1]
    \BODY
    \end{proof}
    \end{adjustwidth}\vspace{0.5em}
}

\theoremstyle{definition}
\newtheorem{theorem}{Theorem}[section]
\newtheorem*{theorem*}{Theorem}%
\newtheorem*{proposition*}{Proposition}
\newtheorem{proposition}[theorem]{Proposition}
\newtheorem{lemma}[theorem]{Lemma}
\newtheorem*{lemma*}{Lemma}
\newtheorem{fact}[theorem]{Fact}
\newtheorem*{fact*}{Fact}
\newtheorem{corollary}[theorem]{Corollary}

\theoremstyle{definition}
\newtheorem{definition}[theorem]{Definition}
\newtheorem*{definition*}{Definition}

\newtheorem*{claim*}{Claim}
\newtheorem*{corollary*}{Corollary}
\newtheorem*{example*}{Example}

\theoremstyle{remark}
\newtheorem*{remark}{Remark}

\newcommand{\BC}{\mathbb C}

\newcommand{\BZ}{\mathbb Z}
\newcommand{\BN}{\mathbb N}

\newcommand{\id}{\textnormal{id}}
\renewcommand{\d}{\textnormal{d}}
\DeclareMathOperator{\img}{Img}

\newcommand{\nc}{\newcommand}

\nc{\on}{\operatorname}
\nc{\spec}{\on{Spec}}

\DeclareMathOperator{\Ker}{Ker}
\renewcommand{\ker}{\Ker}
\DeclareMathOperator{\Hom}{Hom}
\renewcommand{\hom}{\Hom}

\DeclarePairedDelimiter\pr{(}{)}

\DeclareMathOperator{\Spec}{Spec}
\renewcommand{\spec}{\Spec}

\newcommand{\ol}[1]{\overline{#1}}

\DeclareMathOperator{\Sym}{Sym}

\DeclareMathOperator{\rad}{Rad}

\newcommand{\pidl}{\mathfrak p}

\nc{\wt}[1]{\widetilde{#1}}
\nc{\wh}[1]{\widehat{#1}}
\nc{\pd}{\partial}
\nc{\cal}[1]{\mathcal{#1}}
\renewcommand{\frak}[1]{\mathfrak{#1}}
\nc{\scr}[1]{\mathscr{#1}}
\nc{\tsubseteq}[1]{\overset{\textnormal{#1}}{\subseteq}}
\nc{\tsupseteq}[1]{\overset{\textnormal{#1}}{\supseteq}}
\nc{\tsubset}[1]{\overset{\textnormal{#1}}{\subset}}
\nc{\tsupset}[1]{\overset{\textnormal{#1}}{\supset}}
\nc{\vphi}{\varphi}
\nc{\beau}{\displaystyle}
\nc{\existss}{\exists\;}
\nc{\foralls}{\forall\;}

\DeclarePairedDelimiter\set{\{}{\}}
\nc{\ov}[1]{\overrightarrow{#1}}
\nc{\cdotsc}{\cdots\hspace{-0.1pt}}

\nc{\D}{\textnormal{D}}

\nc{\tn}[1]{\textnormal{#1}}
\nc{\lto}{\longrightarrow}
\nc{\lmto}{\longmapsto}

\DeclareMathOperator{\gal}{Gal}

\nc{\sm}{\setminus}

\nc{\scomp}{\textsf{c}}
\DeclareMathOperator{\cnt}{cnt}

\nc{\ncnt}{\lnot\cnt}
\nc{\ucong}{\overset{!}{\cong}}
\nc{\tcong}[1]{\overset{#1}{\cong}}
\nc{\coker}{\on{Coker}}

\nc{\eps}{\varepsilon}
\nc{\BS}{\mathbb{S}}
\nc{\BT}{\mathbb{T}}
\nc{\SC}{\mathscr{C}}

\DeclareMathOperator{\tor}{Tor}
\nc{\wa}[1]{\wideparen{#1}}
\nc{\qidl}{\frak{q}}
\nc{\nidl}{\frak n}

\nc{\te}[1]{\text{#1}}
\nc{\imps}{\quad\quad\ \:\,}
\nc{\comp}{\complement}
\nc{\van}{\cal V}
\nc{\ide}{\cal I}
\nc{\nfrac}[2]{\nicefrac{#1}{#2}}
\nc{\snfrac}[2]{\scalebox{1.45}{$\nfrac{#1}{#2}$}}
\nc{\inj}{\hookrightarrow}
\nc{\surj}{\twoheadrightarrow}
\nc{\linj}{\longhookrightarrow}
\nc{\lsurj}{\longtwoheadrightarrow}
\nc{\four}{\cal F}
\nc{\urot}[1]{\mathbin{\rotatebox[origin=c]{90}{$#1$}}}
\nc{\drot}[1]{\mathbin{\rotatebox[origin=c]{-90}{$#1$}}}
\nc{\vrot}[2]{\mathbin{\rotatebox[origin=c]{#2}{$#1$}}}

\nc{\wb}[1]{\widebreve{#1}}
\DeclareMathOperator{\Id}{Id}
\nc{\BH}{\mathbb{H}}

\nc{\T}{\textnormal{T}}
\nc{\N}{\vrot{\T}{180}\!}
\nc{\blt}{\bullet}
\nc{\lbd}{\lambda}

\renewcommand{\tt}[1]{\texttt{#1}}
\renewcommand{\sf}[1]{\mathsf{#1}}
\nc{\sq}{{\mathop{\square}}}
\nc{\onl}[1]{\operatorname*{#1}}
\nc{\bigger}[1]{\scalebox{1.5}{$#1$}}
\nc{\tbt}[4]{\begin{pmatrix}#1&#2\\#3&#4\end{pmatrix}}
\nc{\thbth}[9]{\begin{pmatrix}#1&#2&#3\\#4&#5&#6\\#7&#8&#9\end{pmatrix}}
\nc{\tbo}[2]{\begin{pmatrix}#1\\#2\end{pmatrix}}
\nc{\thbo}[3]{\begin{pmatrix}#1\\#2\\#3\end{pmatrix}}
\nc{\tsq}[3]{\tensor{{#1}}{^{#2}_{#3}}}
\nc{\rnk}{\on{rnk}}
\nc{\Ad}{\on{Ad}}
\nc{\ad}{\on{ad}}
\nc{\glie}{\frak{g}}
\nc{\SO}{\on{SO}}
\nc{\sk}{\on{Sk}}
\nc{\snakeanchor}{\ar[draw=none]{d}[name=X, anchor=center]{}}
\nc{\snakearrow}[1]{\ar[rounded corners,
	to path={ -- ([xshift=2ex]\tikztostart.east)
		|- (X.center) \tikztonodes
		-| ([xshift=-2ex]\tikztotarget.west)
		-- (\tikztotarget)}]{dll}{#1}}
\nc{\cosnakearrow}[1]{\ar[rounded corners,
	to path={ -- ([xshift=-2ex]\tikztostart.west)
		|- (X.center) \tikztonodes
		-| ([xshift=2ex]\tikztotarget.east)
		-- (\tikztotarget)}]{urr}{#1}}
\nc{\colim}{\on*{colim}}
\nc{\sing}{\on{Sing}}
\makeatletter
\newcommand{\dircolim@}[2]{%
  \vtop{\m@th\ialign{##\cr
    \hfil$#1\operator@font colim$\hfil\cr
    \noalign{\nointerlineskip\kern1.5\ex@}#2\cr
    \noalign{\nointerlineskip\kern-\ex@}\cr}}%
}
\newcommand{\dircolim}{%
  \mathop{\mathpalette\dircolim@{\rightarrowfill@\textstyle}}\nmlimits@
}
\makeatother
\makeatletter
\newcommand{\xdashrightarrow}[2][]{\ext@arrow 0359\rightarrowfill@@{#1}{#2}}
\newcommand{\xdashleftarrow}[2][]{\ext@arrow 3095\leftarrowfill@@{#1}{#2}}
\newcommand{\xdashleftrightarrow}[2][]{\ext@arrow 3359\leftrightarrowfill@@{#1}{#2}}
\def\rightarrowfill@@{\arrowfill@@\relax\relbar\rightarrow}
\def\leftarrowfill@@{\arrowfill@@\leftarrow\relbar\relax}
\def\leftrightarrowfill@@{\arrowfill@@\leftarrow\relbar\rightarrow}
\def\arrowfill@@#1#2#3#4{%
  $\m@th\thickmuskip0mu\medmuskip\thickmuskip\thinmuskip\thickmuskip
   \relax#4#1
   \xleaders\hbox{$#4#2$}\hfill
   #3$%
}
\makeatother
\nc{\lcm}{\on{lcm}}
\nc{\Stab}{\on{Stab}}
\nc{\Nbhd}{\on{Nbhd}}
\nc{\Lie}{\on{Lie}}
\nc{\gl}{\on{\mathfrak{gl}}}
\nc{\Gr}{\on{Gr}}
\nc{\act}{\on{act}}
\nc{\ev}{\on{ev}}
\nc{\Sq}{\on{Sq}}
\nc{\ext}{\on{Ext}}
\nc{\Ext}{\on{Ext}}
\nc{\Dist}{\on{Dist}}
\nc{\dlie}{\frak{d}}
\nc{\hlie}{\frak{h}}
\nc{\pt}{\te{pt}}
\nc{\qedcheck}{\vspace{-1em}\begin{flushright}$\checkmark\quad$\end{flushright}\vspace{-0.5em}}
\nc{\lemproof}[1]{\smallskip\noindent\underline{#1:}}
\nc{\F}{\te{F}}
\nc{\jota}{\jmath}
\nc{\Tor}{\tor}
\nc{\Fl}{\on{Fl}}
\nc{\prob}{\on{prob}}
\nc{\didl}{\frak{d}}
\nc{\Rep}{\on{\sf{Rep}}}
\nc{\Irrep}{\on{\sf{irRep}}}
\nc{\Res}{\on{Res}}
\nc{\Ind}{\on{Ind}}
\nc{\Der}{\on{Der}}
\renewcommand{\sl}{\on{\frak{sl}}}
\nc{\so}{\on{\frak{so}}}

\nc{\Img}{\img}
\nc{\Sidl}{\frak{S}}
\nc{\Slie}{\frak{S}}
\nc{\iita}{\imath}
\nc{\hidl}{\hlie}
\nc{\blie}{\frak{b}}
\nc{\bidl}{\frak{b}}
\nc{\nlie}{\nidl}
\nc{\larr}{\leftarrow}
\nc{\darr}{\downarrow}
\nc{\nwarr}{\nwarrow} 
\nc{\searr}{\searrow} 
\nc{\Tot}{\on{Tot}}
\nc{\Wt}{\on{Wt}}
\nc{\sfcong}[1]{\overset{\sf{#1}}{\cong}}
\nc{\vtheta}{\vartheta}
\nc{\JH}{\on{JH}}
\nc{\Typ}{\on{Typ}}
\nc{\alie}{\frak{a}}

\nc{\llra}{\longleftrightarrow}
\nc{\Llra}{\Longleftrightarrow} 
\nc{\lla}{\longleftarrow} 
\nc{\lra}{\longrightarrow} 
\nc{\Lra}{\Longrightarrow}  
\newcommand\sqplus{\mathbin{\ooalign{$\sqcup$\cr%
   \hfil\raise0.42ex\hbox{$\scriptscriptstyle+$}\hfil\cr}}}
\newcommand\sqminus{\mathbin{\ooalign{$\sqcup$\cr%
   \hfil\raise0.42ex\hbox{$\scriptscriptstyle-$}\hfil\cr}}}
\newcommand\uminus{\mathbin{\ooalign{$\cup$\cr%
   \hfil\raise0.42ex\hbox{$\scriptscriptstyle-$}\hfil\cr}}}
   
\DeclarePairedDelimiter\biggpr{\bigg(}{\bigg)}

\nc{\CO}{\cal O}
\nc{\nor}{\on{nor}}
\nc{\CP}{\cal P}
\nc{\CF}{\cal F}
\nc{\CG}{\cal G}
\nc{\sh}{\te{sh}}
\nc{\Supp}{\on{Supp}}
\nc{\Sh}{\on{\sf{Sh}}}
\nc{\rv}{\rvert}
\nc{\bigrv}{\big\rvert}
\nc{\Bigrv}{\Big\rvert}
\nc{\biggrv}{\bigg\rvert}
\nc{\Biggrv}{\Bigg\rvert}
\nc{\Cl}{\on{Cl}}
\nc{\existsu}{\exists!\;}
\nc{\disc}{\on{disc}}
\nc{\Proj}{\on{Proj}}
\nc{\op}{\te{\,op}}
\nc{\actson}{\mathrel{\vrot{\circlearrowright}{270}}}
\nc{\Mod}{\sf{Mod}}
\nc{\Grp}{\sf{Grp}}
\nc{\Set}{\sf{Set}}
\nc{\Alg}{\sf{Alg}}

\nc{\Ring}{\sf{Ring}}
\nc{\Ann}{\on{Ann}}
\DeclarePairedDelimiter\wan{\langle}{\rangle}
\nc{\osubseteq}{\overset{\te{o}}{\subseteq}}
\nc{\osubset}{\overset{\te{o}}{\subset}}
\nc{\osupseteq}{\overset{\te{o}}{\supseteq}}
\nc{\osupset}{\overset{\te{o}}{\supset}}
\nc{\CV}{\cal V}
\nc{\CI}{\cal I}
\nc{\glu}{\on{glu}}
\nc{\jidl}{\frak j}
\nc{\Span}{\on{Span}}

\nc{\oscong}[1]{\overset{#1}{\cong}}
\nc{\ossubseteq}[1]{\overset{#1}{\subseteq}}
\nc{\ossupseteq}[1]{\overset{#1}{\supseteq}}
\nc{\ossubset}[1]{\overset{#1}{\subset}}
\nc{\ossupset}[1]{\overset{#1}{\supset}}
\nc{\Gal}{\gal}
\nc{\uad}{\:\;}
\nc{\CQ}{\cal Q}
\nc{\stspace}{\ \ \,}
\nc{\ractson}{\mathrel{\vrot{\circlearrowleft}{90}}}
\nc{\len}{\on{len}}
\nc{\CC}{\cal C}
\nc{\bsfrac}[2]{\scalebox{1.4}[1.5]{$\sfrac{#1}{#2}$}}
\newcommand\quot[2]{
        \mathchoice
            {
                \text{\raise1ex\hbox{$#1$}\Big/\lower1ex\hbox{$#2$}}%
            }
            {
                #1\,/\,#2
            }
            {
                #1\,/\,#2
            }
            {
                #1\,/\,#2
            }
    }
\nc{\acts}{\actson}
\nc{\racts}{\ractson}
\nc{\calD}{\cal D}
\nc{\Fib}{\on{Fib}}
\nc{\cartsq}{\arrow[dr, phantom, "\ulcorner", very near start]}
\nc{\Glu}{\on*{Glu}}
\nc{\bnfrac}[2]{\raisebox{-0.175em}{\scalebox{1.4}[1.5]{$\nfrac{#1}{#2}$}}}
\nc{\mfrac}[2]{\raisebox{-0.15em}{\scalebox{1.32}[1.4]{$\nfrac{#1}{#2}$}}} 
\nc{\Sch}{\sf{Sch}}
\nc{\tlto}[1]{\overset{\te{#1}}{\lto}}
\nc{\tlmto}[1]{\overset{\te{#1}}{\lmto}}
\nc{\tlsurj}[1]{\overset{\te{#1}}{\lsurj}}
\nc{\tlinj}[1]{\overset{\te{#1}}{\linj}}
\nc{\oslto}[1]{\overset{#1}{\lto}}
\nc{\oslmto}[1]{\overset{#1}{\lmto}}
\nc{\oslsurj}[1]{\overset{#1}{\lsurj}}
\nc{\oslinj}[1]{\overset{#1}{\linj}}
\DeclarePairedDelimiter\pparen{(\!(}{)\!)}
\nc{\cornk}{\on{cornk}}
\nc{\ute}[1]{\uad\te{#1}}
\nc{\CL}{\cal L}
\nc{\val}{\on{val}}
\nc{\warn}[1]{\textcolor{red}{#1}}
\nc{\tsp}{\ \ }
\nc{\tte}[1]{\tsp\te{#1}}
\nc{\QCoh}{\sf{QCoh}}
\nc{\qcoh}{\sf{QCoh}}
\nc{\onqcoh}{\on{\qcoh}}
\renewcommand{\sh}{\on{Sh}}
\nc{\Weil}{\on{Weil}}
\nc{\vrho}{\varrho}
\nc{\cD}{\cal D}
\nc{\shom}{\on{\cal H \te{om}}}
\nc{\Rho}{\te{P}}
\nc{\CM}{\cal M}
\nc{\Pic}{\on{Pic}}
\nc{\sep}{\te{sep}}
\nc{\BG}{\mathbb{G}}
\nc{\BGm}{{\BG_{\te m}}}

\nc{\Lbd}{\Lambda}
\nc{\Cone}{\on{Cone}}
\nc{\cone}{\Cone}
\nc{\Cd}{\cal D}
\nc{\CH}{\cal H}
\nc{\CX}{\cal X}
\nc{\heart}{\heartsuit}
\nc{\CA}{\cal A}
\nc{\BGa}{\BG_\te a}
\nc{\Lla}{\Longleftarrow}
\nc{\suptrieq}{\mathbin{\unrhd}}
\nc{\R}{\te{R}}

\nc{\rhom}{\on{RHom}}
\nc{\lotimes}{\mathbin{\overset{\te{L}}{\otimes}}}
\nc{\bdd}{\te{bdd}}
\nc{\twt}{\on{Wt}}
\nc{\weit}{\on{Wt}}
\nc{\typ}{\on{Typ}}
\nc{\fun}{\on{\sf{Fun}}}
\nc{\fib}{\on{Fib}}
\nc{\Cofib}{\on{Cofib}}
\nc{\cofib}{\on{Cofib}}
\nc{\calid}{\cal{I}\te{d}}
\nc{\ladj}{{\vrot{\vdash}{90}}}
\nc{\radj}{{\vrot{\dashv}{90}}}
\nc{\Ch}{\on{\sf{Ch}}}
\nc{\ohlie}{\ol{\hlie^*}}
\nc{\Lan}{\on{Lan}}
\nc{\Ran}{\on{Ran}}
\nc{\CB}{\cal B}
\makeatletter
\newcommand{\subalign}[1]{%
  \vcenter{%
    \Let@ \restore@math@cr \default@tag
    \baselineskip\fontdimen10 \scriptfont\tw@
    \advance\baselineskip\fontdimen12 \scriptfont\tw@
    \lineskip\thr@@\fontdimen8 \scriptfont\thr@@
    \lineskiplimit\lineskip
    \ialign{\hfil$\m@th\scriptstyle##$&$\m@th\scriptstyle{}##$\hfil\crcr
      #1\crcr
    }%
  }%
}
\makeatother
\nc{\veps}{\epsilon}
\nc{\aut}{\on{Aut}}
\nc{\saut}{\on{\cal{A}\te{ut}}}
\nc{\Fgt}{\on{Fgt}}
\nc{\fgt}{\on{Fgt}}
\nc{\ul}[1]{\underline{#1}}
\nc{\SL}{\on{SL}}
\nc{\radu}{\rad_\te{u}}
\nc{\tlie}{\frak{t}}
\nc{\Rad}{\on{Rad}}
\nc{\CE}{\cal E}

\nc{\CU}{\cal U}
\nc{\sproj}{\on{\cal P\te{roj}}}
\nc{\chom}{\shom}
\nc{\cproj}{\sproj}
\nc{\ao}{\te{ao}}
\nc{\CProj}{\cproj}
\nc{\gari}{{g_\te{ari}}}
\nc{\Coh}{\on{\sf{Coh}}}
\nc{\red}{\te{red}}
\nc{\Bl}{\on{Bl}}
\nc{\Hilb}{\on{Hilb}}
\nc{\FX}{\frak{X}}
\usepackage{CJKutf8}
\nc{\AlgGrp}{\sf{AlgGrp}}
\nc{\subnoreq}{\trianglelefteq}
\nc{\subnor}{\triangleleft}
\nc{\supnoreq}{\trianglerighteq}
\nc{\supnor}{\triangleright}
\nc{\Var}{\sf{Var}}
\renewcommand{\D}{\sf D}
\nc{\longlongline}{\noindent ---------------------------------------------------------------------------------------------------------------------------------------------}
\nc{\vareps}{\epsilon}
\usepackage{array}

\makeatletter
\newcommand{\dirlim@}[2]{%
  \vtop{\m@th\ialign{##\cr
    \hfil$#1\operator@font lim$\hfil\cr
    \noalign{\nointerlineskip\kern1.5\ex@}#2\cr
    \noalign{\nointerlineskip\kern-\ex@}\cr}}%
}
\newcommand{\dirlim}{%
  \mathop{\mathpalette\dirlim@{\leftarrowfill@\textstyle}}\nmlimits@
}
\makeatother
\makeatletter

\renewcommand{\R}{\sf{R}}
\renewcommand{\rhom}{\on{\sf{R}\te{Hom}}}
\renewcommand{\lotimes}{\mathbin{\overset{\sf{L}}{\otimes}}}
\nc{\Ass}{\on{Ass}}

\usepackage{changepage}
\nc{\midwedge}{{\textstyle\bigwedge}}
\renewcommand{\calid}{\on{\cal{I}\te{d}}}

\usepackage{accents}

\renewcommand{\shom}{\on{\mathcal H\hspace*{-0.5pt}\textsl{om}}}

\usepackage{bbm}
\newcommand{\bk}{\mathbbm{k}}
\makeatletter 
\def\l@subsection{\@tocline{2}{0pt}{2pc}{6pc}{}} 
\makeatother

\newcommand{\from}{\leftarrow}
\renewcommand{\acts}{{\actson}}
\renewcommand{\racts}{{\ractson}}

\usetikzlibrary{decorations.markings}

\makeatletter
\tikzcdset{
  open/.code     = {\tikzcdset{hook, circled};},
  closed/.code   = {\tikzcdset{hook, slashed};},
  open'/.code    = {\tikzcdset{hook', circled};},
  closed'/.code  = {\tikzcdset{hook', slashed};},
  circled/.code  = {\tikzcdset{markwith = {\draw (0,0) circle (.375ex);}};},
  slashed/.code  = {\tikzcdset{markwith = {\draw[-] (-.4ex,-.4ex) -- (.4ex,.4ex);}};},
  markwith/.code ={
    \pgfutil@ifundefined%
    {tikz@library@decorations.markings@loaded}%
    {\pgfutil@packageerror{tikz-cd}{You need to say %
      \string\usetikzlibrary{decorations.markings} to use arrows with markings}{}}{}%
    \pgfkeysalso{/tikz/postaction = {
      /tikz/decorate,
      /tikz/decoration={markings, mark = at position 0.5 with {#1}}}
    }
  },
}
\makeatother

\DeclareFontFamily{OT1}{pzc}{}
\DeclareFontShape{OT1}{pzc}{m}{it}%
{<-> s * [1.15] pzcmi7t}{}
\DeclareMathAlphabet{\mathpzc}{OT1}{pzc}{m}{it} 

\usepackage{ytableau}
\ytableausetup{smalltableaux}

\newcommand{\simlto}{\oslto\sim}

\nc{\lmfrac}[2]{\raisebox{-0.15em}{\scalebox{1.32}[1.4]{$\reflectbox{\nfrac{\reflectbox{\ensuremath{#1}}}{\reflectbox{\ensuremath{#2}}}}$}}} 
\nc{\lrfrac}[3]{\raisebox{0em}{{$\lmfrac{#2}{#1}\hspace{-0.4em}\mfrac{\ }{#3}$}}}

\newcommand{\Hed}{\on{Hed}}

\DeclareMathAlphabet{\mathdutchcal}{U}{dutchcal}{m}{n}
\SetMathAlphabet{\mathdutchcal}{bold}{U}{dutchcal}{b}{n}
\DeclareMathAlphabet{\mathdutchbcal}{U}{dutchcal}{b}{n}

\renewcommand{\midwedge}{\raisebox{0.15em}{${\textstyle\bigwedge}$}}

\newcommand{\Db}{\sf D^{\te{b}}}
\nc{\BK}{\mathbb K}

\makeatletter 
\def\l@subsection{\@tocline{2}{0pt}{2pc}{6pc}{}} 
\makeatother 

\nc{\Nabla}{\nabla}
\nc{\sic}{\sf{St}\infty\sf{Cat}}
\nc{\gr}{\on{gr}}
\nc{\mr}[1]{\mathring{#1}}
\nc{\tsl}[1]{\textsl{#1}}
\nc{\sym}{\te{sym}}
\renewcommand{\ext}{\on{ext}}
\nc{\nat}{\natural}
\nc{\specseqimplies}{\uad\implies\uad}

\nc{\tX}{\te X}\nc{\tH}{\te H}\nc{\tY}{\te Y}
\nc{\NB}{{\te{N}\CB}}
\nc{\paperprogressmarker}{\warn{PAPER\ WRITEUP\ PROGRESS\ MARKER\ HERE}}
\nc{\RHom}{\rhom}
\nc{\kidl}{\frak{k}}
\nc{\klie}{\frak{k}}
\newcommand{\longsquigglyrightarrow}{\xymatrix{{}\ar@{~>}[r]&{}}} 

\input xy
\usepackage[all]{xy}
\xyoption{line}
\xyoption{arrow}
\xyoption{color}
\SelectTips{cm}{}
\newcommand{\hackcenter}[1]{
 \xy (0,0)*{#1}; \endxy}

\newcommand{\tikzcupcap}{\hackcenter{\begin{tikzpicture}[scale=0.375]
    \draw (0,0) to[out=60,in=120] (2,0);
    \draw (0,2) to[out=-60,in=-120] (2,2);
\end{tikzpicture}}}

\newcommand{\tikzthrough}{\hackcenter{\begin{tikzpicture}[scale=0.375]
    \draw (0,0) to[out=60,in=-60] (0,2);
    \draw (2,0) to[out=120,in=-120] (2,2);
\end{tikzpicture}}}

\title{BGG Resolutions, Koszulity, and Stratifications, Part I:\\ the nilBrauer Algebra} 
\pgfplotsset{compat=1.14}
\begin{document}

\maketitle

\vspace*{-0.25in}
\centerline{Fan Zhou}
\centerline{\href{mailto:fz2326@columbia.edu}{{\tt{\textcolor{black}{fz2326@columbia.edu}}}}}
\centerline{\small{February 9, 2024}}

\begin{abstract}
    In this paper we homologically construct a (functorial) BGG resolution of the finite-dimensional simple module of the nilBrauer algebra by using infinity-categorical methods following the reconstruction-from-stratification philosophy, e.g. appearing in Ayala-MazelGee-Rozenblyum. To do so, we prove a fact of independent interest, that half of the nilBrauer algebra is Koszul. This BGG resolution categorifies a character formula of Brundan-Wang-Webster. More generally, we have a (functorial) ``BGG spectral sequence'' which converges to any desired module; this spectral sequence is secretly a resolution when the desired module is finite-dimensional. This spectral sequence also categorifies the character formulae of Brundan-Wang-Webster for any (possibly infinite-dimensional) simple module. We expect the methods used here for producing BGG resolutions to be applicable to other (graded) triangular-based algebras also, especially diagrammatic ones.
\end{abstract}




\textcolor{white}{$\Db$}\vspace{-1em}

\tableofcontents

\section{Introduction}\label{section:introduction}
\subsection{Introduction and previous work}
In 2023, Brundan-Wang-Webster (\cite{brundan2023nilbrauer1},\cite{brundan2023nilbrauer2}) defined the nilBrauer algebra, which is locally unital but not unital, proved that it is a (graded) triangular-based algebra in the sense of \cite{brundan2023graded}, and showed that its representation theory categorified the split $\iota$-quantum group of rank 1, $U^{\iota}_q(\sl_2)$. One of the character formulas they obtained for nilBrauer modules was that (\cite[Equation 2.32]{brundan2023nilbrauer2})
\[[L_n]=\sum_{k=0}^\infty(-1)^k\frac{q^{-k(1+2\delta_{n\not\equiv t})}}{(1-q^{-4})(1-q^{-8})\cdots(1-q^{-4k})}[\ol\Delta_{n+2k}],\tag{BWW}\label{eqn:bww}\] 
where $L_n$ is the simple module of lowest weight $n$ and $\ol\Delta_{n+2k}$ is the ``proper standard module'', or ``small Verma module'', of lowest weight $n+2k$. They show that this character formula of modules categorifies a change-of-basis (dual PBW to dual canonical) formula for the $\iota$-quantum group, where $[L_n]$ corresponds to the dual $\iota$-canonical basis and $[\ol\Delta_n]$ corresponds to the dual PBW basis. A natural question to ask is then whether this character formula, reminiscent of the Weyl character formula for category $\CO$, can be further categorified as a `BGG resolution'.


The classical construction of the BGG resolution in category $\CO$ was done by extensively studying the internal structure of the Verma modules. In \cite{dhillon2019bernsteingelfandgelfand}, Dhillon takes a different approach, of a more categorical or homological nature. Rather than study the internal structure of Verma modules, which we will henceforth call `Vermas', Dhillon used information about the Ext groups between Vermas and simples to obtain BGG resolutions for highest weight categories, such as a block of category $\CO$, where this Ext information is provided in the form of Kostant's theorem on Lie algebra cohomology. Here in this paper we slightly generalize this approach to `stratified' settings, namely settings where the notion of Verma modules is separated between ``big Vermas'' and ``small Vermas''; in the literature, these modules are typically called ``standard modules'' and ``proper standard modules'', respectively. The reason we take this departure from the classical is that the distinction of big versus small Vermas makes the story more complicated, and in the nilBrauer case the maps between standard modules are not as simple as a scalar times a canonical embedding; moreover our BGG result shows that the terms are nontrivial extensions of proper standard modules, which also make the explicit construction of a map difficult.

Dhillon's approach fits philosophically within the reconstruction-from-stratification framework of \\ \cite{ayala2022stratified}, where the terminology of `noncommutative stacks' was used instead. In that paper, a framework is developed for a series a statements along the vein of `reconstruction', namely that when a (presentable) stable $\infty$-category $\CC$ is `stratified' in a way so that there are recollements at each level, the original category $\CC$ can be `reconstructed' from the `strata' (the quotient categories). Ayala-MazelGee-Rozenblyum gave several categorical levels of such statements; one such statement, made across Remarks 1.3.12 and 1.3.13 of \cite{ayala2022stratified}, says that one can filter the identity functor $\Id_\CC$ with the graded pieces coming from the information about the recollements. This is the statement we adopt here. Our conventions are slightly different from those of \cite{ayala2022stratified}, but we restate all relevant definitions here and provide a proof of the statement regarding the filtration of the identity functor. 

The matter of constructing BGG resolutions for representations of a diagrammatic algebra or category has been taken up before in \cite{khovanov4171760}. The ``Chebyshev algebra'', or perhaps ``doubly degenerate Temperley-Lieb algebra'', defined in that paper also fits the setting of triangular-based algebras from \cite{brundan2023graded}, and due to the semisimplicity of the ``Cartan algebras'', there is no distinction between big and small Vermas. Because the diagrammatic Chebyshev algebra is relatively simple, the authors of \cite{khovanov4171760} were able to construct morphisms between Vermas and directly show that a certain complex of Vermas was a resolution for a simple module $L_n$. They also computed the dimensions of various Ext groups and showed that the entire Chebyshev algebra is Koszul under the grading given by the total number of caps and cups. Let us remark that, by following the philosophy of the proof of Koszulity in the present paper, one can see that the ``upper half'' of the Chebyshev algebra is also Koszul; this, combined with the fact that every simple module over the Chebyshev algebra is one-dimensional, will recover the BGG result of \cite{khovanov4171760}. The dimensions of the multiplicity spaces in the resolution can also be computed from a dimension count of Ext groups. 

As we discuss below, our main results in this paper are that the lower half of the nilBrauer algebra, namely $\NB^-$, is Koszul, and that the simple representations of the nilBrauer algebra have BGG spectral sequences (which are resolutions when the simple is finite-dimensional). This BGG spectral sequence is a categorification of the character formula (\ref{eqn:bww}), and that it is a resolution when the simple is finite-dimensional is proved by using the Koszulness of $\NB^-$. The Koszulness of $\NB^-$ is also of independent interest, for the following reason. According to \cite{kapranovschechtman2020shuffle}, braided ((co)connected) bialgebras which are Koszul can be realized as factorizable sheaves which are moreover shifted-perverse on the space $\Sym^n\BC$. As we prove $\NB^-$ is Koszul, this provides a geometric construction of the lower half of the nilBrauer; cf. the geometric construction in \cite{kapranovschechtman2020shuffle} for the upper half of a quantum group, $U_q^+\glie$.

\subsection{Main Results}
In this paper, our first main result is that half of the nilBrauer algebra (the $A^-$ and $A^+$ alluded to below) is Koszul, and we show how this implies that the finite-dimensional trivial module is Kostant (and therefore admits a BGG resolution). In the process, we define an analogue of Lie algebra cohomology with respect to $\nlie^+$. 

Our second main result is the construction of a BGG resolution which categorifies the character formula (\ref{eqn:bww}) for $n=0$ and $t=0$ ($L_0$ for $t=0$ is the only simple module which is finite-dimensional), much like how the BGG resolution in category $\CO$ of the finite-dimensional simples categorifies the Weyl character formula. Here a BGG resolution means a resolution in which each term is an extension of proper standard modules $\ol\Delta_\lbd$; we also completely determine the structure of this extension, namely the way in which these proper standard modules fit together. (That they should be proper standard modules rather than standard modules is perhaps hinted at by the appearance of proper standard modules in the character formula (\ref{eqn:bww}).) More precisely, the $n$-th term of this BGG resolution $C^\blt_\te{BGG}(L_0)$ (which is homologically positioned so that it is equal to $L_0[0]$ in the derived category) has a decreasing filtration such that the $k$-th graded piece is
\[\gr^kC^{-n}_\te{BGG}(L_0)=\ol\Delta_{2n}\otimes_\BC q^{-n}\BC[\tsl p_2,\tsl p_4,\cdotsc, \tsl p_{2n}]_{\deg_\sym=k},\] 
where $\deg_\sym$ is the symmetric degree, defined\footnote{This is not to be confused with the quantum degree, which we do not decorate because it came first, in which $\deg \tsl p_i=2i$.} so that $\deg_\sym \tsl p_i=1$. 

As a generalization of our second main result, for all values of $n$ and $t$, we produce a (functorial\footnote{This functoriality is thanks to some $\infty$-category magicks.}) ``BGG spectral sequence'' converging to $L_n$ such that by taking the sum $$\sum_{p,q}(-1)^{p+q}[E_1^{p,q}]=\sum_{p,q}(-1)^{p+q}[E_\infty^{p,q}],$$ the spectral sequence categorifies the character formula above also. This spectral sequence deserves the name of BGG because each term on the first page $E_1^{p,q}$ is an (somewhat complicated) extension of proper standard modules $\ol\Delta_{\lbd}$, and the infinity page $E_\infty^{p,q}$ is $L_n$ sitting at $p=q=0$. Hence we provide a more complete categorification of the canonical-to-standard change-of-basis formula for the split $\iota$-quantum group of rank 1, $U_q^{\iota}(\sl_2)$. We however disclaim that we did not completely determine the structure of these extensions in this general case.

We expect the weight stratification machinery and the slogan relating BGG resolutions to the Koszulity of the half algebra to work in great generality for (graded) triangular-based algebras (\cite{brundan2023graded}), especially diagrammatic ones. The nilBrauer algebra is only locally graded-finite-dimensional, has infinite homological dimension, and fails to be Noetherian, so in some sense it is already very poorly-behaved; that our techniques work for this algebra seems to suggest it should work for other difficult algebras also. 

Our construction of the BGG resolution for $L_0$ is homologically powered by our determination of the Ext groups $\Ext(\Delta(\theta),L_0)$, where $\Delta(\theta)\coloneqq A^{\ge \theta}e^\theta=\bigoplus_{\lbd\in\theta}\Delta_\lbd^{\ol{l_\lbd(\theta)}}$ (in our case of $A=\NB$ there is only one term in this direct sum). We remark that the determination of other such Ext groups $\Ext(\Delta(\theta),L_\lbd)$ for general $\lbd\in\BN$ appears to be difficult; in particular, the problem to determine Kazhdan-Lusztig coefficients for the nilBrauer algebra appears to be difficult. The determination of such coefficients would shed more light on the structure of the extensions of proper standard modules in our BGG spectral sequence.

\subsection{Key Ideas and Outline}
There are two key ideas used in this paper. The first is that whenever there is a stratification of recollement situations of module categories, e.g. whenever there is a weight theory, one can put a `filtration' (suitably defined; technically speaking this requires some $\infty$-category language) on the identity functor, with the graded pieces consisting of homological information. This homological information, in the case of graded triangular bases, comes in the form of certain Ext groups (``standard Ext groups''), to be thought of as weight spaces of a cohomology. There is then a spectral sequence converging to the identity functor, and if one is lucky then this spectral sequence is actually a resolution when one plugs in `nice' (Kostant) simple modules to the identity functor. This deserves to be called `the BGG resolution', and indeed this gives the classical BGG resolution in category $\CO$. More generally, even if this spectral sequence is not a resolution, it is still a categorification of the appropriate character formula. The usage of the $\infty$-category language here makes this spectral sequence functorial. This idea is not original; this general machinery of stratifications of $\infty$-categories or derived categories was for example developed in \cite{ayala2022stratified}, where the language of `noncommutative stacks' is used, and \cite{dhillon2019bernsteingelfandgelfand}, where Dhillon uses similar ideas to recover the BGG resolution in (a block of) category $\CO$. In particular, we were taught this idea by Dennis Gaitsgory and Charles Fu.

The second key idea is, in the form of a slogan, ``Koszulity implies Kostantness implies BGG resolutions''. A longer slogan would be
\[\textbf{``Koszulity of half of $A$ (namely $A^-$) is crucial to finding BGG resolutions.''}\]
More precisely, due to the `diagonal Ext groups' formulation of Koszulness, showing that the `lower nilalgebra' or `half-algebra' $A^-$ (the analogue of $\nidl^+\subset \glie$ from Lie theory) is Koszul will show that the `trivial module' ($L_0=\bk e^0$ for $t=0$ in the case of nilBrauer) is Kostant, i.e. has concentrated standard Ext groups. This concentration ensures that the machinery from the first step actually gives a resolution. More generally, one would look for modules which are Koszul over the Koszul algebra $A^+$. 

The paper is organized as follows. In Section \ref{section:setting} we recall some set-up from \cite{brundan2023graded}, \cite{brundan2021semiinfinite}, and \cite{brundan2023nilbrauer1}. In Section \ref{section:reconstruction} we formulate and give proofs of some known facts regarding the reconstruction of the identity functor from stratification data. In Section \ref{section:nilcohomo} we form the `nilalgebra' $A^-$ (in accordance to ideas already present in the literature, e.g. \cite{brundan2021semiinfinite}) and define the ``nilcohomology'' $H^\blt(A^-:M)$ in analogy to Lie algebra cohomology $H^\blt(\nidl^-:M)$. In Section \ref{section:koszul} we show that the nilBrauer nilalgebra $\NB^-$ defined in Section \ref{section:nilcohomo} is Koszul. In Section \ref{section:bgg} we use the reconstruction machine set-up in Section \ref{section:reconstruction} and the Koszulness result from Section \ref{section:koszul} to show our main result, that the one-dimensional representation of $\NB$ is Kostant, and to construct the corresponding BGG resolution. 


\subsection{Future Work}
As mentioned above, we expect this machinery to work in general for many (graded) triangular-based diagrammatic algebras. In future work we plan to attack other such algebras such as KLR and the Lauda categorification of $U_q\sl_2$. Another `character formula' which we hope to some day categorify in this way is the Jacobi-Trudi determinant formula -- though this has already been categorified as a `BGG resolution' by many authors, they do not realize such resolutions as resolutions by Vermas in a highest weight/stratified category. 

The triangular-based formalism of \cite{brundan2023graded} was built to handle path algebras of categories, or rather ``Hom algebras of categories''. Another possible direction is to expand this machinery (and by extension the techniques of this paper) to work for ``Ext algebras of categories''. An example of this is the Ext algebra of Soergel bimodules.

\subsection{Acknowledgements} 
We would like to express our deepest gratitude to our advisors Mikhail Khovanov and Joshua Sussan for their guidance and invaluable help and support throughout this project, as well as Jon Brundan for first introducing us to this nilBrauer setting and the general theory of (graded) triangular bases, and also for the invaluable conversation we had at the WARTHOG workshop, as well as the emails that followed. 
We would also like to thank Ben Webster for a very helpful conversation about big Vermas at the MIT Summer School; the WARTHOG workshop for planting the seed of Koszulness in my mind; Catharina Stroppel for a very helpful conversation on Koszulness at an AIM workshop; Yixuan Li for his trick on determining semisimplicity using a Hom dimension count; Matthew Hase-Liu for his help in proving the combinatorial identity in Section \ref{subsec:kostant} (in particular the proof using finite differences is due to him), as well as his help in improving the proof of the general theorem on the graded pieces of the identity functor; Kevin Chang for conversations regarding \cite{kapranovschechtman2020shuffle} and geometric aspects of Koszul algebras; and Dennis Gaitsgory and Charles Fu, who first taught us a variation of the reconstruction-from-stratification philosophy. 
I would also like to thank my grandparents, who were with me when these results were obtained, for their support. 

\subsection{Conventions}
For most of this paper we use lowest-weight notation. We may sometimes simply say `algebra' when we mean to say locally unital algebra. Our convention on graded vector spaces is that $q^n$ shifts the grading \textit{down} by $n$, namely $(q^nV)_d=V_{d+n}$, and our definition of the graded dimension is
\[\dim_q V=\sum_n q^{-n}\dim_\bk V_n.\] 
We will frequently drop the $q$ subscript on $\dim$ and simply understand $\dim V$ to be a power series in $q$.

Unless otherwise stated, we will be working with left modules. 

Our notation for weight theory is that the weight labelling a module (e.g. the highest or lowest weight of a Verma) is a \textit{subscript}, e.g. $\Delta_\lbd\in\CO$ is the Verma with highest weight $\lbd$; and the weight space of a module is indicated by a superscript, e.g. $\Delta_\lbd{}^{\mu}$ is the $\mu$-weight space of $\Delta_\lbd$. 

In the case of a triangular-based algebra $A$ with weight set $\Theta$ and simples labelled by $\Lbd$ in bijection with the distinguished idempotents, we let $L_\lbd(\theta)$ and $P_\lbd(\theta)$ denote the simple and projective objects over the Cartan algebra, $A^\theta=e^\theta A^{\ge\theta} e^\theta$. 

For a unital bounded-below graded algebra $A$ with simples and their projective covers labelled by $\lbd\in\Lbd$, the projective cover of a module $M$ is given by 
\[P_M=\bigoplus_{\lbd\in\Lbd}P_\lbd^{\oplus\ol{\dim_q\Hom_A(M,L_\lbd)}},\]
where $\ol{\sq}$ refers to the bar-involution, sending $q\lmto q^{-1}$; in particular the algebra itself as a module over itself is $A=\bigoplus_{\lbd} P_\lbd^{\ol{\dim_q L_\lbd}}$. Dually one has that the injective hull of $M$ has $Q_M=\bigoplus_\lbd Q_\lbd^{\oplus\dim_q\Hom_A(L_\lbd,M)}$. The dual of a graded vector space is $V^\dag\coloneqq \bigoplus_{k\in\BZ}(1^kV_{-k})^*$.

We also denote by $\sq$ an arbitrary input; for example, the functor which is tensoring by $M$ would be denoted $\CF=M\otimes\sq$. 

In a spectral sequence, our convention is that $E_r^{p,q}$ has differentials of degree $(r,1-r)$. 

Our convention is that $0\in\BN$.

\section{The Setting}\label{section:setting}
In this paper we are concerned with locally unital graded algebras, which are infinite in two senses: that there may be infinitely many orthogonal idempotents (and in particular there is no unit, $1_A$), and that for two such idempotents $e_1,e_2$ the literal dimension $\dim_\bk e_1 A e_2$ may be infinite. Such orthogonal idempotents are also called local units. Such algebras, infinitely large as they are, are difficult to study without some type of finiteness handle, so we consider locally unital graded algebras which are ``locally graded-finite-dimensional'', which means that $\dim_\bk (e_1 A e_2)_n<\infty$ for every $n,e_1,e_2$, i.e. that $\dim_q e_1 Ae_2 \in\BN\pparen{q,q^{-1}}$ for all local units $e_1,e_2$. (Here `local' refers to the fact that we truncate $A$ by $e_1,e_2$, and `graded-finite-dimensional' means that each piece of this grading is finite-dimensional.) If moreover $\dim_\bk (e_1 A e_2)_n=0$ for $n\gg 0$, i.e. if $\dim_q e_1 A e_2\in\BN\pparen{q}$, then we say $A$ is ``locally bounded above''; and similarly if $\dim_\bk (e_1 A e_2)_n=0$ for $n\ll 0$, i.e. if $\dim_q e_1 A e_2\in\BN\pparen{q^{-1}}$, then we say $A$ is ``locally bounded below''. 

If one were so inclined, one could phrase this in a categorical language -- interpreting this locally unital algebra as the path algebra of a category, the properties that there may be infinitely many orthogonal idempotents becomes that there may be infinitely many objects (but no more than a set's worth) in this category, and the property that $\dim_q e_1 A e_2\in\BN\pparen{q,q^{-1}}$ becomes that this category is enriched over graded-finite-dimensional vector spaces. Similarly that this algebra is locally bounded above (resp. below) becomes that this category is enriched over graded-finite-dimensional vector spaces which are bounded above (resp. below).

\subsection{Triangular-based algebras}\label{subsec:triangular}
Let $A$ be a locally graded-finite-dimensional locally unital algebra as above. Let $I$ index the local units (orthogonal and homogeneous) of $A$, $\Phi\subseteq I$ label the ``distinguished idempotents'', $\Theta$ be a (lower-finite)\footnote{Perhaps this could be weakened to locally finite.} poset of ``weights'', and $\varpi\colon \Phi\lto\Theta$ be a map (with finite fibers) from distinguished idempotents to the weight poset. Diagrammatically this is:
\begin{center}
    \begin{tikzcd}
I                                                           &        \\
\Phi \arrow[u, hook] \arrow[r, "\varpi", two heads] & \Theta
\end{tikzcd}
\end{center}
We will let symbols like $i,j\in I$ and $\alpha,\beta\in\Phi$ and $\theta,\phi,\psi\in\Theta$. In the interest of appealing to weight notation later, we will denote the local unit labelled by $i$ as $1^i$. 

Let $A$ have a ``(graded) triangular basis'' in the sense of Brundan, which is to say that 
\begin{DEF}
    $A$ is ``graded triangular based'' if there are (homogeneous) sets $\te X(i,\alpha)\subseteq 1^i A 1^\alpha,\ \te H(\alpha,\beta)\subseteq 1^\alpha A 1^\beta,\ \te Y(\beta,j)\subseteq 1^\beta A 1^j$ such that 
\begin{enumerate}
    \item products of these things in these sets give a basis for $A$, i.e. $$\set*{xhy:(x,h,y)\in\bigcup_{i,j,\alpha,\beta}\te X(i,\alpha)\times\te H(\alpha,\beta)\times \te Y(\beta,j)}$$ forms a basis of $A$;
    \item $\te X(\alpha,\alpha)=\tY(\alpha,\alpha)=\{1^\alpha\}$;
    \item for $\alpha\neq\beta$, 
    \begin{align*}
        \tX(\alpha,\beta)\neq\emptyset&\implies \varpi(\alpha)>\varpi(\beta),\\
        \tH(\alpha,\beta)\neq\emptyset&\implies \varpi(\alpha)=\varpi(\beta),\\
        \tY(\alpha,\beta)\neq\emptyset&\implies \varpi(\alpha)<\varpi(\beta);
    \end{align*}
    \item for each $i\in I-\Phi$, there are only finitely many $\alpha\in\Phi$ such that $\tX(i,\alpha)\cup \tY(\alpha,i)\neq\emptyset$.
\end{enumerate}
\end{DEF}

Let 
\[e^\theta\coloneqq\sum_{\alpha\in\varpi^{-1}\theta}1^\alpha.\] 
Let $A^{\ge\theta}$ be defined as
\[A^{\ge\theta}\coloneqq \mfrac{A}{\wan{e^\phi\colon\phi\not\ge\theta}},\]
let e.g. $A^{\le\theta},A^{>\theta}$ be defined similarly, and let 
\[A^\theta=e^\theta A^{\ge\theta}e^\theta.\] 
This plays the role of Cartan in the sense that modules over it are induced to form standard objects. Note that
\[A^{>\theta}=\mfrac{A^{\ge\theta}}{A^{\ge\theta}e^\theta A^{\ge\theta}};\] 
hence, by the general theory of Cline-Parshall-Scott \cite{Cline1988}, one expects a recollement $\D\,\Mod\,A^{>\theta}\lto\D\,\Mod\,A^{\ge\theta}\lto \D\,\Mod\,A^\theta$. For $M$ a module over $A$, we write $M^\theta\coloneqq e^\theta M$ in analogy with weight space notation. We often write simply $\alpha\in\theta$ as shorthand for $\alpha\in\varpi^{-1}(\theta)$, so that $e^\theta=\sum_{\alpha\in\theta} 1^\alpha$. 

Let $\Lbd_\theta$ label the simple modules $L_\lbd(\theta)$ of $A^\theta$ as well as their projective covers $P_\lbd(\theta)$ and injective hulls $Q_\lbd(\theta)$, and let $\Lbd=\bigsqcup_{\theta}\Lbd_\theta$ be the disjoint union of these label sets. Then one can show that $\Lbd$ labels all simples of $A$. Frequently it is the case (for example if the semisimplification of the Cartan is $A^\theta/\jidl\cong\prod_{\alpha\in\varpi^{-1}(\theta)}\bk$) that $\Lbd$ can be identified with $\Phi$; in this case we will instead use the symbols $\lbd,\mu\in\Lbd$ for the distinguished idempotents. Here is a picture of the situation:
\begin{center}
    \begin{tikzcd}
I                                                                                       & \Lambda \arrow[d, two heads] \\
\Phi \arrow[u, hook] \arrow[r, "\varpi", two heads] \arrow[ru, no head, dashed] & \Theta                      
\end{tikzcd}
\end{center}

There is moreover often a ``split triangular decomposition'', which is to say that 
\begin{DEF}
$A$ has a ``split triangular decomposition'' if there are $A^-,A^0,A^+$ locally unital algebras such that
\begin{enumerate}
    \item $A^\flat=A^-A^0$ and $A^\sharp=A^0A^+$ are subalgebras;
    \item the multiplication map $A^-\otimes_\BK A^0\otimes_\BK A^+\lto A$ is a linear isomorphism;
    \item $1^i A^\flat 1^i=1^i A^\sharp 1^i=\bk 1^i$ for each $i$ and
    \begin{align*}
        1^i A^0 1^j\neq\emptyset&\implies \varpi(i)=\varpi(j),\\
        1^j A^- 1^i,1^iA^+1^j\neq \bk 1^i &\implies \varpi(i)>\varpi(j).
    \end{align*}
\end{enumerate}
\end{DEF}
We will see later that nilBrauer falls into this set-up.

\subsection{Recollement}
Recall the following recollement diagrams from \cite{brundan2023graded}:

\begin{center}
\begin{tikzcd}
&  & \perp &  &  &  & \perp &  &  \\[-1em]
\D^-\Mod\, A^{>\theta} \arrow[rrrr, "\iota_\theta"] &  &       &  & \D^-\Mod\, A^{\ge\theta} \arrow[rrrr, "\jota^\theta=e^\theta \sq "] \arrow[llll, swap,"\iota_\theta^*=A^{>\theta}\lotimes_{A^{\ge\theta}} \sq", bend right=60] \arrow[bend left=60]{llll}{\iota_\theta^!=\bigoplus_i\rhom_{A^{\ge\theta}}(A^{>\theta}1^i, \sq)} &  &       &  & \D^-\Mod\, A^\theta \arrow[llll, swap,"\jota^\theta_!=A^{\ge\theta}e^\theta\otimes_{A^\theta} \sq", bend right=60] \arrow[bend left=60]{llll}{\jota^\theta_*=\bigoplus_i\Hom_{A^\theta}(e^\theta A^{\ge\theta}1^i, \sq)} \\[-1em]
  &  & \perp &  & & &\perp &  &
\end{tikzcd}
\end{center}
and 
\begin{center}
\begin{tikzcd}
&  & \perp &  & \\[-1em]
\D^-\Mod\, A^{\ge\theta} \arrow[rrrr, "\ol\iota_\theta"] &  &       &  & \D^-\Mod\,A \arrow[bend left=60]{llll}{\ol\iota_\theta^!=\bigoplus_i\rhom_A(A^{\ge\theta}1^i,\sq)}\arrow[llll, "\ol\iota_\theta^*=A^{\ge\theta}\lotimes_A\sq", bend right=60,swap] \\[-1em]
&  & \perp &  & 
\end{tikzcd}
\end{center}
As noted earlier, this is an artifact of the Cline-Parshall-Scott theory of algebraic recollements. Define the ``(co)standard modules'' and ``proper (co)standard modules'' $\Delta$ and $\ol\Delta$ (resp. $\Nabla$ and $\ol\Nabla$) by:
\[\Delta_\lbd=\jota^\theta_! P_\lbd(\theta),\qquad \ol\Delta_\lbd=\jota^\theta_! L_\lbd(\theta),\] 
\[\Nabla_\lbd=\jota^\theta_* Q_\lbd(\theta),\qquad\ol\Nabla_\lbd=\jota^\theta_*L_\lbd(\theta).\]

Let us rewrite a couple of these functors in another slightly more familiar form. 
\begin{align*}
    \jota^\theta&=\hom_A(A^{\ge\theta}e^\theta,\sq)
\end{align*}
is evident, as such maps are determined by where $e^\theta$ is sent (note well the domain of this functor). Now let there be a duality coming from an antiautomorphism of $A$. On the level of $A$-modules this is denoted $\sq^\dag$, while on the level of $A^\theta$-modules we will denote this as $\sq^*$. Recall from \cite[Section 5]{brundan2023graded} that duality exchanges the shriek and the star, so that one has
\[\ol\iota_\theta^* M=(\ol\iota_\theta^! M^\dag)^\dag=\pr*{\bigoplus_i \hom_A(A^{\ge\theta}1^i,M^\dag)}^\dag,\]
where the left $A^{\ge\theta}$-action on the Hom comes from the right action on $\bigoplus_i A^{\ge\theta}1^i$. Then 
\begin{align*}
    \jota^\theta \ol\iota_\theta^*M&=e^\theta\pr*{\bigoplus_i \hom_A(A^{\ge\theta}1^i,M^\dag)}^\dag\\
    &=\pr*{\bigoplus_i e^\theta\hom_A(A^{\ge\theta}1^i,M^\dag)}^*\\
    &=\pr*{\bigoplus_i \hom_A(A^{\ge\theta}1^ie^\theta,M^\dag)}^*\\
    &=\hom_A\pr*{\bigoplus_{\alpha\in\theta} A^{\ge\theta}1^\alpha,M^\dag}^*\\
    &=\hom_A(A^{\ge\theta}e^\theta,M^\dag)^*
\end{align*}
Note well that there is a left $A^\theta$-action on this coming the from right $e^\theta A^{\ge\theta}e^\theta$-action on $A^{\ge\theta}e^\theta$; it is with respect to this action that the outermost dual is taken.

Now suppose there is an identification of $\Phi$ with $\Lbd$; this is for example the case with nilBrauer, where in fact $I=\Phi=\Lbd=\Theta=\BN$. The decomposition $\sum_{\lbd\in\theta} 1^\lbd$ of the unit $e^\theta$ of $A^\theta$ gives a decomposition into projectives
\[A^\theta e^\theta=\bigoplus_{\lbd\in\theta} P_\lbd(\theta)^{\oplus \ol{\dim L_\lbd(\theta)}},\] 
where the bar is included for the graded case. For shorthand we will let $l_\lbd(\theta)=\dim L_\lbd(\theta)$. Then
\[A^{\ge\theta}e^\theta= \bigoplus_{\lbd\in\theta}\Delta_\lbd^{\ol{l_\lbd(\theta)}}.\]
As shorthand, we shall let $$\Delta(\theta)\coloneqq \bigoplus_{\lbd\in\theta}\Delta_\lbd^{\ol{l_\lbd(\theta)}}.$$

In conclusion, from the above discussion, one has (recall $\jota^\theta_!$ is exact in the triangular-based setting)
\begin{LEM}
\[\jota^\theta_!=\bigoplus_{\lbd\in\theta} \Delta_\lbd^{\ol{l_\lbd(\theta)}}\otimes_{A^\theta}\sq=\Delta(\theta)\otimes_{A^\theta}\sq\]
as well as
\[\jota^\theta\ol\iota_\theta^*=\rhom_A\pr*{\bigoplus_{\lbd\in\theta}\Delta_\lbd^{\ol{l_\lbd(\theta)}},\sq^\dag}^*=\rhom_A(\Delta(\theta),\sq^\dag)^*.\]
\end{LEM}
Again note well where the $A^\theta$-action on $\jota^\theta\ol\iota_\theta^*$ comes from -- it comes from all the standard modules bunched together, and if one took RHom from each individual standard module then one loses this action.

\subsection{The nilBrauer algebra}
Let us briefly recall some theory from \cite{brundan2023nilbrauer1} and \cite{brundan2023nilbrauer2}. The nilBrauer category, which we denote as $\NB=\NB_t$ in an abuse of notation, is a strict graded monoidal category depending on a parameter $t$ with one generating object $B$, whose identity endomorphism is represented diagrammatically by an unlabelled vertical string, also called a propagating string/strand:
\begin{center}
    $\Id_B=$\hackcenter{\begin{tikzpicture}[scale=0.375]
        \draw (0,0)--(0,2);
    \end{tikzpicture}}\ .
\end{center}
The morphisms in this category are generated by the following four morphisms,
\begin{align*}
    {}\hackcenter{\begin{tikzpicture}[scale=0.375]
        \draw (0,0)--(0,2);
        \fill (0,1) circle (5pt);
    \end{tikzpicture}}\ ,&\te{ of degree }2,\\
    {}\hackcenter{\begin{tikzpicture}[scale=0.375]
        \draw (0,0)--(2,2);
        \draw (2,0)--(0,2);
    \end{tikzpicture}}\ ,&\te{ of degree }-2,\\
    {}\hackcenter{\begin{tikzpicture}[scale=0.375]
        \draw (0,0) arc (180:0:1);
    \end{tikzpicture}}\ ,&\te{ of degree }0,\\
    {}\hackcenter{\begin{tikzpicture}[scale=0.375]
        \draw (0,0) arc (180:360:1);
    \end{tikzpicture}}\ ,&\te{ of degree }0,
\end{align*}
subject to the following relations, where $t=0$ or $t=1$:
\begin{align*}
    {}\hackcenter{\begin{tikzpicture}[scale=0.375]
        \draw (0,0)--(2,2)--(0,4);
        \draw (2,0)--(0,2)--(2,4);
    \end{tikzpicture}}&=0, &{}\hackcenter{\begin{tikzpicture}[scale=0.375]
        \draw (0,0)--(4,4);
        \draw (0,4)--(4,0);
        \draw (2,0)--(0,2)--(2,4);
    \end{tikzpicture}}&=\hackcenter{\begin{tikzpicture}[scale=0.375]
        \draw (0,0)--(4,4);
        \draw (0,4)--(4,0);
        \draw (2,0)--(4,2)--(2,4);
    \end{tikzpicture}},\\
    {}\hackcenter{\begin{tikzpicture}[scale=0.375]
        \draw (0,0) arc (0:360:1);
    \end{tikzpicture}}&=t \Id_{\mathbf{1}},
    &
    {}\hackcenter{\begin{tikzpicture}[scale=0.375]
        \draw (0,0)--(0,3);
        \draw (0,3) arc(180:0:1);
        \draw (2,3)--(2,1);
        \draw (2,1) arc(180:360:1);
        \draw (4,1)--(4,4);
    \end{tikzpicture}}&=\hackcenter{\begin{tikzpicture}[scale=0.375]
        \draw (0,0)--(0,4);
    \end{tikzpicture}}=\hackcenter{\begin{tikzpicture}[scale=0.375]
        \draw (0,4)--(0,1);
        \draw (0,1) arc(180:360:1);
        \draw (2,1)--(2,3);
        \draw (2,3) arc(180:0:1);
        \draw (4,3)--(4,0);
    \end{tikzpicture}},\\
    {}\hackcenter{\begin{tikzpicture}[scale=0.375]
        \draw (0,0)--(2,2);
        \draw (0,2) arc(180:0:1);
        \draw (2,0)--(0,2);
    \end{tikzpicture}}&=0,
    &
    {}\hackcenter{\begin{tikzpicture}[scale=0.375]
        \draw (0,0) arc(180:0:1);
        \draw (1,0)--(0,2);
    \end{tikzpicture}}&={}\hackcenter{\begin{tikzpicture}[scale=0.375]
        \draw (0,0) arc(180:0:1);
        \draw (1,0)--(2,2);
    \end{tikzpicture}},\\
    {}\hackcenter{\begin{tikzpicture}[scale=0.375]
    \draw (0,0)--(2,2);
    \draw (2,0)--(0,2);
    \fill (0.5,1.5) circle (5pt);
\end{tikzpicture}}
-
\hackcenter{\begin{tikzpicture}[scale=0.375]
    \draw (0,0)--(2,2);
    \draw (2,0)--(0,2);
    \fill (1.5,0.5) circle (5pt);
\end{tikzpicture}}
&=
\hackcenter{\begin{tikzpicture}[scale=0.375]
    \draw (0,0) to[out=60,in=-60] (0,2);
    \draw (2,0) to[out=120,in=-120] (2,2);
\end{tikzpicture}}
-
\hackcenter{\begin{tikzpicture}[scale=0.375]
    \draw (0,0) to[out=60,in=120] (2,0);
    \draw (0,2) to[out=-60,in=-120] (2,2);
\end{tikzpicture}},
    &
    {}\hackcenter{\begin{tikzpicture}[scale=0.375]
        \draw (0,0) arc(180:0:1);
        \fill (1,0) ++(135:1) circle (5pt);
    \end{tikzpicture}}&=-{}\hackcenter{\begin{tikzpicture}[scale=0.375]
        \draw (0,0) arc(180:0:1);
        \fill (1,0) ++(45:1) circle (5pt);
    \end{tikzpicture}}.
\end{align*}
From these defining relations, one can then show that the following relations are also satisfied:
\begin{align*}
    {}\hackcenter{\begin{tikzpicture}[scale=0.375]
        \draw (0,0) arc(180:360:1);
        \draw (1,0)--(0,-2);
    \end{tikzpicture}}&={}\hackcenter{\begin{tikzpicture}[scale=0.375]
        \draw (0,0) arc(180:360:1);
        \draw (1,0)--(2,-2);
    \end{tikzpicture}},
    &
    {}\hackcenter{\begin{tikzpicture}[scale=0.375]
        \draw (0,0)--(0,1)--(2,3);
        \draw (2,3) arc(180:0:1);
        \draw (4,3)--(4,1);
        \draw (4,1) arc(0:-180:1);
        \draw (2,1)--(0,3)--(0,4);
    \end{tikzpicture}}&=0={}\hackcenter{\begin{tikzpicture}[scale=0.375]
        \draw (4,0)--(4,1)--(2,3);
        \draw (2,3) arc(0:180:1);
        \draw (0,3)--(0,1);
        \draw (0,1) arc(180:360:1);
        \draw (2,1)--(4,3)--(4,4);
    \end{tikzpicture}},\\
    {}\hackcenter{\begin{tikzpicture}[scale=0.375]
        \draw (0,0)--(2,2);
        \draw (2,0)--(0,2);
        \draw (0,0) arc(180:360:1);
    \end{tikzpicture}}&=0,
    &
    {}\hackcenter{\begin{tikzpicture}[scale=0.375]
        \draw (0,0) arc(180:0:2);
        \draw (0,3) arc(180:360:2);
    \end{tikzpicture}}&=0,\\
    \hackcenter{\begin{tikzpicture}[scale=0.375]
    \draw (0,0)--(2,2);
    \draw (2,0)--(0,2);
    \fill (0.5,0.5) circle (5pt);
\end{tikzpicture}}
-
\hackcenter{\begin{tikzpicture}[scale=0.375]
    \draw (0,0)--(2,2);
    \draw (2,0)--(0,2);
    \fill (1.5,1.5) circle (5pt);
\end{tikzpicture}}
&=
\hackcenter{\begin{tikzpicture}[scale=0.375]
    \draw (0,0) to[out=60,in=-60] (0,2);
    \draw (2,0) to[out=120,in=-120] (2,2);
\end{tikzpicture}}
-
\hackcenter{\begin{tikzpicture}[scale=0.375]
    \draw (0,0) to[out=60,in=120] (2,0);
    \draw (0,2) to[out=-60,in=-120] (2,2);
\end{tikzpicture}},
    &
    {}\hackcenter{\begin{tikzpicture}[scale=0.375]
        \draw (0,0) arc(180:360:1);
        \fill (1,0) ++(225:1) circle (5pt);
    \end{tikzpicture}}&=-{}\hackcenter{\begin{tikzpicture}[scale=0.375]
        \draw (0,0) arc(180:360:1);
        \fill (1,0) ++(-45:1) circle (5pt);
    \end{tikzpicture}}.
\end{align*}

We can consider the path algebra of this category, treating objects as idempotents and morphisms as elements of the algebra; the algebra so obtained is also denoted $\NB$, the ``nilBrauer algebra''. That this nilBrauer algebra fits into the triangular-based setting recalled in Section \ref{subsec:triangular} is the main result of \cite{brundan2023nilbrauer1}. Let us briefly describe how the nilBrauer fits that setting.

For the nilBrauer, one has $I=\Phi=\Theta=\Lbd=\BN$; in particular the weight idempotents $e^\theta$ correspond to the number of strands at the top or bottom of a diagram. Then the sets $\tX,\ \tH,\ \tY$ are defined by letting $\tX$ be the set of diagrams with no bubbles and no caps and no crossing propagating strands, where dots can only appear at the left end of cups; $\tH$ be the set of diagrams with neither cups nor caps, where propagating strands can cross but two strands can cross at most once, dots can only appear either at the bottom-left of propagating strands (away from crossings) or in bubbles, and bubbles can only appear to the right of the diagram; and $\tY$ be the set of diagrams with no bubbles and no cups and no crossing propagating strands, where dots can only appear at the left end of caps. 

Note that in $\tX$ (resp. $\tY$), cups (resp. caps) are allowed to intersect each other, though necessarily at most once due to the nilBrauer relations. For diagrams in any of the sets defined above, the number of strands at the top and at the bottom of the diagram must have the same parity. We also remark that the choice of putting dots (on cups/caps/propagating strands/bubbles) towards the left and bubbles towards the right is arbitrary; one can make any choice of distinguished locations to place the dots. 

Having seen that the nilBrauer fits the triangular-based formalism, one can then study its representation theory and form objects such as (big and small) Vermas, projectives, and simples. That the Grothendieck group of the nilBrauer algebra is isomorphic to the split $\iota$-quantum group of rank 1, $U^\iota_q(\sl_2)$, and that the projectives (resp. standard modules) go to the $\iota$-canonical basis (resp. PBW basis) under this isomorphism, is one of the main results of \cite{brundan2023nilbrauer2}, as recalled in Section \ref{section:introduction}. It is this categorification result which we further investigate in this paper, by categorifying the character formula in Equation (\ref{eqn:bww}). To do this, let us first establish some tools from $\infty$-category theory.

\section{Reconstruction from Stratification}\label{section:reconstruction}
\subsection{Rough picture}
This section rephrases ideas already present in the literature, for instance in \cite{ayala2022stratified}. The authors of that paper likely primarily had applications to algebraic geometry in mind; we will here use it to do representation theory, in much the same way as \cite{dhillon2019bernsteingelfandgelfand}. However, the setup/notation in that paper is slightly different from what we need, so we depart from their conventions in favor of that taught to us by Dennis Gaitsgory and Charles Fu. 

Roughly speaking, the general fact is that in appropriate recollement situations, namely where there is a big category stratified in such a way that each piece of the stratification is part of a recollement, the identity functor on the big category ($\Db\Mod\, A$ in this case) admits a filtration (appropriately defined) whose graded pieces can be described in terms of the recollement functors:
\begin{align*}
    \gr^\theta\calid&=\ol\iota_\theta\jota^\theta_!\jota^\theta\ol\iota_\theta^*.
    \intertext{In our case of locally unital graded algebras, we then know from the previous discussion that this is the same as}
    &=\Delta(\theta)\otimes_{A^\theta}\rhom_A\pr{\Delta(\theta),\sq^\dag}^*.
\end{align*}

If we can put an additional linear ordering on $\Theta$, for example via some sort of a length function (as is the case for usual category $\CO$), then there is a spectral sequence whose first page consists of (direct sums of) terms of the form above, and which converges to (the graded pieces of) the cohomology of the input. In the case of category $\CO$ one inputs the finite-dimensional simples to recover the BGG resolution; this works thanks to the fine control on the Ext groups $\Ext^\blt(\Delta,L)$ afforded by Kostant (in the guise of $\nidl^+$-cohomology, or `nilcohomology'). Note how the above generalizes the spectral sequence of standard objects appearing in Dhillon's work (the setting there is in an actual highest weight category). 

In the nilBrauer setting we do in fact have a linear ordering on $\Theta$ (per block, i.e. within even/odd), so one can show (and we will prove this later in this section)
\begin{THM}\label{thm:nbspecseq}
    For any $M\in\D^-(\Mod^\vtheta\,\NB)$ (here $\vtheta=0$ or $1$, namely even or odd, denotes the homological block $M$ lies in), one has
    \[\gr^k\calid_\vtheta M=\Delta(\vtheta+2k)\otimes_{\NB^{\vtheta+2k}} \rhom_{\NB}^\blt\big(\Delta(\vtheta+2k),M^\dag\big)^*;\] 
    and there is a spectral sequence
    \[E_1^{p,q}=H^{p+q}\gr^{-p}\calid_\vtheta M\specseqimplies E_\infty^{p,q}=\gr^{-p}H^{p+q} M\]
    which in this case reads as
    \[E_1^{p,q}=\Delta(\vtheta-2p)\otimes_{\NB^{\vtheta-2p}}\Ext^{-(p+q)}_\NB\big(\Delta(\vtheta-2p),M^\dag\big)^*\specseqimplies E_\infty^{p,q}=\gr^{-p}H^{p+q}M.\]
\end{THM}
So in order to obtain a BGG spectral sequence or resolution an appropriate simple module $L$, it remains to compute these Ext groups for $M=L$ as well as the $A^\theta$-action on them. In order to obtain a resolution, one can only hope that the Ext groups involved are homologically appropriately concentrated.

\subsection{With more pixels}

We use highest weight notation in this section on general machinery, and will return to lowest weight notation when we return to the setting of nilBrauer or other specific diagram algebras. Let $\sic$ denote the $\infty$-category of stable $\infty$-categories, and let $\Lbd$ be a poset with a unique final object (our convention on arrows is that final means maximal). 
\begin{DEF}
    A ``filtration'' on $\CC\in\sic$ is a functor
    \[\CF\colon\Lbd\lto\sic\]
    such that all arrows go to fully faithful embeddings and the final object goes to $\CC$. For $\lbd\in\Lbd$, we denote $\CC^{\le\lbd}=\CF(\lbd)$ and $\CC^{<\lbd}$ as the smallest stable subcategory containing all $\CC^{\le\mu}$ for all $\mu<\lbd$. Let $\CC^{=\lbd}=\CC^{\le\lbd}/\CC^{<\lbd}$ be the Verdier quotient.
\end{DEF}
Our convention is to write $\mu\to\lbd$ for $\mu<\lbd$. Let us name
\begin{align*}
    \ol\iota_\lbd&\colon \CC^{\le\lbd}\lto\CC,\\
    \iota_\mu^{<\lbd}&\colon \CC^{\le\mu}\lto\CC^{<\lbd},\\
    \iota_\mu^\lbd&\colon \CC^{\le\mu}\lto\CC^{\le\lbd},\\
    \iota_\lbd&\colon \CC^{<\lbd}\lto\CC^{\le\lbd},\\
\jota^\lbd&\colon \CC^{\le\lbd}\lto\CC^{=\lbd};
\end{align*}
for our purposes let us require both the inclusions of each $\CC^{\le\lbd}$ as well as the Verdier quotient functors to have both adjoints. This is for example the case with recollement situations. Let the adjoints, like in the recollement case, be named so that $\ol\iota_\lbd^*$ and $\jota^\lbd_!$ are left adjoints. 

In another direction, one can filter an object of an $\infty$-category by a poset $\Lbd$:
\begin{DEF}
    A ``$\Lbd$-filtered object of $\CC$'' is a functor
    \[\CX\in\sf{Fun}(\Lbd,\CC).\] 
    Letting $\Lbd^0$ be the 0-skeleton of $\Lbd$, the ``associated graded'' of this filtered object is 
    \begin{align*}
        \gr\colon\sf{Fun}(\Lbd,\CC)&\lto\sf{Fun}(\Lbd^0,\CC)\\
        \CX&\lmto \gr\CX
    \end{align*}
    defined by
    \[(\gr\CX)(\lbd)=\gr^\lbd\CX=\Fib\pr*{\CX(\lbd)\lto\lim_{\mu\from\lbd}\CX(\mu)},\]
    assuming such limits exist. 

    An object $X\in\CC$ is said to have a $\Lbd$-filtration if after adding an initial/minimal element $\lbd_{-\infty}$ to $\Lbd$ there is a functor $\CX\colon\Lbd\cup\{\lbd_{-\infty}\}\lto\CC$ such that $\CX(\lbd_{-\infty})=X$.
\end{DEF}
\begin{remark}
    Note that the definition of a filtered object here is sort of opposite to one might expect; rather than inclusions, one simply has maps $\CX(\mu)\lto\CX(\lbd)$ for each $\mu\to\lbd$ (i.e. $\mu\le\lbd$), where the smallest element $\lbd_{-\infty}$ corresponds to the actual object $X$, so that $X$ has maps to all the filtered pieces rather than the filtered pieces `including' into $X$; and rather than being a quotient, the associated graded pieces are sort of kernels. Of course in the derived category these things are the same up to shift anyway, but roughly speaking this is why there is a minus sign on the gr in Theorem \ref{thm:nbspecseq}.
\end{remark}

In the previous subsection we alluded to a general fact that said that in appropriate stratification/recollement situations, the identity functor on the big category admits a filtration. More precisely, the statement is
\begin{THM}[\cite{ayala2022stratified}, folklore]\label{thm:folklorefiltration}
    Let $\CC$ admit a filtration by $\Lbd$ such that the inclusions $\iota_\lbd$ and quotients $\jota^\lbd$ have both adjoints. Then the identity functor $\Id_\CC\in\on{\sf{End}}\CC$ admits a $\Lbd^\op$-filtration with terms
    \[\calid(\lbd)=\ol\iota_\lbd\ol\iota_\lbd^*,\] 
    and the associated graded can be computed as
    \[\gr^\lbd\calid=\ol\iota_\lbd\jota^\lbd_!\jota^\lbd\ol\iota_\lbd^*,\] 
    provided that $\Lbd$ is locally-finite and ``eventually disjoint-totally-ordered'' (which means that for any $\lbd$, there exists $\mu\le\lbd$ such that $\{\nu:\nu\le\mu\}$ is totally ordered).

    If there is a ``dimension/length function'' $\ell\colon\Lbd\lto\BZ^\op$ and $\CC$ moreover has a t-structure, then there is a spectral sequence
    \[E_1^{p,q}=\bigoplus_{\lbd\in\ell^{-1}(-p)}\pi^{p+q}(\gr^\lbd\calid)\specseqimplies E_\infty^{p,q}=\gr^{-p}\pi^{p+q}(\Id_\CC),\]
    where $\pi^\blt$ are the homotopy groups associated to the t-structure. 
\end{THM}
\vspace{-0.5em}
\begin{adjustwidth}{2em}{0pt}
    \begin{proof}
    In the definition above we see that to endow $\Id_\CC$ with a filtration is to define a functor $\calid\colon \Lbd^\op\cup\{\lbd_{-\infty}\}\lto \sf{Fun}(\CC,\CC)$; as $\Lbd$ has a unique maximal element $\lbd_\infty$ by assumption, it is not necessary to union this $\{\lbd_{-\infty}\}$. We will use the same symbol $\lbd_\infty$ to denote the unique minimal element of $\Lbd^\op$. We may well define $\calid\colon \Lbd^\op\lto\fun(\CC,\CC)$ by $\calid(\lbd)=\ol\iota_\lbd \ol\iota_\lbd^*$ and $\calid(\lbd_{\infty})=\Id$. Given $\lbd\overset\op\to\mu$ in $\Lbd^\op$, i.e. $\mu\to\lbd$ in $\Lbd$ (i.e. $\mu<\lbd$), we should also construct an arrow $\ol\iota_\lbd \ol\iota_\lbd^*\to \ol\iota_\mu \ol\iota_\mu^*$. This is given by the map
    \begin{align*}
        \id_{\ol\iota_\mu^* X}\in \hom_{\CC^{\le\mu}}(\ol\iota_\mu^* X,\ol\iota_\mu^* X)&\cong \hom_\CC(X,\ol\iota_\mu \ol\iota_\mu^* X)\\
        &\cong\hom_\CC(X,\ol\iota_\lbd \iota_\mu^\lbd\ol\iota_\mu^* X)\\
        &\cong \hom_{\CC^{\le\lbd}}(\ol\iota_\lbd^* X,\iota_\mu^ \lbd\ol\iota_\mu^* X)\\
        &\cong \hom_\CC(\ol\iota_\lbd \ol\iota_\lbd^* X, \ol\iota_\mu \ol\iota_\mu^* X).
    \end{align*}
    It is straightforward to check that the map so defined is a natural transformation of functors. Similarly, we define the arrow out of the initial object $\lbd_{\infty}\overset{\te{op}}\to \lbd$ to give rise to a morphism $\Id\to \ol\iota_\lbd \ol\iota_\lbd^*$ given by the unit of the adjunction. 

    One has the usual triangle from recollement
    \[\jota^\lbd_!\jota^\lbd\lto\Id_{\CC^{\le\lbd}}\lto\iota_\lbd\iota_\lbd^*\oslto{+1},\] 
    which one sandwiches to get
    \[\ol\iota_\lbd\jota^\lbd_!\jota^\lbd\ol\iota_\lbd^*\lto \ol\iota_\lbd\ol\iota_\lbd^*\lto \ol\iota_\lbd\iota_\lbd\iota_\lbd^*\ol\iota_\lbd^*\oslto{+1};\] 
    hence it suffices to show that
    \[\ol\iota_\lbd\iota_\lbd\iota_\lbd^*\ol\iota_\lbd^*=\lim_{\mu\to\lbd}\ol\iota_\mu\ol\iota_\mu^*.\] 

    It is maybe worth saying that in the diagram with respect to which the limit is taken above, the morphisms are those between the $\ol\iota_\mu\ol\iota_\mu^*$ constructed earlier. There is a map 
    \[\ol\iota_\lbd\iota_\lbd\iota_\lbd^*\ol\iota_\lbd^*\lto\ol\iota_\mu\ol\iota_\mu^*\]
    constructed as follows: take the unit 
    \[\Id\lto\iota_\mu^{<\lbd}\iota_\mu^{<\lbd}{}^*,\] 
    sandwich it with $\iota_\lbd$ and $\iota_\lbd^*$ on the left and right to obtain
    \[\iota_\lbd\iota_\lbd^*\lto \iota_\lbd\iota_\mu^{<\lbd}\iota_\mu^{<\lbd}{}^*\iota_\lbd^*=\iota_\mu^\lbd\iota_\mu^\lbd{}^*,\] 
    and now sandwich with $\ol\iota_\lbd$ and $\ol\iota_\lbd^*$ on the left and right to obtain the desired map
    \[\ol\iota_\lbd\iota_\lbd\iota_\lbd^*\ol\iota_\lbd^*\lto \ol\iota_\lbd\iota_\mu^\lbd\iota_\mu^\lbd{}^*\ol\iota_\lbd^*=\ol\iota_\mu\ol\iota_\mu^*.\]
    This then gives a map
    \[\ol\iota_\lbd\iota_\lbd\iota_\lbd^*\ol\iota_\lbd^*\lto\lim_{\mu\to\lbd}\ol\iota_\mu\ol\iota_\mu^*.\]

    Also note that, in the description above, as limits commute with the right adjoints $\iota$, to show that the map $\ol\iota_\lbd\iota_\lbd\iota_\lbd^*\ol\iota_\lbd^*\lto\lim_{\mu\to\lbd}\ol\iota_\mu\ol\iota_\mu^*$ is an isomorphism, it suffices to show
    \[\Id_{\CC^{<\lbd}}\simlto \lim_{\mu\to\lbd}\iota_\mu^{<\lbd}\iota_\mu^{<\lbd}{}^*\] 
    is an isomorphism. This is proved by using induction on the poset, hence the conditions on $\Lbd$. 

    First let us do the base case(s). The trivial base case is when $\lbd$ admits no smaller $\mu$ (i.e. `minimal'), in which case $\CC^{<\lbd}$ is the trivial category, so that the claim is trivially true. The next base case is when $\lbd$ admits `one layer' of smaller $\mu$, namely there are minimal (and therefore mutually incomparable) $\mu_i<\lbd$. Then $\CC^{<\lbd}=\prod_i \CC^{\le\mu_i}$, and as $\iota_{\mu}^{<\lbd}{}^*$ can be easily seen to be the projection, again the claim is trivially true. 

    Now we do the inductive step. Let $\mu_i$ be the maximal incomparable elements lower than $\lbd$, possible since $\Lbd$ is locally-finite. We suppress the subscript when we talk about just one. Let $\CC^{\ll\lbd}\coloneqq\wan{\CC^{<\mu_i}}_i$ be the category generated by the $\CC^{<\mu_i}$. Consider the diagram:
    \begin{center}
        \begin{tikzcd}
        \CC^{\ll\lbd} \arrow[r, "\iota_{<\lbd}", hook]                                   & \CC^{<\lbd} \arrow[r, "\jota^{<\lbd}"]                                   & \CC^{<\lbd}/\CC^{\ll\lbd} \arrow[d, dashed, bend right=49, shift right] \\
        \CC^{<\mu} \arrow[u, "\iota_{\mu}^{\ll\lbd}", hook] \arrow[r, "\iota_\mu", hook] & \CC^{\le\mu} \arrow[r, "\jota^\mu"] \arrow[u, "\iota_\mu^{<\lbd}", hook] & \CC^{=\mu} \arrow[u, dashed]                                           
        \end{tikzcd}
    \end{center}
    Note that each row is a recollement. To save space we have not included the usual left adjoints. We claim there is an upward dotted map $\wt\iota_\mu^{<\lbd}\colon \CC^{=\mu}\lto \CC^{<\lbd}/\CC^{\ll\lbd}$ and a left adjoint $\wt\iota_\mu^{<\lbd}{}^*$ such that
    \[\wt\iota_\mu^{<\lbd}\jota^\mu=\jota^{<\lbd}\iota_\mu^{<\lbd},\qquad \jota^\mu_!\wt\iota_\mu^{<\lbd}{}^*=\iota_\mu^{<\lbd}{}^*\jota^{<\lbd}_!,\] 
    namely that the diagrams
    \begin{center}
\begin{tikzcd}
\CC^{<\lbd} \arrow[r, "\jota^{<\lbd}"]                                   & \CC^{<\lbd}/\CC^{\ll\lbd}                            \\
\CC^{\le\mu} \arrow[r, "\jota^\mu"] \arrow[u, "\iota_\mu^{<\lbd}", hook] & \CC^{=\mu} \arrow[u, "\wt\iota_\mu^{<\lbd}", dashed,swap]
\end{tikzcd},\qquad 
\begin{tikzcd}
\CC^{<\lbd} \arrow[d, "\iota_\mu^{<\lbd}{}^*"'] & \CC^{<\lbd}/\CC^{\ll\lbd} \arrow[l, "\jota^{<\lbd}_!"'] \arrow[d, "\wt\iota_\mu^{<\lbd}{}^*", dashed] \\
\CC^{\le\mu}                                    & \CC^{=\mu} \arrow[l, "\jota^\mu_!"']                                                                 
\end{tikzcd}
\end{center}
commute. Let us construct them as thus:
\begin{align*}
    \wt\iota_\mu^{<\lbd}&=\jota^{<\lbd}\iota_\mu^{<\lbd}\jota^\mu_*,\\ 
    \wt\iota_\mu^{<\lbd}{}^*&=\jota^\mu\iota_\mu^{<\lbd}{}^*\jota^{<\lbd}_!.
\end{align*}
Let us check that so-defined the diagrams commute. Indeed, one has from recollement and $\iota_!\iota^!\to\Id\to\jota_*\jota^*\to[+1]$ that
\begin{align*}
    \wt\iota_\mu^{<\lbd}\jota^\mu=\jota^{<\lbd}\iota_\mu^{<\lbd}\jota^\mu_*\jota^\mu=\on{Cofib}(\underbrace{\jota^{<\lbd}\iota_\mu^{<\lbd}\iota_\mu}_{=\jota^{<\lbd}\iota_{<\lbd}\iota_\mu^{\ll\lbd}=0}\iota_\mu^!\lto\jota^{<\lbd}\iota_\mu^{<\lbd})=\jota^{<\lbd}\iota_\mu^{<\lbd},
\end{align*}
which checks the first square, and similarly
\begin{align*}
    \jota^\mu_!\jota^\mu\iota_\mu^{<\lbd}{}^*\jota^{<\lbd}_!=\Fib(\iota_\mu^{<\lbd}{}^*\jota^{<\lbd}_!\lto \iota_\mu\underbrace{\iota_\mu^*\iota_\mu^{<\lbd}{}^*\jota^{<\lbd}_!}_{=\iota_\mu^{\ll\lbd}{}^*\iota_{<\lbd}^*\jota^{<\lbd}_!=0})=\iota_\mu^{<\lbd}{}^*\jota^{<\lbd}_!
\end{align*}
checks the second square. So the maps we advertise do exist. Now, the different $\CC^{=\mu_i}$ form minimal mutually incomparable subcategories inside $\CC^{<\lbd}/\CC^{\ll\lbd}$, so by the base case of the induction we know
\[\Id_{\CC^{<\lbd}/\CC^{\ll\lbd}}\simlto\lim_\mu \wt\iota_\mu^{<\lbd}\wt\iota_{\mu}^{<\lbd}{}^*.\] 
But now note, since $\Id\simlto\jota^\mu\jota^\mu_!$, we know
\begin{align*}
    \lim_\mu \wt\iota_\mu^{<\lbd}\wt\iota_{\mu}^{<\lbd}{}^*=\lim_\mu \wt\iota_\mu^{<\lbd}\jota^\mu\jota^\mu_!\wt\iota_\mu^{<\lbd}{}^*=\lim_\mu \jota^{<\lbd}\iota_\mu^{<\lbd}\iota_\mu^{<\lbd}{}^*\jota^{<\lbd}_!=\jota^{<\lbd}\lim_\mu(\iota_\mu^{<\lbd}\iota_\mu^{<\lbd}{}^*)\jota^{<\lbd}_!,
\end{align*}
so that combined with the above we know the map
\[\jota^{<\lbd}_!\simlto \lim_\mu\iota_\mu^{<\lbd}\iota_\mu^{<\lbd}{}^*\jota^{<\lbd}_!\] 
is an equivalence. Moreover, by induction, we know
\[\iota_{<\lbd}\simlto \lim_\mu \iota_\mu^{<\lbd}\iota_\mu^{<\lbd}{}^*\iota_{<\lbd}\] 
is also an equivalence. Together, this shows that
\[\Id_{\CC^{<\lbd}}\simlto\lim_\mu \iota_\mu^{<\lbd}\iota_\mu^{<\lbd}{}^*\] 
is an equivalence due to the following general lemma.

The general lemma is that in a recollement situation $\CC^<\to\CC^\le\to\CC^=$, if $\CF$ is an endofunctor of $\CC^\le$ which is the identity on the image of $\iota_*=\iota_!$ and $\jota_!$, then $\CF$ must be the identity. This is clearly true because the triangle $\jota_!\jota^!\to\Id\to\iota_*\iota^*$ gives
\[\CF\jota_!\jota^!\lto\CF\lto\CF\iota_*\iota^*\oslto{+1},\] 
and since $\CF\jota_!=\jota_!$ and $\CF\iota_*=\iota_*$ we conclude by the uniqueness of cones.

Now, if the poset $\Lbd$ is down-finite, then this induction process will terminate. If it is instead locally finite and eventually disjoint-totally-ordered, then eventually the induction process will turn into proving this for a totally ordered set, where it is easily seen to be true, since if $\Lbd$ is totally ordered and locally finite then for any $\lbd$ there is a single maximal $\mu<\lbd$, whereupon $\CC^{<\lbd}=\CC^{\le\mu}$.
            \end{proof}
            \end{adjustwidth}\vspace{0.5em}

\begin{remark}
    The above setup, with slight differences, has appeared for instance in \cite{ayala2022stratified}, Remarks 1.3.12 and 1.3.13 on pages 22--23. Our set-ups are slightly different (for example they require subcategories to have a tower of two right adjoints), but philosophically they should be roughly the same. 
\end{remark}
\begin{remark}
    Also note that if $\Lbd$ is down-finite (namely that for any $\lbd$ the set $\{\mu:\mu\le\lbd\}$ is finite), then $\Lbd$ automatically satisfies the conditions of the theorem. This is for instance the case with (a block of) category $\CO$.
\end{remark}

Let us say a few words about the maps in the spectral sequence in Theorem \ref{thm:folklorefiltration} above. The maps on the first page $E_1$ are coming from the functoriality of the (co)cone; indeed, one has a $\BZ$-filtration
\[\cdots\lto \bigoplus_{\ell(\lbd)=i}\calid(\lbd)\lto\bigoplus_{\ell(\lbd)=i+1}\calid(\lbd)\lto\cdots,\] 
where the morphisms are as constructed in the proof above, so that applying the functoriality of the fiber to
\[\bigoplus_{\ell(\lbd)=i}\calid(\lbd)\lto\bigoplus_{\ell(\lbd)=i+1}\calid(\lbd)\lto\bigoplus_{\ell(\lbd)=i+2}\calid(\lbd)\]
we obtain the exact triangle
\[\bigoplus_{\ell(\lbd)=i}\gr^\lbd\calid \lto \Fib\pr*{\bigoplus_{\ell(\lbd)=i}\calid(\lbd)\to\bigoplus_{\ell(\lbd)=i+2}\calid(\lbd)} \lto \bigoplus_{\ell(\lbd)=i+1}\gr^\lbd\calid \oslto{+1};\]
rotating this gives a map
\[\bigoplus_{\ell(\lbd)=i+1}\gr^\lbd\calid\lto \bigoplus_{\ell(\lbd)=i}\gr^\lbd\calid[1],\]
which under the action of taking $\pi$ gives the differential maps of $E_1$.

\subsection{Application to triangular-based algebras}
Now let us return to the setting of the triangular based algebra $A$, which has lowest weight notation. Note that the recollements described in Section \ref{section:setting} give rise to a $\Theta^\op$-filtration of the stable $\infty$-category $\D^-\on\Mod A$ where
\[(\D^-\on\Mod A)^{\le\theta}=\D^-\on\Mod A^{\ge\theta}.\] 
It is easy to see that
\[(\D^-\on\Mod A)^{<\theta}=\D^-\on\Mod A^{>\theta}=\D^-\on\Mod(\mfrac{A^{\ge\theta}}{A^{\ge\theta}e^\theta A^{\ge\theta}}),\] 
so that the Verdier quotient is $\D^-\on\Mod A^\theta$. Then, by Theorem \ref{thm:folklorefiltration}, there exists a $\Theta$-filtration on the identity functor whose associated graded is
\begin{align*}
    \gr^\theta\calid=\ol\iota_\theta\jota^\theta_!\jota^\theta\ol\iota_\theta^*=\Delta(\theta)\otimes_{A^\theta}\rhom_A\pr*{\Delta(\theta),\sq^\dag}^*,
\end{align*}
where we recall the shorthand $\Delta(\theta)\coloneqq \bigoplus_{\lbd\in\theta}\Delta_\lbd^{\ol{l_\lbd(\theta)}}$. Moreover, \ref{thm:folklorefiltration} tells us that provided there is a length function $\ell\colon \Theta\lto\BZ$, there is a spectral sequence (functorial in the input $\sq$)
\[E_1^{p,q}=\bigoplus_{\ell(\theta)=-p}\Delta(\theta)\otimes_{A^\theta}\Ext^{-(p+q)}_A(\Delta(\theta),\sq^\dag)^* \specseqimplies E_\infty^{p,q}=\gr^{-p}H^{p+q}(\sq).\]

To recover Theorem \ref{thm:nbspecseq} for the nilBrauer algebra, consider the poset of weights $\Theta$ for $\NB$ as two disjoint totally ordered sets, namely $0< 2< 4<\cdots$ and $1< 3< 5< \cdots$, equipped with a length function in which $\ell(2k)=\ell(2k+1)=k$. Then specializing to a block (either even or odd) recovers \ref{thm:nbspecseq}.

\begin{remark}
    One could just as well consider $\Theta$ for nilBrauer to be totally ordered as $0<1<2<\cdots$; in this case the spectral sequence would look different but be equivalent.
\end{remark}

To describe this spectral sequence in more detail for the nilBrauer setting, it is then crucial to investigate the Ext groups appearing on the first page.

\section{Standard Ext Groups, Part 1: Nilcohomology}\label{section:nilcohomo}
In order to obtain a BGG resolution of $L$ using this philosophy of reconstructing the identity functor, one needs to compute these Ext groups, $\Ext_A^\blt(\Delta(\theta),L)$. Either way, the Ext groups between standard modules and simples are interesting alone, being the analogue of the Kazhdan-Lusztig problem. As they are Ext groups out of standard modules rather than proper standard modules, we shall call them ``standard Ext groups'', or ``big Ext groups'' (after `big Verma'). We attack this by emulating Lie algebra nilcohomology.

Recall that $\Ext^\blt(\Delta_\lbd,\sq)=H^\blt(\nlie^+:\sq)^\lbd$ in category $\CO$, where the key was to write this as $\rhom_{\nlie^+}(\BC_0,\sq)$, which admits an action of $\hlie$. This in turn can be computed with the Koszul resolution of $\BC_0$, which importantly is a finite resolution (because it is built from wedge powers of $\nlie^+$, which is finite-dimensional). This is the story we will try to emulate.

\subsection{The nilalgebra and nilcohomology}
In the below story, we have written $A$ instead of $\NB$ to indicate that this procedure is general. Let $\BK=\bigoplus \bk e^\theta$, and let the ``nilalgebra'' be defined as
\[A^-\coloneqq \BK \tY=\BK\oplus\mathop{\bigoplus\bigoplus}_{\psi\not\ge\theta} e^\psi \BC \te Y_+ e^\theta,\]
where the subscript in $\tY_+$ indicates that the diagram $e^\theta$ is not included. Note that in the nilBrauer case one may as well take the bottom index to run over $\psi<\theta$ of the same parity. This is the analogue of $U\nlie^+$. Of course one must check on a case-by-case basis that, so defined, this $A^-$ is actually a subalgebra of $A$, i.e. that it is closed under multiplication. This is typically not difficult, and we check it now for the nilBrauer:
\begin{LEM}
    So defined, $\NB^-$ is a subalgebra of $\NB$.
\end{LEM}
\vspace{-0.5em}
\begin{PRF}
    It is important to note here that, due to the way the PBW basis theorem is stated for nilBrauer, the diagrams of $\BC\tY_+$ have dots in `preferred' locations: it is our choice to have dots on the far left rather than the right of each cap. One needs to check that $\NB^-$, so defined, is closed under multiplication: namely that when two diagrams, each with dots only at the left base of each cap, are stacked with each other, one obtains a linear combination of diagrams also with dots only at the left base of each cap. This is true because the dots in the top diagram must connect with propagating strings of the bottom diagram, but propagating strings of the bottom diagram can only cross with caps by assumption (and not other propagating strings). Then the relation to observe is simply (here the orange lines are not part of the diagrammatics and are visual aids placed to make clear where the diagram starts/ends)
    \[\hackcenter{\begin{tikzpicture}[scale=0.375]
    \draw[thick,orange] (0,0) -- (4,0);
    \draw[thick,orange] (0,4) -- (4,4);
   \draw (2.5,0) arc (0:180:1);
    \draw (1.5,0)--(1.5,2);
    \draw (3.5,0)--(3.5,2);
    \draw (3.5,2) arc (0:180:1);
    \fill (1.5,0.5) circle (5pt);
\end{tikzpicture}}
-
\hackcenter{\begin{tikzpicture}[scale=0.375]
    \draw[thick,orange] (0,0) -- (4,0);
    \draw[thick,orange] (0,4) -- (4,4);
    \draw (2.5,0) arc (0:180:1);
    \draw (1.5,0)--(1.5,2);
    \draw (3.5,0)--(3.5,2);
    \draw (3.5,2) arc (0:180:1);
    \fill (1.5,1.5) circle (5pt);
\end{tikzpicture}}
=
\hackcenter{\begin{tikzpicture}[scale=0.375]
    \draw[thick,orange] (0,0) -- (4,0);
    \draw[thick,orange] (0,4) -- (4,4);
    \draw (1.8,0) arc (0:180:0.8);
    \draw (3.8,0) arc (0:180:0.8);
\end{tikzpicture}}
-
\hackcenter{\begin{tikzpicture}[scale=0.375]
    \draw[thick,orange] (0,0) -- (4,0);
    \draw[thick,orange] (0,4) -- (4,4);
    \draw (2.5,0) arc (0:180:0.5);
    \draw (0.5,0)--(0.5,2);
    \draw (1.5,2) arc (0:180:0.5);
    \draw (2.5,2) arc (0:-180:0.5);
    \draw (3.5,2) arc (0:180:0.5);
    \draw (3.5,0)--(3.5,2);
\end{tikzpicture}},
\]
which is a consequence of the most important defining relation of nilBrauer, that
\[\hackcenter{\begin{tikzpicture}[scale=0.375]
    \draw (0,0)--(2,2);
    \draw (2,0)--(0,2);
    \fill (0.5,1.5) circle (5pt);
\end{tikzpicture}}
-
\hackcenter{\begin{tikzpicture}[scale=0.375]
    \draw (0,0)--(2,2);
    \draw (2,0)--(0,2);
    \fill (1.5,0.5) circle (5pt);
\end{tikzpicture}}
=
\hackcenter{\begin{tikzpicture}[scale=0.375]
    \draw (0,0)--(2,2);
    \draw (2,0)--(0,2);
    \fill (0.5,0.5) circle (5pt);
\end{tikzpicture}}
-
\hackcenter{\begin{tikzpicture}[scale=0.375]
    \draw (0,0)--(2,2);
    \draw (2,0)--(0,2);
    \fill (1.5,1.5) circle (5pt);
\end{tikzpicture}}
=
\hackcenter{\begin{tikzpicture}[scale=0.375]
    \draw (0,0) to[out=60,in=-60] (0,2);
    \draw (2,0) to[out=120,in=-120] (2,2);
\end{tikzpicture}}
-
\hackcenter{\begin{tikzpicture}[scale=0.375]
    \draw (0,0) to[out=60,in=120] (2,0);
    \draw (0,2) to[out=-60,in=-120] (2,2);
\end{tikzpicture}};
\]
by rearranging the first diagram as
\[\hackcenter{\begin{tikzpicture}[scale=0.375]
    \draw[thick,orange] (0,0) -- (4,0);
    \draw[thick,orange] (0,4) -- (4,4);
   \draw (2.5,0) arc (0:180:1);
    \draw (1.5,0)--(1.5,2);
    \draw (3.5,0)--(3.5,2);
    \draw (3.5,2) arc (0:180:1);
    \fill (1.5,0.5) circle (5pt);
\end{tikzpicture}}=\hackcenter{\begin{tikzpicture}[scale=0.375]
    \draw[thick,orange] (0,0) -- (4,0);
    \draw[thick,orange] (0,4) -- (4,4);
    \draw (3.5,0) arc (0:180:1);
    \draw (2.5,0)--(2.5,2);
    \draw (0.5,0)--(0.5,2);
    \draw (2.5,2) arc (0:180:1);
    \fill (2.5,0) ++(135:1) circle (5pt); 
\end{tikzpicture}},\]
one can rearrange this equation into
\[
\hackcenter{\begin{tikzpicture}[scale=0.375]
    \draw[thick,orange] (0,0) -- (4,0);
    \draw[thick,orange] (0,4) -- (4,4);
    \draw (2.5,0) arc (0:180:1);
    \draw (1.5,0)--(1.5,2);
    \draw (3.5,0)--(3.5,2);
    \draw (3.5,2) arc (0:180:1);
    \fill (1.5,1.5) circle (5pt);
\end{tikzpicture}}
=
\hackcenter{\begin{tikzpicture}[scale=0.375]
    \draw[thick,orange] (0,0) -- (4,0);
    \draw[thick,orange] (0,4) -- (4,4);
    \draw (3.5,0) arc (0:180:1);
    \draw (2.5,0)--(2.5,2);
    \draw (0.5,0)--(0.5,2);
    \draw (2.5,2) arc (0:180:1);
    \fill (2.5,0) ++(135:1) circle (5pt); 
\end{tikzpicture}}
-\hackcenter{\begin{tikzpicture}[scale=0.375]
    \draw[thick,orange] (0,0) -- (4,0);
    \draw[thick,orange] (0,4) -- (4,4);
    \draw (1.8,0) arc (0:180:0.8);
    \draw (3.8,0) arc (0:180:0.8);
\end{tikzpicture}}
+
\hackcenter{\begin{tikzpicture}[scale=0.375]
    \draw[thick,orange] (0,0) -- (4,0);
    \draw[thick,orange] (0,4) -- (4,4);
    \draw (2.5,0) arc (0:180:0.5);
    \draw (0.5,0)--(0.5,2);
    \draw (1.5,2) arc (0:180:0.5);
    \draw (2.5,2) arc (0:-180:0.5);
    \draw (3.5,2) arc (0:180:0.5);
    \draw (3.5,0)--(3.5,2);
\end{tikzpicture}}.
\]
Hence, when multiplying a diagram with another diagram, the dots on the top diagram can be moved to be at the left base of the resulting product diagram, at the cost of error terms which (e.g. by induction) are also in $\BC\tY$. This shows $\BC\tY$ is closed under multiplication, as desired.
\end{PRF}

To return to the general picture, this $A^-$ then has an ideal 
\[I\coloneqq A^-_+\coloneqq \mathop{\bigoplus\bigoplus}_{\psi\not\ge\theta}e^\psi \BC \te Y_+ e^\theta,\] 
which is analogous to the positive degree elements of $U\nlie^+$. Once it is shown that $A^-$ is a subalgebra of $A$, it is evidently automatic that $I$ is an ideal. Note that $A^-/I=\BK$, and that there is an action
\[A^-\actson \bk e^\theta\]
where all of $I$ acts by zero and $\BK$ acts in the obvious way by multiplication. We can consider the obvious modification of the bar resolution in order to resolve $\bk e^\theta$ with free $A^-$-modules (here all tensors are over $\BK$):
\[\cdots\lto A^-\otimes I^{\otimes k}e^\theta\lto\cdots\lto A^-\otimes I\otimes Ie^\theta\lto A^-\otimes Ie^\theta\lto A^- e^\theta\lto \bk e^\theta,\]
where the maps are (for brevity we write $(a_1,\cdotsc,a_k)$ for $a_1\otimes\cdots\otimes a_k$)
\begin{align*}
    \pd_k\colon A^-\otimes_\BK I^{\otimes k} e^\theta&\lto A^-\otimes_\BK I^{\otimes k-1}e^\theta\\
    x\otimes (a_1,\cdotsc, a_ke^\theta)&\lmto xa_1\otimes (a_2,\cdotsc,a_ke^\theta)+\sum_{i=1}^{k-1} (-1)^i x\otimes (a_1,\cdotsc,a_ia_{i+1},\cdotsc, a_ke^\theta).
\end{align*}
This can be shown to be exact by using the obvious modification of the usual homotopy: send $x\otimes a_1\otimes\cdots\otimes a_k$ to $1\otimes \pi_+(x)\otimes a_1\otimes\cdots\otimes a_k$, where $\pi_+\colon A^-\lto I$ is projection away from $A^-_0=\BK$.

Now, note 
\[\hom_A(A^{\ge\theta}e^\theta,M)\cong \hom_{A^-}(\bk e^\theta, M)\tag{1}\label{eqn:vermahom}\] 
This is because a map on the LHS is determined by where it sends $e^\theta$, and since the map must commute with all of $A$, it must send $e^\theta$ to something that is killed by anything passing through $e^\psi$ for $\psi\not\ge\theta$. A similar analysis of the RHS shows these two are isomorphic, and both are isomorphic to $(M^{A^-})^\theta$, where
\[M^{A^-}\coloneqq \{v\in M:I\cdot v=0\}\] 
is the submodule of vectors killed by $I$. This is all in analogy\footnote{The fact that the plus sign and minus sign are switched is a reflection of the fact that the nilBrauer representation theory is lowest weight rather than highest weight.} to the observation from Lie theory of $\hom_\CO(\Delta_\lbd,M)=\hom_{\nlie^+}(\BC_\lbd, M)=M^{\nidl^+}$. Note well that the left $A^\theta$-action on the LHS of Equation (\ref{eqn:vermahom}), which came from the right $A^\theta$-action on $A^{\ge\theta}e^\theta$, becomes the left action on the RHS inherited from $M$. The issue with this is that the $A^\theta$-action on the RHS becomes `lost' when we take derived functors, in the sense that the above complex, $A^-\otimes I^{\otimes\blt}e^\theta$, does not have a right $A^\theta$-action. Secretly it must still be there, but one cannot compute it using this complex. Hence we may define (the $\theta$-weight space of) the ``nilcohomology'' as
\[H^\blt(A^-:M)^\theta\coloneqq \rhom_A(A^{\ge\theta}e^\theta,M)\oscong{\sf{Vec}}\rhom_{A^-}(\bk e^\theta,M)\]
which can be computed as the cohomology of the complex
\[\hspace{-18em}0\lto \hom_{A^-}(A^-e^\theta,M)\lto \hom_{A^-}(A^-\otimes Ie^\theta,M)\lto\cdots \] 
\[\iddots\] 
\[\hspace{15em}\cdots\lto \hom_{A^-}(A^-\otimes I^{\otimes k}e^\theta,M)\lto \hom_{A^-}(A^-\otimes I^{\otimes k+1}e^\theta,M)\lto \cdots\] 
This complex can be rewritten using tensor-Hom as
\[\hspace{-20em}0\lto \hom_{\BK}(\bk e^\theta,M)\lto \hom_{\BK}(Ie^\theta,M)\lto\cdots\] 
\[\hspace{-3em}\iddots\] 
\[\hspace{18em}\cdots\lto \hom_{\BK}(I^{\otimes k}e^\theta,M)\lto \hom_{\BK}(I^{\otimes k+1}e^\theta,M)\lto \cdots\]
where each differential $\d^k(f)=f(\pd_{k+1}\sq)$ can be simply found to be 
\begin{align*}
    \d^k\colon \hom_\BK(I^{\otimes k}e^\theta,M)&\lto \hom_\BK(I^{\otimes k+1}e^\theta,M)\\
    f&\lmto \biggpr{ (a_1,\cdotsc,a_{k+1}e^\theta)\mapsto a_1f(a_2,\cdotsc,a_{k+1}e^\theta)+\sum_{i=1}^k(-1)^i f(a_1,\cdotsc,a_ia_{i+1},\cdotsc,a_{k+1}e^\theta)  }.
\end{align*}

    This story can be made to even more closely resemble the Lie algebra one by taking a more `uniform' approach -- consider the action $A^-\actson \BK$ instead of the action on $\bk e^\theta$; then $\rhom_{A^-}(\BK,M)$ has a left action of $\BK$, and we can take the weight spaces with respect to this action to get each 
    \[H^\blt(A^-:M)^\theta=e^\theta\Ext_{A^-}^\blt(\BK,M)=\Ext_A^\blt(A^{\ge\theta}e^\theta,M)=\Ext_A^\blt(\Delta(\theta),M).\]
    This $\rhom_{A^-}(\BK,M)$ then deserves the name ``nilcohomology''; we will call this
    \[H^\blt(A^-:M)\coloneqq \Ext_{A^-}^\blt(\BK,M).\]
    This is computed with the complex
    \[\cdots\lto A^-\otimes I^{\otimes k}\lto\cdots\lto A^-\otimes I\otimes I\lto A^-\otimes I\lto A^- \lto \BK ,\]
    which `weight-specializes' to the complex $A^-\otimes I^{\otimes\blt}e^\theta$ considered above. More explicitly, 
    \[H^\blt(A^-:M)\simeq \hom_\BK(I^{\otimes\blt},M),\] 
    where the differentials are identical to $\d^k$ computed above except erasing each appearance of $e^\theta$. Worth saying is that, also in analogy with Lie theory, $\rhom_{A^-}^\blt(\BK,M)$ is the right derived functor of the invariant functor $M\mapsto M^{A^-}$ which picks out the set of vectors killed by $I$. So again one has a cohomology theory obtained from requiring a `nilalgebra' to act by zero.
\begin{remark}
    A warning though: this incarnation of the nilcohomology, using the bar complex above, loses the left $A^\theta$-action on $H^\blt(A^-:M)^\theta$ due to the lack of an obvious right $A^\theta$-action on $A^-\otimes I^{\otimes k}$. We will fix this later in this Section.
\end{remark}

From this resolution one can immediately say something weak about where the standard Ext groups vanish. As the tensors are over $\BK$, they must be permeable to the idempotents $e^\theta$, so that a tensor $a_1\otimes\cdots\otimes a_k\in I\otimes\cdots\otimes I$ must have the beginning and the ending idempotents of each consecutive $a_i$ match up. In other words, it must be that $a_i=a_i e^{\theta_i}$ and $a_{i+1}=e^{\theta_i}a_{i+1}$ for each $i$. In the nilBrauer case, this means for instance that 
\[\text{$I^{\otimes k}e^\theta=0$ if $k>\theta/2$},\] 
since elements of $I$ decrease the weight by at least 2. In this respect this story resembles that of Lie algebra nilcohomology even further. This means in particular we have the (very weak) bound
\begin{LEM}
    \[\Ext_\NB^{>\theta/2}(\NB^{\ge\theta}e^\theta,M)=0,\]
    so in particular
    \[\Ext_\NB^{>\theta/2}(\Delta_\theta,M)=0.\]
\end{LEM}
The above is true for all $M$, and the results can be improved depending on the lowest weight of $M$, since the Homs must moreover respect the $\BK$-action on the first tensor factor (i.e. maps must send $a_1\otimes\cdots\otimes a_k$ to the weight space $M^\psi$ if $a_1=e^\psi a_1$):
\begin{LEM}
    If $M$ has lowest weight $\psi$ of the same parity as $\theta$ (all such $\rhom$'s vanish if they are not of the same parity), then
    \[\Ext_\NB^{>\frac{\theta-\psi}{2}}(\NB^{\ge\theta}e^\theta,M)=0,\] 
    and in particular 
    \[\Ext_\NB^{>\frac{\theta-\psi}{2}}(\Delta_\theta,M)=0.\] 
\end{LEM}
These cohomology vanishing results indicate that the main obstruction to $\on\Mod \NB$ having finite global dimension is the fact that there are infinitely many standard objects (corresponding to infinitely many idempotents). If one were to `truncate' the algebra $\NB$ in a stupid manner, for example simply by deleting all diagrams that pass through more than $\theta$ strands, then one would obtain an algebra whose module category had finite global dimension.

It is perhaps worth pointing out, as a special case of the constructions above, that for the module $A^-\actson \BK$, the nilcohomology $H^\blt(A^-:\BK)=\rhom_{A^-}^\blt(\BK,\BK)$ can be computed as 
\begin{align*}
    \Ext^\blt_{A^-}(\BK,\BK)&=H^\blt\big( \hom_{A^-}(A^-\otimes I^{\otimes\blt},\BK)\big)\\
    &=H^\blt\big( \hom_\BK(I^{\otimes\blt},\BK)\big)\\
    &=H^{-\blt}(I^{\otimes\blt}). 
\end{align*} 
Also worth noting is that, in the case that $A^-$ is Koszul, we know that
\[\Ext^\blt_{A^-}(\BK,\BK)\cong A^{-,!,\op}.\] 
Using the right actions of $\BK\ractson\BK$ for both entries to obtain left and right $\BK$-actions on the Ext group, we then have
\[e^\phi A^{-,!}e^\theta=e^\theta\Ext^\blt_{A^-}(\BK,\BK)e^\phi=\Ext^\blt_{A^-}(\bk e^\theta,\bk e^\phi).\] 

\subsection{Then the nilcohomology assumes his multi-armed form}\label{subsec:multiarm}
We had remarked earlier that our earlier construction of nilcohomology using $A^-$ fails to reproduce the $A^\theta$-action. In this section we remedy this. This is accomplished by defining another subalgebra $A^\flat$, which is roughly a `bigger' version of $A^-$; computing nilcohomology with this $A^\flat$ will reproduce the $A^\theta$-action. One might wonder why the earlier construction was necessary at all then -- the reason is because the former is easier to handle, and indeed one at least has a chance of proving that $A^-$ is Koszul, as it at least has a semisimple degree 0 part. The algebra we will now define does not admit Koszul methods as simply (as far as we can tell). 

Define the locally unital subalgebra of $A$
\[A^\flat\coloneqq \bigoplus_{\theta\le\phi} e^\theta A e^\phi.\] 
From this definition, it is clear that this is a subalgebra of $A$. Note that cup-caps are allowed as elements of this subalgebra in the case of $A=\NB$. Also define the locally unital subalgebra
\[A^\nat\coloneqq \bigoplus_\theta e^\theta A e^\theta,\]
which is a quotient of $A^\flat$ by the ideal
\[A^\flat_+\coloneqq \bigoplus_{\theta<\phi} e^\theta A e^\phi.\] 
Contrast $A^\nat$ with 
\[A^\circ\coloneqq \bigoplus_\theta A^\theta,\] 
which is not a subalgebra, but rather a subquotient. There is a quotient map
\[0\lto \kidl=e^\theta\wan{e^\phi:\phi\not\ge\theta}e^\theta\lto A^\nat\lto A^\circ\lto 0,\]
and also
\[0\lto I^\flat\lto A^\flat\lto A^\circ\lto 0,\] 
where 
\[I^\flat\coloneqq \wan{\kidl,A^\flat_+}.\] 

In the case of nilBrauer, we can put an explicit grading on $\NB^\flat$. Let
\[\NB^\flat_i\coloneqq \bigoplus_\theta e^{\theta-2i}\NB e^\theta,\] 
which is to say the $i$-th graded part are the diagrams with $i$ more caps than cups. Under this grading $\NB^\flat_0=\NB^\nat$, and 
\[\NB^\flat_+=\bigoplus_{i\ge 1} \NB^\flat_i.\]


Then one makes a similar observation to before, that
\[\hom_A(A^{\ge\theta}e^\theta,M)=\hom_{A^-}(\bk e^\theta,M)=\hom_{A^\flat}(A^\theta,M).\]
The `uniform' (with respect to weights) statement of this would be
\[\hom_A\pr*{\bigoplus_\theta A^{\ge\theta}e^\theta,M}=\hom_{A^-}(\BK,M)=\hom_{A^\flat}(A^\circ,M).\] 
The benefit of this is that now $\hom_{A^\flat}(A^\theta,M)$ admits a left $A^\theta$-action coming from the right action on the first input. Note that we can also think of this action as an $A^\nat$-action coming from $A^{\circ}\ractson A^\nat$, which moreover factors to an $A^{\circ}$-action. To compute the derived functors, we utilize also a bar resolution. This time consider the resolution
\[\cdots\lto A^\flat\otimes_{\BK} I^{\flat,\otimes_{A^\nat}k}\lto\cdots\lto A^\flat\otimes_{\BK} I^\flat\otimes_{A^\nat}I^\flat\lto A^\flat\otimes_{\BK} I^\flat\lto A^\flat\lto A^{\circ}, \]
which would weight-specialize to
\[\cdots\lto A^\flat\otimes_{\BK} I^{\flat,\otimes_{A^\nat}k}e^\theta \lto\cdots\lto A^\flat\otimes_{\BK} I^\flat\otimes_{A^\nat}I^\flat e^\theta\lto A^\flat\otimes_{\BK} I^\flat e^\theta\lto A^\flat e^\theta\lto A^\theta. \]
This resolution is by $(A^\flat,A^\nat)$-bimodules which are free on the left (due to the left-most tensor being over $\BK$), and moreover there is a right $A^\nat$-action coming from $I^\flat\ractson A^\nat$. Note well that the first tensor in each term is over $\BK$, while subsequent tensors are over $A^\nat$. Then each term of the complex which computes Ext
\[\hom_{A^\flat}(A^\flat,M)\lto \hom_{A^\flat}(A^\flat\otimes_{\BK} I^\flat,M)\lto \hom_{A^\flat}(A^\flat\otimes_{\BK} I^\flat\otimes_{A^\nat} I^\flat,M)\lto\cdots\] 
admits a left $A^\nat$-action, and because this $A^\nat$-action is respected by the complex and because it coincides with the $A^\nat$-action on $\hom_{A^\flat}(A^{\circ},M)$ which factors through $A^{\circ}$, when one takes (co)homology of this complex to get the Ext groups, this $A^\nat$-action on the Ext groups must also factor through to a $A^{\circ}$-action. Specializing to a weight, we get an incarnation of nilcohomology which describes the action of the Cartans $A^\theta$. When we wish to emphasize that this $A^\flat$-resolution of $A^{\circ}$ is what we are using, we may add a subscript $H^\blt_\flat(A^-:M)$ to nilcohomology. 


\section{Standard Ext Groups, Part 2: Koszulness}\label{section:koszul}
Recall that a Kostant module is one for which the nilcohomology/standard Ext groups are appropriately concentrated; for example, the finite-dimensional simple modules of category $\CO$ are Kostant modules. We wish to show $L_0=\bk e^0$ (at $t=0$) is a Kostant module. It turns out this will be the same as showing the algebra $A^-=\NB^-$ is Koszul. This section is devoted to the proof that $\NB^-$ is Koszul. But first we must do more set up; for instance, what is even the grading with respect to which Koszulness is proved?

\subsection{The Koszul Grading}
Recall $\NB^-=\BK\oplus\beau\mathop{\bigoplus\bigoplus}_{\psi<\theta} e^\psi \BC \tY_+ e^\theta$. Endow an additional grading on this algebra, appropriately called the ``cap grading'': 
\[\deg_{\te{cap}}(\te{diagram})=\#\te{ of caps}.\] 
We use the name `Koszul grading' interchangeably. Let
\[\NB^-_i\coloneqq \te{$\BC$-span of all cap-degree-$i$ $\tY$-diagrams},\] 
so that
\[\NB^-=\bigoplus_{i=0}^\infty \NB^-_i.\] 
Note that we have put the cap grading in the subscript; the cap grading induces a grading on Ext groups, which we shall later also denote in the subscript.
\begin{remark}
    Warning: there are now two gradings at play, the cap grading and the quantum grading (or `$q$-grading'). The former counts the number of caps, while the latter counts crossings and dots. There is also another sort of `grading' coming from the weights/idempotents $e^\theta$; note that once you know the number of caps and the weight you start at, the weight you end at is determined, so perhaps this is only half a grading. Later there will also be a homological grading. So this story is sort of triply-and-a-half graded. We will also later consider a `dot filtration', counting the total number of dots. 
\end{remark}

Maybe now is a good time to remind ourselves that the most important relation in $\NB$ is
\[\hackcenter{\begin{tikzpicture}[scale=0.375]
    \draw (0,0)--(2,2);
    \draw (2,0)--(0,2);
    \fill (0.5,1.5) circle (5pt);
\end{tikzpicture}}
-
\hackcenter{\begin{tikzpicture}[scale=0.375]
    \draw (0,0)--(2,2);
    \draw (2,0)--(0,2);
    \fill (1.5,0.5) circle (5pt);
\end{tikzpicture}}
=
\hackcenter{\begin{tikzpicture}[scale=0.375]
    \draw (0,0)--(2,2);
    \draw (2,0)--(0,2);
    \fill (0.5,0.5) circle (5pt);
\end{tikzpicture}}
-
\hackcenter{\begin{tikzpicture}[scale=0.375]
    \draw (0,0)--(2,2);
    \draw (2,0)--(0,2);
    \fill (1.5,1.5) circle (5pt);
\end{tikzpicture}}=\tikzthrough-\tikzcupcap.\]
As remarked earlier, it is our convention that dots appear on the far left of each cap. Hence we only need to slide dots across crossings when both strands of the crossing are parts of caps; for example, in 
\[\hackcenter{\begin{tikzpicture}[scale=0.375]
    \draw[thick,orange] (0,0) -- (4,0);
    \draw[thick,orange] (0,4) -- (4,4);
    \draw (3.5,0) arc (0:180:1);
    \draw (2.5,0)--(2.5,2);
    \draw (0.5,0)--(0.5,2);
    \draw (2.5,2) arc (0:180:1);
    \fill (2.5,0) ++(135:1) circle (5pt); 
\end{tikzpicture}}
-
\hackcenter{\begin{tikzpicture}[scale=0.375]
    \draw[thick,orange] (0,0) -- (4,0);
    \draw[thick,orange] (0,4) -- (4,4);
    \draw (2.5,0) arc (0:180:1);
    \draw (1.5,0)--(1.5,2);
    \draw (3.5,0)--(3.5,2);
    \draw (3.5,2) arc (0:180:1);
    \fill (1.5,1.5) circle (5pt);
\end{tikzpicture}}
=
\hackcenter{\begin{tikzpicture}[scale=0.375]
    \draw[thick,orange] (0,0) -- (4,0);
    \draw[thick,orange] (0,4) -- (4,4);
    \draw (1.8,0) arc (0:180:0.8);
    \draw (3.8,0) arc (0:180:0.8);
\end{tikzpicture}}
-
\hackcenter{\begin{tikzpicture}[scale=0.375]
    \draw[thick,orange] (0,0) -- (4,0);
    \draw[thick,orange] (0,4) -- (4,4);
    \draw (2.5,0) arc (0:180:0.5);
    \draw (0.5,0)--(0.5,2);
    \draw (1.5,2) arc (0:180:0.5);
    \draw (2.5,2) arc (0:-180:0.5);
    \draw (3.5,2) arc (0:180:0.5);
    \draw (3.5,0)--(3.5,2);
\end{tikzpicture}},
\]
If one of the strands is a propagating string, we keep the dot on the left. It is clear this is not affected by composition of diagrams. 

The first important observation is that so defined, this graded algebra is quadratic. Indeed, every diagram that has multiple caps can be written as the vertical composition of multiple diagrams each with a single cap. So there is a surjective map
\[{\bigotimes}^\blt\NB^-_1\lsurj \NB^-\] 
given by multiplication. More importantly, the kernel $\qidl$ of this map
\[0\lto\qidl\lto{\bigotimes}^\blt\NB^-_1\lsurj \NB^-\lto 0\] 
is a `quadratic ideal', meaning that it is generated by a subspace of $\NB^-_1\otimes\NB^-_1$. For example, up to dots, the relations inside $\NB^-$ look like
\begin{align*}
    {}\hackcenter{\begin{tikzpicture}[scale=0.375]
        \draw[thick,orange] (0,0)--(4,0);
        \draw[thick,orange] (0,4)--(4,4);
        \draw (0.5,0)--(0.5,1);
        \draw (1.5,0)--(1.5,1);
        \draw (0.5,1) arc (180:0:0.5);
        \draw (2.5,0)--(2.5,2);
        \draw (3.5,0)--(3.5,2);
        \draw (2.5,2) arc(180:0:0.5);
    \end{tikzpicture}
    }
    &=
    \hackcenter{\begin{tikzpicture}[scale=0.375]
        \draw[thick,orange] (0,0)--(4,0);
        \draw[thick,orange] (0,4)--(4,4);
        \draw (0.5,0)--(0.5,2);
        \draw (1.5,0)--(1.5,2);
        \draw (0.5,2) arc (180:0:0.5);
        \draw (2.5,0)--(2.5,1);
        \draw (3.5,0)--(3.5,1);
        \draw (2.5,1) arc(180:0:0.5);
    \end{tikzpicture}
    },\\
    {}\hackcenter{\begin{tikzpicture}[scale=0.375]
        \draw[thick,orange] (0,0)--(4,0);
        \draw[thick,orange] (0,4)--(4,4);
        \draw (0.5,0)--(0.5,1);
        \draw (2.5,0)--(2.5,1);
        \draw (0.5,1)arc(180:0:1);
        \draw (1.5,0)--(1.5,2.5);
        \draw (3.5,0)--(3.5,2.5);
        \draw (1.5,2.5)arc(180:0:1);
    \end{tikzpicture}
    }
    &=\hackcenter{\begin{tikzpicture}[scale=0.375]
        \draw[thick,orange] (0,0)--(4,0);
        \draw[thick,orange] (0,4)--(4,4);
        \draw (0.5,0)--(0.5,2.5);
        \draw (2.5,0)--(2.5,2.5);
        \draw (0.5,2.5)arc(180:0:1);
        \draw (1.5,0)--(1.5,1);
        \draw (3.5,0)--(3.5,1);
        \draw (1.5,1)arc(180:0:1);
    \end{tikzpicture}
    },\\
    {}\hackcenter{\begin{tikzpicture}[scale=0.375]
        \draw[thick,orange] (0,0)--(4,0);
        \draw[thick,orange] (0,4)--(4,4);
        \draw (0.5,0)arc(180:0:1.5);
        \draw (1.5,0)--(1.5,2.5);
        \draw (2.5,0)--(2.5,2.5);
        \draw (1.5,2.5)arc(180:0:0.5);
    \end{tikzpicture}
    }&=0.
\end{align*}
Translated into the form of a quadratic ideal, elements in $\qidl$ would look like
\begin{align*}
    {}\hackcenter{\begin{tikzpicture}[scale=0.25]
        \draw[thick,orange] (0,0)--(4,0);
        \draw[thick,orange] (0,4)--(4,4);
        \draw (0.5,0)arc(180:0:1.5);
    \end{tikzpicture}
    }
    \otimes 
    \hackcenter{\begin{tikzpicture}[scale=0.25]
        \draw[thick,orange] (0,0)--(4,0);
        \draw[thick,orange] (0,4)--(4,4);
        \draw (0.5,0)--(0.5,1);
        \draw (1.5,0)--(1.5,1);
        \draw (0.5,1)arc(180:0:0.5);
        \draw (2.5,0)--(2.5,4);
        \draw (3.5,0)--(3.5,4);
    \end{tikzpicture}
    }
    -
     \hackcenter{\begin{tikzpicture}[scale=0.25]
        \draw[thick,orange] (0,0)--(4,0);
        \draw[thick,orange] (0,4)--(4,4);
        \draw (0.5,0)arc(180:0:1.5);
    \end{tikzpicture}
    }
    \otimes
    \hackcenter{\begin{tikzpicture}[scale=0.25]
        \draw[thick,orange] (0,0)--(4,0);
        \draw[thick,orange] (0,4)--(4,4);
        \draw (0.5,0)--(0.5,4);
        \draw (1.5,0)--(1.5,4);
        \draw (2.5,1)arc(180:0:0.5);
        \draw (2.5,0)--(2.5,1);
        \draw (3.5,0)--(3.5,1);
    \end{tikzpicture}
    }
    &\in\qidl,\\  
    {}\hackcenter{\begin{tikzpicture}[scale=0.25]
        \draw[thick,orange] (0,0)--(4,0);
        \draw[thick,orange] (0,4)--(4,4);
        \draw (0.5,0)arc(180:0:1.5);
    \end{tikzpicture}
    }
    \otimes
    \hackcenter{\begin{tikzpicture}[scale=0.25]
        \draw[thick,orange] (0,0)--(4,0);
        \draw[thick,orange] (0,4)--(4,4);
        \draw(0.5,0)--(0.5,1); 
        \draw(2.5,0)--(2.5,1);
        \draw (0.5,1)arc(180:0:1);
        \draw (1.5,0)--(1.5,4);
        \draw (3.5,0)--(3.5,4);
    \end{tikzpicture}
    }
    -
    \hackcenter{\begin{tikzpicture}[scale=0.25]
        \draw[thick,orange] (0,0)--(4,0);
        \draw[thick,orange] (0,4)--(4,4);
        \draw (0.5,0)arc(180:0:1.5);
    \end{tikzpicture}
    }
    \otimes
    \hackcenter{\begin{tikzpicture}[scale=0.25]
        \draw[thick,orange] (0,0)--(4,0);
        \draw[thick,orange] (0,4)--(4,4);
        \draw(0.5,0)--(0.5,4); 
        \draw(2.5,0)--(2.5,4);
        \draw (1.5,1)arc(180:0:1);
        \draw (1.5,0)--(1.5,1);
        \draw (3.5,0)--(3.5,1);
    \end{tikzpicture}
    }
    &\in\qidl,\\
    {}\hackcenter{\begin{tikzpicture}[scale=0.25]
        \draw[thick,orange] (0,0)--(4,0);
        \draw[thick,orange] (0,4)--(4,4);
        \draw (0.5,0)arc(180:0:1.5);
    \end{tikzpicture}
    }
    \otimes
    \hackcenter{\begin{tikzpicture}[scale=0.25]
        \draw[thick,orange] (0,0)--(4,0);
        \draw[thick,orange] (0,4)--(4,4);
        \draw (0.5,0)arc(180:0:1.5);
        \draw (1.5,0)--(1.5,4);
        \draw (2.5,0)--(2.5,4);
    \end{tikzpicture}
    }
    &\in\qidl.
\end{align*}
The other relations involving dots are the obvious ones obtained from these by using the nilBrauer relations (e.g. the dot-crossing commutation relation). Our convention is to prefer that dots in diagrams of $\NB^-$ are situated on the left of caps. For instance, it is clear that
\[\hackcenter{\begin{tikzpicture}[scale=0.375]
        \draw[thick,orange] (0,0)--(4,0);
        \draw[thick,orange] (0,2)--(4,2);
        \draw[thick,orange] (0,4)--(4,4);
        \draw (0.5,0)--(0.5,1);
        \draw (1.5,0)--(1.5,1);
        \draw (0.5,1) arc (180:0:0.5);
        \fill (1,1) ++(135:0.5) circle (5pt);
        \node at (0.2,1) {\tiny $a$};
        \draw (2.5,0)--(2.5,3);
        \draw (3.5,0)--(3.5,3);
        \draw (2.5,3) arc(180:0:0.5);
        \fill (3,3) ++(135:0.5) circle (5pt);
        \node at (2.2,3) {\tiny $b$};
    \end{tikzpicture}
    }
    -
    \hackcenter{\begin{tikzpicture}[scale=0.375]
        \draw[thick,orange] (0,0)--(4,0);
        \draw[thick,orange] (0,2)--(4,2);
        \draw[thick,orange] (0,4)--(4,4);
        \draw (0.5,0)--(0.5,3);
        \draw (1.5,0)--(1.5,3);
        \draw (0.5,3) arc (180:0:0.5);
        \fill (1,3) ++(135:0.5) circle (5pt);
        \node at (0.2,3) {\tiny $a$};
        \draw (2.5,0)--(2.5,1);
        \draw (3.5,0)--(3.5,1);
        \draw (2.5,1) arc(180:0:0.5);
        \fill (3,1) ++(135:0.5) circle (5pt);
        \node at (2.2,1) {\tiny $b$};
    \end{tikzpicture}
    }\in\qidl,\]
where the orange line in the middle of the diagram is introduced as short-hand for the tensor symbol. 

Let us remark that for instance
\[\hackcenter{\begin{tikzpicture}[scale=0.375]
    \draw (0,0) arc(180:0:1);
    \fill (1,0) ++(135:1) circle (5pt);
\end{tikzpicture}}
+
\hackcenter{\begin{tikzpicture}[scale=0.375]
    \draw (0,0) arc(180:0:1);
    \fill (1,0) ++(45:1) circle (5pt);
\end{tikzpicture}}
\not\in\qidl,
\]
because the relation that $\hackcenter{\begin{tikzpicture}[scale=0.375]
    \draw (0,0) arc(180:0:1);
    \fill (1,0) ++(135:1) circle (5pt);
\end{tikzpicture}}=-\hackcenter{\begin{tikzpicture}[scale=0.375]
    \draw (0,0) arc(180:0:1);
    \fill (1,0) ++(45:1) circle (5pt);
\end{tikzpicture}}$ is already true in $\NB^-_1$. 

That the ideal $\qidl$ is quadratic is intuitively obvious and an artifact of the setting of diagrammatic algebras; the specific relations of nilBrauer are not important for this to be true. 
\begin{LEM}
    The ideal $\qidl$ is quadratic.
\end{LEM}
\vspace{-0.5em}
\begin{PRF}
    One way to see this is to appeal to \cite[Chapter 1, Corollary 1.5.3]{polishchuk2005quadratic}:
    \begin{LEM}[\cite{polishchuk2005quadratic}]
        $A$ is quadratic iff $\Ext^i_A{}_j(\BK,\BK)=0$ for $j>i$ and $i=1,2$; $A$ is one-generated iff $\Ext^i_A{}_j(\BK,\BK)=0$ for $j>i$ and $i=1$.
    \end{LEM}
    We have already seen that $\NB^-$ is one-generated, so it suffices to consider $i=2$. Using the resolutions constructed in Section \ref{section:nilcohomo}, one can compute $\Ext_{\NB^-}^2(\BK,\BK)$ as $\ker/\img$, where $\ker$ and $\img$ are the subspaces of $\hom_\BK(I\otimes I,\BK)$ which are
    \[\ker=\{f:f(ab\otimes c)=f(a\otimes bc)\},\quad \img=\{f:f(a\otimes b)=-g(ab)\te{ for }g\in I^*\}.\]
    The claim is that any $f\in\Ker/\Img$ must vanish on tensors $x\otimes y$ which have $\deg_\te{cap}xy\ge 3$, i.e. that any $f\in\Ker$ must be constant on inputs which multiply ($\mu\colon x\otimes y\mapsto xy$) to the same algebra element of cap-degree $\ge3$. This is true for essentially the same reason the Eckmann-Hilton argument works, which is also why this only works for cap-degree $\ge3$. 
    
    One can for example observe that there is a basis of $\NB^-$ in which the caps are visually arranged so that, for example, the highest cap is at the far left and the lowest cap is at the far right. This arrangement is mostly cosmetic; the difference is that dots live on the top-left of caps, so members of this basis have dots which are arranged diagonally from the top-left to the bottom-right. Then, given any $\sum x_i\otimes y_i$ multiplying to a fixed element $z$ of degree $\ge3$, we can rewrite it according to this basis and keep pushing caps across the tensor sign until we arrive at something of form $\sum a_i\otimes z_i$, where $\deg_\te{cap}a_i=1$ and $a_iz_i$ are diagrams in this basis. This process does not change the value of $f\in\Ker$, and the final result is independent of the form of $\sum x_i\otimes y_i$, as long as it multiplies to $z$. 
\end{PRF}
One can now form the quadratic dual of $\NB^-$, which would also be called the Koszul dual if $\NB^-$ were Koszul.

\subsection{The Koszul claim}
At this point, the claim is then that $\NB^-$ is not only quadratic but moreover Koszul. 
Firstly, let us recall (here $A^!$ refers to the quadratic dual of an algebra)
\begin{DEF}
    A graded (locally unital) algebra $A$ with $A_0=\BK$ is called ``Koszul'' if the following equivalent conditions hold:
    \begin{itemize}
        \item $\Ext_{A}^{i}{}_j=0$ for $i\neq j$;
        \item $A$ is one-generated and the algebra $\Ext_A^\blt(\BK,\BK)$ is generated by $\Ext^1_A(\BK,\BK)$;
        \item $A$ is quadratic and $\Ext_A^\blt(\BK,\BK)\cong A^!$;
        \item the algebra $\Ext_A^\blt(\BK,\BK)$, equipped with the cap grading, is one-generated.
    \end{itemize}
\end{DEF}
We then claim that
\begin{THM}
    The quadratic algebra $\NB^-$ is Koszul. 
\end{THM}
\begin{remark}
    We should remark that in the argument spanning the following pages, we make use of two important facts, 4.7.1 and 2.4.1, from \cite{polishchuk2005quadratic}; a review of the proofs of these facts there reveals that the same proofs go through for locally graded-finite-dimensional algebras.
\end{remark}

Let us briefly sketch the plan for the proof. The overall plan is to prove $\NB^-$ is Koszul by putting a filtration (the `dot filtration') on it and showing the associated graded is Koszul; by general Koszul theory this then shows the original $\NB^-$ is Koszul. To show the associated graded is Koszul we will use the characterization of Koszulness using distributive lattices. To do this we will find a distributing basis by dividing-and-conquering according to so-called `intersection/nesting patterns'. And lastly for each fixed intersection/nesting pattern, the problem is reduced to a more classical problem which can be solved using general Koszul theory.


\subsection{Filtrations on quadratic algebras}
First we need to define what it means to be a filtered quadratic algebra.
\begin{DEF}
    A ``graded ordered semigroup'' $\Gamma$ is a (not necessarily commutative) semigroup with unit $e$ equipped with a homomorphism 
    \[g\colon \Gamma\lto\BN\]
    and total orders on the fibers $g^{-1}(n)=\Gamma_n$ such that 
    \[\alpha<\beta\implies\alpha\gamma<\beta\gamma,\ \gamma\alpha<\gamma\beta\uad\foralls\alpha,\beta\in\Gamma_n,\ \gamma\in\Gamma_k\] 
    and $g^{-1}(0)=\{e\}$.
\end{DEF}
Here the symbol $\Gamma$ is not to be confused with the ring of Schur $q$-functions from earlier. We remark that for us it is best to think of the semigroup $\Gamma$ additively; for instance we will write $e$ as 0 and use $+$.
\begin{DEF}
    A ``$\Gamma$-valued filtration'' on a graded algebra $A$ is a collection of subspaces $\F_\alpha A_n\subseteq A_n$ for every $\alpha\in\Gamma_n$ and $n\ge 0$ such that
    \begin{enumerate}
        \item $\alpha\le\beta\implies \F_\alpha A_n\subseteq \F_\beta A_n$,
        \item $\F_{\infty_n} A_n=A_n$ for $\infty_n$ the maximal element of $\Gamma_n$, and
        \item $\F_\alpha A_n\cdot \F_\beta A_m\subseteq \F_{\alpha+\beta}A_{n+m}$.
    \end{enumerate}
    Then the ``associated graded'' of $A$ is
    \[\gr^\F A=\bigoplus_\alpha \F_\alpha A_n/\F_{\alpha'}A_n,\] 
    where $\alpha'\lessdot \alpha$.

    This filtration is said to be `one-generated' if $\F_\alpha A_n=\sum_{i_1+\cdots+ i_n\le \alpha}\F_{i_1}A_1\cdots \F_{i_n}A_1$. 
\end{DEF}
Note that $\F$ is one-generated if and only if $\gr^\F A$ is one-generated. 

Then a fact from general Koszul theory is
\begin{THM}[\cite{polishchuk2005quadratic}4.7.1]
    If $A$ is a quadratic algebra equipped with a one-generated $\Gamma$-valued filtration, then Koszulness of $\gr^\F A$ implies Koszulness of $A$. 
\end{THM}
We will thus give a filtration, called the ``dot filtration'', on $\NB^-$, and prove that its associated graded is Koszul. The benefit of this filtration is that the dot-sliding rules are vastly simplified in the associated graded.




Let $\ol\BN=\BN\cup\{\infty\}$. Consider the (additive) semigroup
\[\Gamma=(0,0)\sqcup \bigsqcup_{\substack{i>0\\j\in\ol\BN }}(i,j),\] 
where $(0,0)$ is the unit and the operation is $(i,j)+(i',j')=(i+i',j+j')$. This comes with a homomorphism
\begin{align*}
    g\colon\Gamma&\lto\BN\\
    (i,j)&\lmto i,
\end{align*}
with each fiber $g^{-1}(i)\cong\ol\BN$ ($i\neq 0$) coming with the usual total order $0\le 1\le 2\le\cdots\le\infty$. Then we endow $\NB^-$ with a $\Gamma$-valued filtration $\F=\F^\te{dot}$ by letting
\[\F_j \NB_i^-\coloneqq \F_{(i,j)}\NB_i^-=\te{diagrams with $\le j$ total dots on $i$ caps}.\] 
Of course, if there are no caps at all, then there is nowhere to put dots. This filtration is moreover one-generated since given any diagram with $i$ caps and $j$ dots (WLOG dots placed on caps), one can write it as the product of $i$ diagrams, each with a single cap and some number of dots on it, such that the total number of dots is $j$. So for instance 
\[\F_0\NB_i^-\subset \F_1\NB_i^-\subset\F_2\NB_i^-\subset\cdots\subset\F_\infty\NB_i^-=\NB_i^-.\]

The associated graded, $\gr^\F\NB^-$, is much simpler. Indeed, modulo diagrams with at most a strictly smaller number of dots, dots slide perfectly across crossings; for example: 
\begin{align*}
    {}\hackcenter{\begin{tikzpicture}[scale=0.375]
        \draw[thick,orange] (0,0)--(4,0);
        \draw[thick,orange] (0,4)--(4,4);
        \draw (2.5,0) arc (0:180:1);
    \draw (1.5,0)--(1.5,2);
    \draw (3.5,0)--(3.5,2);
    \draw (3.5,2) arc (0:180:1);
    \fill (1.5,0.5) circle (5pt);
    \end{tikzpicture}
    }
    -
    \hackcenter{\begin{tikzpicture}[scale=0.375]
    \draw[thick,orange] (0,0) -- (4,0);
    \draw[thick,orange] (0,4) -- (4,4);
    \draw (2.5,0) arc (0:180:1);
    \draw (1.5,0)--(1.5,2);
    \draw (3.5,0)--(3.5,2);
    \draw (3.5,2) arc (0:180:1);
    \fill (1.5,1.5) circle (5pt);
\end{tikzpicture}}
    =
    \hackcenter{\begin{tikzpicture}[scale=0.375]
        \draw[thick,orange] (0,0)--(4,0);
        \draw[thick,orange] (0,4)--(4,4);
        \draw (0.5,0)--(0.5,1);
        \draw (1.5,0)--(1.5,1);
        \draw (0.5,1) arc (180:0:0.5);
        \draw (2.5,0)--(2.5,2);
        \draw (3.5,0)--(3.5,2);
        \draw (2.5,2) arc (180:0:0.5);
    \end{tikzpicture}
    }
    -
    \hackcenter{\begin{tikzpicture}[scale=0.375]
        \draw[thick,orange] (0,0)--(4,0);
        \draw[thick,orange] (0,4)--(4,4);
        \draw (2.5,0) arc (0:180:0.5);
    \draw (0.5,0)--(0.5,2);
    \draw (1.5,2) arc (0:180:0.5);
    \draw (2.5,2) arc (0:-180:0.5);
    \draw (3.5,2) arc (0:180:0.5);
    \draw (3.5,0)--(3.5,2);
    \end{tikzpicture}
    }
    &=
    0\mod \F_0\NB^-_2,\\
    \implies 
    {}\hackcenter{\begin{tikzpicture}[scale=0.375]
        \draw[thick,orange] (0,0)--(4,0);
        \draw[thick,orange] (0,4)--(4,4);
        \draw[thick,orange] (0,2)--(4,2);
        \draw (3.5,0) arc (0:180:1);
        \draw (0.5,0)--(0.5,2.5);
        \draw (2.5,0)--(2.5,2.5);
        \draw (2.5,2.5) arc (0:180:1);
        \fill (2.5,0) ++(135:1) circle (5pt);
    \end{tikzpicture}
    }
    -
    \hackcenter{\begin{tikzpicture}[scale=0.375]
        \draw[thick,orange] (0,0)--(4,0);
        \draw[thick,orange] (0,4)--(4,4);
        \draw[thick,orange] (0,2)--(4,2);
        \draw (2.5,0) arc (0:180:1);
        \draw (1.5,0)--(1.5,2.5);
        \draw (3.5,0)--(3.5,2.5);
        \draw (3.5,2.5) arc (0:180:1);
        \fill (1.5,2.5) circle (5pt);
    \end{tikzpicture}
    }
    &\in\pidl\coloneqq\te{the quadratic ideal of }\gr^\F\NB^-.
\end{align*}
We will now venture forth to show that $\gr^\F\NB^-$ is Koszul. For shorthand let us write
\[B\coloneqq \gr^\F\NB^-,\]
with quadratic kernel $\pidl$ (not a prime ideal; this is not a notational conflict because we will not be considering prime ideals).

\begin{remark}
    This filtration $\F$ can likely be further refined to control the number of dots on each string, e.g. at most $j_1$ dots on string number 1, etc.; however we do not take this approach here.
\end{remark}

\subsection{Distributive lattices and nesting patterns}
Let us first define what a distributive lattice is. Note that they are defined differently in \cite{polishchuk2005quadratic}; we take here as definition an equivalent condition, see Proposition 1.7.1 of \cite{polishchuk2005quadratic}.
\begin{DEF}[\cite{polishchuk2005quadratic}1.7.1]
    For $V$ a vector space, a collection of its subspaces $U_1,\cdotsc,U_n$ is said to be ``distributive'' (or ``generate a distributive lattice'') if there exists a basis of $V$ such that each $U_i$ is the span of a subset of this basis. Such a basis is called a ``distributing basis''.  
\end{DEF}
The significance of this definition is that it can be used to characterize Koszulness of a quadratic algebra.
\begin{PROP}[\cite{polishchuk2005quadratic}2.4.1]
    A quadratic algebra $A$ with quadratic kernel $\qidl$ is Koszul if and only if for every $n\ge 0$ the collection of subspaces
    \[U_i=A_1^{\otimes i-1}\otimes\qidl\otimes A_1^{\otimes n-i-1}\subset A_1^{\otimes n}\] 
    for $1\le i\le n-1$ is distributive. More precisely, the following are equivalent:
    \begin{itemize}
        \item $\Ext^i_A{}_j(\BK,\BK)=0$ for all $i<j\le n$;
        \item the Koszul complex is acyclic in positive internal degrees $\le n$;
        \item the collection $\{U_1,\cdotsc,U_{n-1}\}$ in $A_1^{\otimes n}$ is distributive.
    \end{itemize}
\end{PROP}
Hence we seek a distributing basis of $B_1^{\otimes n}$ for each $n\ge 0$, distributing the collection of spaces $B_1^{\otimes i-1}\otimes\pidl\otimes B_1^{\otimes n-i-1}$. This is a massive space -- there are infinitely many idempotents, and on top of that even for fixed idempotents the spaces are nongraded-infinite-dimensional, but there are some `divide-and-conquer' reductions we can make to make this job more tractable. 

As diagrams with different ending/beginning idempotents are linearly independent, it suffices to give a distributing basis for $e^{\theta-2n}B_1^{\otimes n} e^\theta$ for each $\theta$.
Moreover, if we have a basis for $\theta=2n$, we can add propagating strings in appropriate places to all the basis diagrams to obtain a basis for any desired $\theta$; hence it suffices to give a distributing basis for $e^0 B_1^{\otimes n}e^{2n}$. As diagrams with different numbers of total dots are linearly independent in $B$ (thanks to the associated graded procedure), it also suffices to give a distributing basis for a fixed number of dots. By adding dots in appropriate places to basis diagrams, we can get any number of total dots we desire; hence it suffices to give a distributing basis for $e^0\mr B_1^{\otimes n}e^{2n}$, where $\mr B$ is the algebra obtained from $B$ by modding out by all diagrams with a positive number of dots:
\[\mr B\coloneqq \bigoplus_i\F_0 \NB_i^-.\]
Most importantly, note that the spaces $B_1^{\otimes i-1}\otimes\pidl\otimes B_1^{\otimes n-i-1}$ respect these reductions, in the sense that each relation in $\pidl$ is homogeneous with respect to the beginning/ending idempotents, as well as with respect to the number of dots. 

Things are already looking much more tractable, but we can make things even easier by looking at nesting patterns.
\begin{DEF}
    The ``nesting pattern'' of a diagram is given by looking at its `shadow'. More precisely, ignore dots (if they exist), color each cap, and project the diagram to the horizontal axis; the result is a 1-dimensional topological subspace of the horizontal axis, consisting of segments of different colors, possibly overlapping, with their intersection inheriting both colors. Let these segments have a `thickness' coming from the number of colors they carry. Then forget the color, and the resulting 1-dimensional subspace with thickness, considered up to isotopy within the 1-dimensional parent space, is the nesting pattern of this diagram.
\end{DEF}  
Here's a quick diagram depicting this:
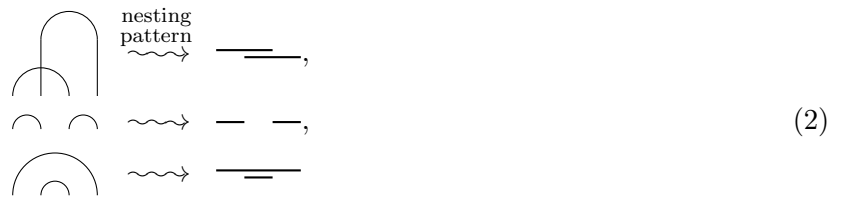
\begin{align*}
    {}\hackcenter{\begin{tikzpicture}[scale=0.375]
        \draw (0.5,0) arc (180:0:1);
        \draw (1.5,0)--(1.5,2);
        \draw (3.5,0)--(3.5,2);
        \draw (1.5,2) arc (180:0:1);
    \end{tikzpicture}}
    &\ \overset{\substack{\te{nesting} \\\te{pattern}}}{\longsquigglyrightarrow}\ 
    \hackcenter{\begin{tikzpicture}[scale=0.375]
        \draw[thick] (1.5,0)--(3.5,0);
        \draw[thick] (0.5,0.25)--(2.5,0.25);
    \end{tikzpicture}},\\\tag{2}\label{diag:pattern}
    {}\hackcenter{\begin{tikzpicture}[scale=0.375]
        \draw (0.5,0) arc (180:0:0.5);
        \draw (2.5,0) arc (180:0:0.5);
    \end{tikzpicture}}
    &\ \longsquigglyrightarrow
    \ \hackcenter{\begin{tikzpicture}[scale=0.375]
        \draw[thick] (0.5,0)--(1.5,0);
        \draw[thick] (2.5,0)--(3.5,0);
    \end{tikzpicture}},\\
    {}\hackcenter{\begin{tikzpicture}[scale=0.375]
        \draw (0.5,0) arc (180:0:1.5);
        \draw (1.5,0) arc (180:0:0.5);
    \end{tikzpicture}}
    &\ \longsquigglyrightarrow\ 
    \hackcenter{\begin{tikzpicture}[scale=0.375]
        \draw[thick] (0.5,0.25)--(3.5,0.25);
        \draw[thick] (1.5,0)--(2.5,0);
    \end{tikzpicture}}
\end{align*}
Here to make the picture clearer we have staggered the shadows of the caps in drawing the nesting pattern; it does not matter which one goes on top; all that matters is where they overlap, if they do. It is clear that diagrams with different nesting patterns are linearly independent.
Also importantly, all relations in $\pidl$ are homogeneous with respect to nesting patterns. 
Note that this is not true without passing from $\NB^-$ to $B=\gr^\F\NB^-$! Hence it suffices to provide a distributing basis of the subspace of $e^0\mr B_1^{\otimes n}e^{2n}$ consisting of a fixed nesting pattern $\te P$; denote this subspace $e^0\mr B_1^{\otimes n}e^{2n}(\te P)$. 

So let us fix a nesting pattern $\te P$. Also fix a numbering of each segment, in such a way that a segment labeled $j$ nested inside a segment labeled $i$ has $i<j$:
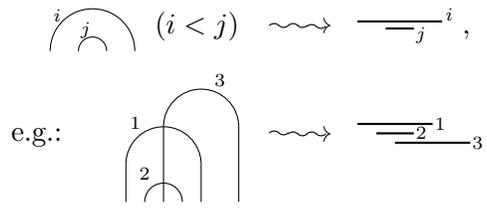
\begin{align*}
    {}\hackcenter{\begin{tikzpicture}[scale=0.375]
        \draw (0.5,0) arc (180:0:1.5);
        \draw (1.5,0) arc (180:0:0.5);
        \node at (0.75,1.25) {\tiny $i$};
        \node at (1.75,0.75) {\tiny $j$};
    \end{tikzpicture}}
    \ \ (i<j)&\ \longsquigglyrightarrow\ 
    \hackcenter{\begin{tikzpicture}[scale=0.375]
        \draw[thick] (0.5,0.25)--(3.5,0.25);
        \draw[thick] (1.5,0)--(2.5,0);
        \node at (3.75,0.5){\tiny $i$};
        \node at (2.75,-0.25){\tiny $j$};
    \end{tikzpicture}},\\
    \te{e.g.:}\qquad 
    \hackcenter{\begin{tikzpicture}[scale=0.5]
        \draw (0.5,0)--(0.5,1);
        \draw (2.5,0)--(2.5,1);
        \draw (0.5,1) arc (180:0:1);
        \draw (1,0) arc (180:0:0.5);
        \draw (1.5,0)--(1.5,2);
        \draw (3.5,0)--(3.5,2);
        \draw (1.5,2) arc (180:0:1);
        \node at (0.75,2.1) {\tiny $1$};
        \node at (1,0.75) {\tiny $2$};
        \node at (3,3.25) {\tiny $3$};
    \end{tikzpicture}}
    &\ \longsquigglyrightarrow\ 
    \hackcenter{\begin{tikzpicture}[scale=0.5]
        \draw[thick] (0.5,0.5)--(2.5,0.5);
        \draw[thick] (1,0.25)--(2,0.25);
        \draw[thick] (1.5,0)--(3.5,0);
        \node at (2.7,0.5) {\tiny $1$};
        \node at (2.2, 0.25) {\tiny $2$};
        \node at (3.7,0) {\tiny $3$};
    \end{tikzpicture}}
\end{align*}
For any two segments, diagram \ref{diag:pattern} exhausts the possibilities of their nesting patterns.
Let each segment $i$ have a formal algebraic symbol $x_i$ attached to it. We think of $x_i$ as `the cap for segment $i$'. Given a monomial in the $x_i$'s, we read it as right to left corresponding to bottom to top in nilBrauer diagrams. For example, if one has a pattern $\te P$ of a segment labeled 2 sitting inside a segment labeled 1, then $x_1x_2$ (cf. $x_2x_1$) corresponds to the following element of $e^0\mr B_1^{\otimes 2}e^4$:
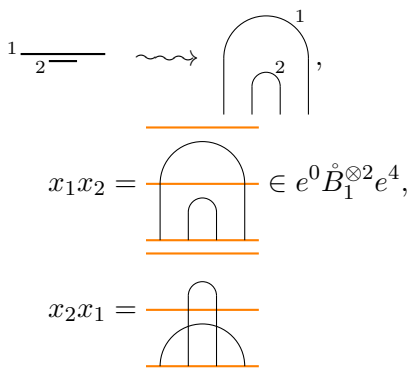
\begin{align*}
    {}\hackcenter{\begin{tikzpicture}[scale=0.375]
        \draw[thick] (0.5,0.25)--(3.5,0.25);
        \draw[thick] (1.5,0)--(2.5,0);
        \node at (0.2,0.45) {\tiny $1$};
        \node at (1.2,-0.2) {\tiny $2$};
    \end{tikzpicture}}
    &\ \longsquigglyrightarrow\ 
    \hackcenter{\begin{tikzpicture}[scale=0.375]
        \draw (0.5,0)--(0.5,2);
        \draw (3.5,0)--(3.5,2);
        \draw (0.5,2) arc (180:0:1.5);
        \draw (1.5,0)--(1.5,1);
        \draw (2.5,0)--(2.5,1);
        \draw (1.5,1) arc (180:0:0.5);
        \node at (3.2,3.5) {\tiny $1$};
        \node at (2.5,1.7) {\tiny $2$};
    \end{tikzpicture}},\\
    x_1x_2&=
    \hackcenter{\begin{tikzpicture}[scale=0.375]
        \draw[thick,orange] (0,0)--(4,0);
        \draw[thick,orange] (0,4)--(4,4);
        \draw[thick,orange] (0,2)--(4,2);
        \draw (0.5,0)--(0.5,2);
        \draw (3.5,0)--(3.5,2);
        \draw (0.5,2) arc (180:0:1.5);
        \draw (1.5,0)--(1.5,1);
        \draw (2.5,0)--(2.5,1);
        \draw (1.5,1) arc (180:0:0.5);
    \end{tikzpicture}}
    \in e^0 \mr B_1^{\otimes 2}e^4,\\
    x_2x_1&=
    \hackcenter{\begin{tikzpicture}[scale=0.375]
        \draw[thick,orange] (0,0)--(4,0);
        \draw[thick,orange] (0,4)--(4,4);
        \draw[thick,orange] (0,2)--(4,2);
        \draw (0.5,0) arc (180:0:1.5);
        \draw (1.5,0)--(1.5,2.5);
        \draw (2.5,0)--(2.5,2.5);
        \draw (1.5,2.5) arc (180:0:0.5);
    \end{tikzpicture}}
\end{align*}
Note that $x_2x_1$ is then an element in $\pidl$. Similarly, if the pattern was two segments neither of which sits entirely inside the other, then $x_1x_2-x_2x_1$ is an element in $\pidl$. Here's a slightly more complicated example:
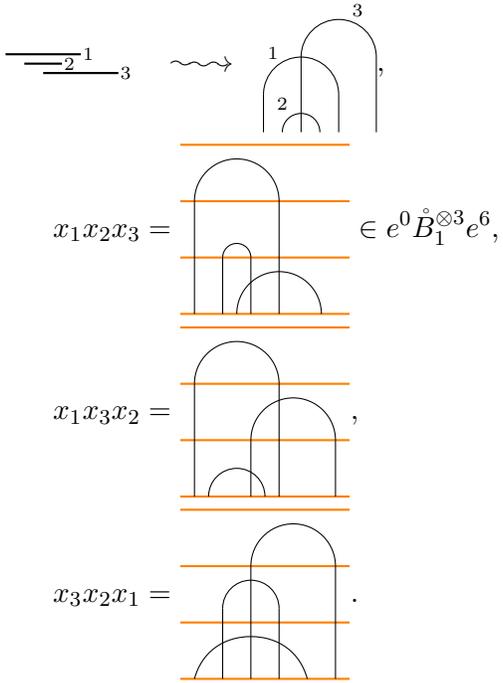
\begin{align*}
    {}\hackcenter{\begin{tikzpicture}[scale=0.5]
        \draw[thick] (0.5,0.5)--(2.5,0.5);
        \draw[thick] (1,0.25)--(2,0.25);
        \draw[thick] (1.5,0)--(3.5,0);
        \node at (2.7,0.5) {\tiny $1$};
        \node at (2.2, 0.25) {\tiny $2$};
        \node at (3.7,0) {\tiny $3$};
    \end{tikzpicture}}
    &\ \longsquigglyrightarrow\ 
    \hackcenter{\begin{tikzpicture}[scale=0.5]
        \draw (0.5,0)--(0.5,1);
        \draw (2.5,0)--(2.5,1);
        \draw (0.5,1) arc (180:0:1);
        \draw (1,0) arc (180:0:0.5);
        \draw (1.5,0)--(1.5,2);
        \draw (3.5,0)--(3.5,2);
        \draw (1.5,2) arc (180:0:1);
        \node at (0.75,2.1) {\tiny $1$};
        \node at (1,0.75) {\tiny $2$};
        \node at (3,3.25) {\tiny $3$};
    \end{tikzpicture}},\\
    x_1x_2x_3&=
    \hackcenter{\begin{tikzpicture}[scale=0.375]
        \draw[thick,orange] (0,0)--(6,0);
        \draw[thick,orange] (0,6)--(6,6);
        \draw[thick,orange] (0,2)--(6,2);
        \draw[thick,orange] (0,4)--(6,4);
        \draw (0.5,0)--(0.5,4);
        \draw (3.5,0)--(3.5,4);
        \draw (0.5,4) arc (180:0:1.5);
        \draw (1.5,0)--(1.5,2);
        \draw (2.5,0)--(2.5,2);
        \draw (1.5,2) arc (180:0:0.5);
        \draw (2,0) arc (180:0:1.5);
    \end{tikzpicture}}
    \in e^0 \mr B_1^{\otimes 3}e^6,\\
    x_1x_3x_2&=
    \hackcenter{\begin{tikzpicture}[scale=0.375]
        \draw[thick,orange] (0,0)--(6,0);
        \draw[thick,orange] (0,6)--(6,6);
        \draw[thick,orange] (0,2)--(6,2);
        \draw[thick,orange] (0,4)--(6,4);
        \draw (0.5,0)--(0.5,4);
        \draw (3.5,0)--(3.5,4);
        \draw (0.5,4) arc (180:0:1.5);
        \draw (1,0) arc (180:0:1);
        \draw (2.5,0)--(2.5,2);
        \draw (5.5,0)--(5.5,2);
        \draw (2.5,2) arc (180:0:1.5);
    \end{tikzpicture}},\\
    x_3x_2x_1&=
    \hackcenter{\begin{tikzpicture}[scale=0.375]
        \draw[thick,orange] (0,0)--(6,0);
        \draw[thick,orange] (0,6)--(6,6);
        \draw[thick,orange] (0,2)--(6,2);
        \draw[thick,orange] (0,4)--(6,4);
        \draw (2.5,0)--(2.5,4);
        \draw (5.5,0)--(5.5,4);
        \draw (2.5,4) arc (180:0:1.5);
        \draw (0.5,0) to[out=75,in=180] (2.5,1.5);
        \draw (2.5,1.5) to[out=0,in=105] (4.5,0);
        \draw (1.5,0)--(1.5,2.5);
        \draw (3.5,0)--(3.5,2.5);
        \draw (1.5,2.5) arc (180:0:1);
    \end{tikzpicture}}.
\end{align*}

More generally, given any two segments $i,j$ in $\te P$, corresponding to two caps $x_i,x_j$, either
\[x_ix_j-x_jx_i\in\pidl\quad\te{or}\quad x_jx_i\in\pidl,\uad i<j;\] 
the latter case corresponds to if segment $j$ lies inside segment $i$, and the former is otherwise. Since the relation of `lying inside of' is transitive, namely $k\subset j,\ j\subset i\implies k\subset i$, we have
\[x_kx_j\in\pidl,\ x_jx_i\in\pidl\implies x_kx_i\in\pidl.\] 
Lastly, we only want each cap to appear once, so we must impose for each $i$
\[x_i^2\in\pidl.\]

We want a distributing basis for $e^0\mr B_1^{\otimes n}e^{2n}(\te P)$. There will be $n$ segments in $\te P$, corresponding to caps $x_1,\cdotsc,x_n$. Consider the algebra $A=A(\te P)$ over $\BC$ with generators $x_1,\cdotsc,x_n$ subject to the quadratic relations
\begin{enumerate}
    \item $x_i^2=0$ for each $i$,
    \item either $x_ix_j=x_jx_i$ or $x_jx_i=0$ for every $i<j$, depending on $\te P$; this will automatically be in such a way that $x_kx_j=0$ and $x_jx_i=0\implies x_kx_i=0$;
\end{enumerate}
clearly, if this were $n$-Koszul\footnote{We will note use the definition of $n$-Koszul in this paper, but for completeness let us recall it is defined as follows (see \cite[Chapter 2, Theorem 4.1]{polishchuk2005quadratic}): a quadratic algebra is $n$-Koszul if $\Ext_A^i{}_j(\BK,\BK)=0$ for all $i<j\le n$, or equivalently the collection $\{A_1^{\otimes i-1}\otimes\qidl\otimes A_1^{\otimes n-i-1}\}_{1\le i\le n-1}$ is distributive in $A_1^{\otimes n}$.}, then one has a basis of $A_1^{\otimes n}$ distributing the collection $\{A_1^{\otimes i-1}\otimes\qidl\otimes A_1^{\otimes n-i-1}\}_i$, and this same basis will be a distributing basis for $e^0\mr B_1^{\otimes n}e^{2n}(\te P)$ under the obvious isomorphism. So we have reduced what we want to a (semi?)commutative algebra question. Thankfully it turns out such algebras $A(\te P)$ are not only $n$-Koszul but Koszul.

\subsection{Reduction to Koszul theory}
At this point we can cite a theorem from general Koszul theory.
\begin{THM}[Froberg, \cite{polishchuk2005quadratic}4.7]
    Consider the quadratic algebra $A$ with generators $x_1,\cdotsc,x_n$ subject to relations from the following list of (distinct) possibilities for each $i\neq j$:
    \begin{enumerate}
        \item $x_ix_j=c_{ij}x_jx_i$ for $c_{ij}\in\bk^\times$;
        \item either $x_ix_j=0$ or $x_jx_i=0$ (but not both);
        \item $x_ix_j=x_jx_i=0$;
        \item no relation;
        \item possibly $x_i^2=0$ or not for each $i$.
    \end{enumerate}
    Let $\theta_i=x_i^*$ be the generators of $A^!$. Then $A$ is Koszul if and only if for any two nonzero monomials $x^\alpha\in A$ and $\theta^\beta\in A^!$, one can find $i$ and $\gamma$ such that either one of the following two is true (where $c\in \bk^\times$):
    \begin{itemize}
        \item $x^\alpha=cx^\gamma x_i$ and $\theta_i \theta^\beta\neq 0$,
        \item $x^\alpha x_i\neq 0$ and $c\theta_i\theta^\gamma=\theta^\beta$. 
    \end{itemize}
\end{THM}
Using this, one can prove that for instance quotients of skew polynomial rings by monomial quadratic ideals are Koszul. The fact we will need is
\begin{PROP}
    Let $A$ be the algebra over $\BC$ with generators $x_1,\cdotsc,x_n$ subject to the relations
    \begin{enumerate}
        \item $x_i^2=0$ for each $i$,
        \item either $x_ix_j=x_jx_i$ or $x_jx_i=0$ for every $i<j$, in such a way that $x_kx_j=0$ and $x_jx_i=0\implies x_kx_i=0$.
    \end{enumerate}
    Then $A$ is Koszul. 
\end{PROP}
This proposition, in turn, is a special case of the following algebraic fact:
\begin{PROP}
    Consider an algebra $A$ which is generated by $x_1,\cdotsc,x_n$ subject to the relations
    \begin{enumerate}
        \item $x_i^2=0$ for each $i$;
        \item for each $i<j$ either:
        \begin{enumerate}
            \item $x_ix_j=c_{ij}x_jx_i$ for $c_{ij}\in\bk^\times$, or
            \item $x_jx_i=0$,
        \end{enumerate}
        such that $x_kx_j=0,\ x_jx_i=0\implies x_kx_i=0$. 
    \end{enumerate}
    Then $A$ is Koszul. 
\end{PROP}
\vspace{-0.5em}
\begin{PRF}
    This is in essence a combinatorial argument. Recall the quadratic dual is $A^!=\mfrac{\otimes^\blt A_1^*}{\qidl^\perp}$; as such, inside $A^!$ the relations satisfied by $\theta_i$ are precisely
    \begin{enumerate}
        \item either $\theta_i\theta_j=-c_{ij}^{-1}\theta_j\theta_i$ or $\theta_i\theta_j=0$ for every $i<j$, in such a way that $\theta_i\theta_j=0$ and $\theta_j\theta_k=0\implies \theta_i\theta_k=0$.
    \end{enumerate}
    The relations $x_i^2=0$ turn into no restrictions on nilpotence of $\theta_i$. Let us state for the record that the only nonzero possibilities for two variables to not commute are $\theta_{\te{big}}\theta_{\te{small}}$ and $x_{\te{small}}x_{\te{big}}$.

    We utilize the characterization of Froberg. Consider any two nonzero monomials $x^\alpha=(\cdots)x_i$ and $\theta^\beta=\theta_j(\cdots)$. Write $(\cdots)x_i$ on the left and $\theta_j(\cdots)$ on the right:
    \[x^\alpha=\underbrace{(\cdots)}_{x^{\alpha'}}x_i\qquad \theta_j\underbrace{(\cdots)}_{\theta^{\beta'}}=\theta^\beta;\]
    we wish to show we can push one of $x_i$ or $\theta_j$ to the other side and still have two nonzero monomials.
    
    A quick lemma:
    \begin{LEM}
        In the above setup, in $(x^{\alpha'})x_i\neq 0$, if $x_j$ is such that $x_jx_i=0$, then $x_j$ cannot be a factor in $x^{\alpha'}$. Dually, in $\theta_j(\theta^{\beta'})\neq 0$, if $\theta_k$ is such that $\theta_j\theta_k=0$, then $\theta_k$ cannot be a factor in $\theta^{\beta'}$. 
    \end{LEM}
    Let us quickly see this. We argue only for the former case, as the other argument is dual. Suppose there was a $x_j$ in $x^{\alpha'}$; let us try to commute it to the right. If it can be commuted rightwards (using relations of form $x_ax_b=c_{ab}x_bx_a$) to meet $x_i$, then one has $x_jx_i$ as a subexpression (up to a nonzero constant) and this monomial must be zero. If it cannot be commuted to meet $x_i$, it must be that it were stuck somewhere along the way as $x_jx_k$ for $j<k$; but for this to be the case, it must be that $x_kx_j=0$, which together with $x_jx_i=0$ implies $x_kx_i=0$. Now $x_k$ is closer to $x_i$ than $x_j$, and we can again try to commute it towards $x_i$. By induction on the distance from $x_i$, this shows the monomial must be zero. 
    
    Now we return to trying to push one of $x_i$ or $\theta_j$ to the other side. We break into cases.

    \textit{Case 1.} Suppose $\theta_i\theta_j=0$, $i<j$, which is to say that $x_jx_i=0$. Then split $\theta_j$ off from $\theta^\beta$ and move it to $x^\alpha$; note $x_ix_j\neq 0$. In order to show $x^\alpha x_j=x^{\alpha'} x_ix_j\neq 0$, we would need to show the introduction of $x_j$ does not allow for any new cancellation, namely that (a) it cannot be the case that there is a $x_j$ in $x^{\alpha'}$ that can be commuted to meet the new $x_j$ and (b) it cannot be the case that there is a $x_k$ such that $x_kx_j=0$ in $x^{\alpha'}$ that can be commuted to meet the new $x_j$. The former possibility is impossible due to the lemma above. The latter case is also not possible because $x_kx_j=0$ and $x_jx_i=0$ gives $x_kx_i=0$, and again by the lemma such a $x_k$ cannot be a factor in $x^{\alpha'}$. 

    \textit{Case 2.} Suppose $\theta_i\theta_j=-c_{ij}^{-1}\theta_j\theta_i$, and that there is no $\theta_k$ in $\theta^{\beta'}$ such that $\theta_i\theta_k=0$. Then we can simply split $x_i$ from $x^\alpha$ and move it to $\theta^\beta$, and there is no possibility of cancellation. 

    \textit{Case 3.} Suppose $\theta_i\theta_j=-c_{ij}^{-1}\theta_j\theta_i$, but there is a $\theta_k$ in $\theta^{\beta'}$ such that $\theta_i\theta_k=0$. Again, by trying to commute $\theta_k$ towards $\theta_j$ (using relations like $\theta_a\theta_b=-c_{ab}^{-1}\theta_b\theta_a$) and inducting on the distance from $\theta_j$, we can get to the point $\theta^\beta=c\theta_j\theta_{k'}(\cdots)$, where $k'$ might be different from $k$ but still satisfies $\theta_i\theta_{k'}=0$ (by using the procedure as in the proof of the lemma). But $\theta_j,\theta_{k'}$ must commute (up to $-c_{k'j}^{-1}$), for otherwise we would have $\theta_{k'}\theta_j=0$, which would imply $\theta_i\theta_j=0$, contradiction. Hence we can write $\theta^\beta=c'\theta_{k'}\theta_j(\cdots)$, where $\theta_i\theta_{k'}=0$; this reduces to Case 1.

    These three cases exhaust all possibilities, and we are done.
\end{PRF}

\begin{remark}
It is remarkable that this fact about the nilBrauer algebra ends up boiling down to a (non?)commutative algebra fact! Maybe worth remarking is that this is almost not true. Indeed, according to this proposition, \[\mfrac{\BC\wan{x_1,x_2,x_3}}{x_i^2=0,\ x_1x_3=x_3x_1,\ x_2x_3=x_3x_2,\ x_2x_1=0}\] is Koszul, but the algebra \[\mfrac{\BC\wan{x_1,x_2,x_3}}{x_1x_3=x_3x_1,\ x_1x_2=x_2x_1=0}\] is not, as remarked by \cite{polishchuk2005quadratic} on page 91.

Also, note that the above proof still works even if you weaken the $x_i^2=0$ for each $i$ to $x_i^2=0$ for some $i$. We don't need this though.
\end{remark}


\section{Standard Ext Groups, Part 3: Kostantness and BGG (Main Results)}\label{section:bgg}
\subsection{Kostantness}\label{subsec:kostant}
Having established Koszulness, we can now show
\begin{THM}
    At $t=0$, the Ext groups $\Ext^\blt_\NB(\NB^{\ge\theta}e^\theta,L_0)$ are concentrated in one degree (namely $n=\theta/2$); i.e., $L_0$ is a Kostant module.
\end{THM}
\vspace{-0.5em}
\begin{PRF}
    Indeed, we can compute these Ext groups $\rhom_\NB(\NB^{\ge\theta}e^\theta,\BC)=\rhom_{A^-}(\bk e^\theta,\BC)$ using the resolution of $\bk e^\theta$ given above. So
    \begin{align*}
        \Ext^\blt_\NB(\NB^{\ge\theta}e^\theta,\BC e^0)&=\Ext^\blt_{A^-}(\bk e^\theta,\BC e^0)\\
        &=H^\blt\hom_{A^-}(A^-\otimes I^{\otimes\blt}e^\theta,\BC e^0)\\
        &=H^\blt \hom_\BK(I^{\otimes\blt}e^\theta,\BC e^0)\\
        &=H^\blt\hom_\BC(e^0I^{\otimes\blt}e^\theta,\BC e^0)\\
        &=H^{-\blt} e^0I^{\otimes\blt}e^\theta.
    \end{align*}
    Here the complex $e^0 I^{\otimes \blt}e^\theta$ has differential coming from the original complex, namely $\pd_k\colon a_1\otimes\cdots \otimes a_k\lmto \sum_{i=1}^{k-1}(-1)^i a_1\otimes\cdots a_ia_{i+1}\cdots\otimes a_k$. Note that for instance $e^0 I^{\otimes 0}e^\theta =e^0\BK e^\theta=0$ if $\theta\neq 0$ (the $\theta=0$ case is trivial anyway since $\Delta_0$ is projective).

    However, as we know $\NB^-$ is Koszul, we immediately know that $\rhom_{\NB^-}(\BK,\BK)$ is diagonal, meaning that it is only nonzero when the homological degree matches with the cap degree. As we've mentioned earlier, $\rhom_{\NB^-}(\BK,\BK)$ has a left and right $\BK$-action corresponding to picking out the weights in the first and second entries, and by picking out the weight space $\rhom_{\NB^-}(\bk e^\theta,\bk e^0)$, we see that the cap degree is concentrated in $n=\theta/2$, so that necessarily the homological degree must also be concentrated in $n$.
\end{PRF}

In fact, we can compute the dimensions of these Ext groups. Before doing so, some setup: Brundan had introduced $T_{f,n}(q)$ as the $q$-generating function of the number of tethered chord diagrams with $f$ free chords and $n$ tethered chords. Rather than use this notation, we will let
\[Y_{i,j}\te{ be the generating function for dotless $\te Y$ diagrams from $j$ to $i$}.\]
The relation between the two is $T_{f,n}=Y_{n,2f+n}$ and $Y_{i,j}=T_{\frac{j-i}{2},i}$. BWW showed that these generating functions satisfy $T_{f,n}=T_{f,n-1}+q^n[n+1]T_{f-1,n+1}$, though we will not appeal to this.
\begin{remark}
    Another way to deduce the following dimension claim is to use the linear independence of the symbols $[\ol\Delta_{\lbd}]$ in the Grothendieck group of the nilBrauer modules, in conjunction with the BWW character formula 
    \[[L_0]=\sum_{n=0}^\infty(-1)^n\frac{q^{-n}}{(1-q^{-4})(1-q^{-8})\cdots(1-q^{-4n})}[\ol\Delta_{2n}].\]
    Then, by the spectral sequence in Theorem \ref{thm:nbspecseq} applied to $L_0$, the alternating sums of the characters of the terms of the first page must be the same as the character of $L_0$, but by Kostantness this alternating sum is simply the alternating sum of the standard Ext groups across $n=\theta/2$, so that we know
    \[[L_0]=\sum_{n=0}^\infty (-1)^n[\jota^{2n}_!\Ext^n_\NB(\Delta(2n),L_0)^*]=\sum_{n=0}^\infty (-1)^n\frac{1}{[2n]!}\dim(\Ext^n_\NB(\Delta(2n),L_0)^*)[\ol\Delta_{2n}],\] 
    where we recall that $\Ext^n_\NB(\Delta(2n),M)^*$ is a $\NB^{2n}$-module and note $\jota^{2n}_! L_{2n}(2n)=\ol\Delta_{2n}$ and $\dim L_{2n}(2n)=[2n]!$. Comparison with the BWW formula gives an alternative way to deduce the following.
\end{remark}
\ \vspace{-0.5em}
\begin{THM}
    Let $t=0$ and $\theta=2n$. Then the dimension of the Ext groups are exactly what one expects for the BGG resolution to hold, namely
    \[\dim_q \Ext^n_\NB(\NB^{\ge\theta}e^\theta,L_0)^*=\frac{q^{-n}[2n]!}{(1-q^{-4})\cdots(1-q^{-4n})}=\frac{(1+q^2+q^4)\cdots (1+\cdots+q^{4n-4})}{(1-q^{-2})^n}.\]
    For convenience in the future we will denote $m_\theta(q)=\frac{q^{-n}}{(1-q^{-4})\cdots(1-q^{-4n})}$.
\end{THM}
\vspace{-0.5em}
\begin{PRF}
    As the complex (here the heart symbol indicates homological degree 0)
    \[0\lto e^0 I^{\otimes n}e^\theta\lto\cdots\lto e^0Ie^\theta\lto \overset{\heart}{e^0\BK e^\theta}\lto 0\]
    is quasiisomorphic to $\Ext^\blt_\NB(\NB^{\ge\theta}e^\theta,L_0)^*=\Ext^n_\NB(\NB^{\ge\theta}e^\theta,L_0)^*[n]$, one has
    \[\sum_{k=0}^n(-1)^k \dim_q e^0I^{\otimes k}e^\theta=(-1)^n\dim_q \Ext^n_\NB(\NB^{\ge\theta}e^\theta,L_0)^*,\] 
    where
    \[\dim_q e^{0}I^{\otimes k}e^{2n}=\frac{1}{(1-q^{-2})^n}\sum_{1\le i_1<\cdots<i_{k-1}\le n-1} Y_{0,2i_1}\cdots Y_{2i_{k-1},2n};\] 
    the factor in front is for the dots (recall propagating strands cannot have dots for $\te Y$).

    So the claim is then saying that
    \[\sum_{k=0}^n(-1)^k\frac{1}{(1-q^{-2})^n}\sum_{1\le i_1<\cdots<i_{k-1}\le n-1} Y_{0,2i_1}\cdots Y_{2i_{k-1},2n}=(-1)^n\frac{q^{-n}[2n]!}{(1-q^{-4})\cdots(1-q^{-4n})}.\]
    A straight-forward computation shows that
    \begin{align*}
        (1-q^{-2})^n\frac{[2n]!}{(1-q^{-4})\cdots(1-q^{-4n})}&=\frac{[2n]!}{(1+q^{-2})(1+q^{-2}+q^{-4}+q^{-6})\cdots(1+\cdots+q^{-4n+2})}\\
        &=\frac{q^{-n}[2n]!}{(1-q^{-4})\cdots(1-q^{-4n})}\cdot \frac{q^{1+2+\cdots+(2n-1)}}{q^{2+6+\cdots+(4n-2)}}\cdot \frac{q^{1+2+\cdots+(2n-1)}}{q^{4+8+\cdots+(4n-4)}}\\
        &=q^n\frac{(1+q^2)(1+q^2+q^4)\cdots (1+\cdots+q^{4n-2})}{(1+q^2)(1+\cdots+q^6)\cdots(1+\cdots+q^{4n-2})}\\
        &=q^n(1+q^2+q^4)(1+q^2+q^4+q^6+q^8)\cdots(1+\cdots+q^{4n-4}),
    \end{align*}
    so that the above claim is equivalent to saying (note the $k=0$ term is trivial)
    \[(1+q^2+q^4)(1+q^2+q^4+q^6+q^8)\cdots(1+\cdots+q^{4n-4})=\sum_{k=0}^n\sum_{1\le i_1<\cdots<i_{k-1}\le n-1}(-1)^{n-k} Y_{0,2i_1}\cdots Y_{2i_{k-1},2n}.\]


    This is morally a PIE (principle of inclusion-exclusion) argument. On the LHS is the generating function for number of ways to go from $2n$ strings to 0 strings in $n$ steps, where each step has exactly one cap whose left-most vertex is at the far left (or alternatively whose right-most vertex is at the far right). On the RHS one starts at $k=n$ with the product of all the ways to go from $2n$ strings to 0 strings in $n$ steps except with no requirement on where the endpoints of the cap are, and each subsequent term consists of ways of combining multiple steps into one. To make this identity easier to show, we can apply a reduction using the idea of nesting patterns from earlier: it is clear that for each nesting pattern, there is exactly one cap diagram counted by the LHS, namely the one for which the caps are arranged like a staircase from left to right. Then it suffices to show this for a fixed nesting pattern, as the terms on the RHS can also be organized according to nesting patterns. WLOG it suffices to show this for the nesting pattern of $n$ disjoint segments, where the desired identity is
    \[1=\sum_{k=1}^n\sum_{1\le i_1<\cdots< i_{k-1}\le n-1}(-1)^{n-k}\binom{i_1}{0}\cdots\binom{n}{i_{k-1}}.\] 
    We can now symbolically manipulate
    \begin{align*}
        \sum_{k=1}^n\sum_{1\le i_1<\cdots< i_{k-1}\le n-1}(-1)^{n-k}\binom{i_1}{0}\cdots\binom{n}{i_{k-1}}&=\sum_{k=1}^n\sum_{\substack{j_1+\cdots+j_k=n\\j_i>0}}(-1)^{n-k}\binom{n}{j_1,\cdotsc,j_k}\\
      \te{(PIE)}\implies  &=\sum_{k=1}^n(-1)^{n-k} \big(  k^n-\binom{k}{1}(k-1)^n+\binom{k}{2}(k-2)^n-\cdots \big)\\
      &=\sum_{k=1}^n (-1)^{n-k}\sum_{i=0}^k(-1)^{k-i}\binom{k}{i}i^n\\
      &=\sum_{k=1}^n\sum_{i=0}^\infty(-1)^{n-i}\binom{k}{i}i^n\\
      &=\sum_{i=0}^\infty (-1)^{n-i}i^n\sum_{k=1}^n\binom{k}{i}\\
      &=\sum_{i=0}^n (-1)^{n-i}\binom{n+1}{i+1}i^n.
    \end{align*}
    So we wish to show $\sum_{i=0}^n (-1)^{n-i}\binom{n+1}{i+1}i^n=1$. In fact, we claim that
    \[\sum_{i=0}^n (-1)^{n-i}\binom{n+1}{i+1}i^k=(-1)^{n-k}\qquad\foralls k\le n.\]
    This follows by applying the following lemma to $P(i)=(i-1)^k$ and using $n+1$ as the top entry in the binomial coefficient. Indeed, then one has 
    \begin{align*}
        0&=\sum_{i+1=0}^{n+1}(-1)^{i+1}\binom{n+1}{i+1}P(i+1)\\
        &=\sum_{i=-1}^n (-1)^{n-i}\binom{n+1}{i+1}i^k\\
        &=(-1)^{n+1}\binom{n+1}{0}(-1)^k+\sum_{i=0}^n (-1)^{n-i}\binom{n+1}{i+1}i^k\\
        \implies (-1)^{n-k}&=\sum_{i=0}^n (-1)^{n-i}\binom{n+1}{i+1}i^k,
    \end{align*}
    as claimed.
\end{PRF}
\ \vspace{-1em}
\begin{LEM}
    For any polynomial $\deg P<n$, one has
    \[\sum_{i=0}^n(-1)^i\binom{n}{i}P(i)=0.\] 
\end{LEM}
\vspace{-0.5em}
\begin{PRF}
    One can prove this using generating function theory. Indeed, consider 
    \[\sum_{i=0}^n\binom{n}{i}P(i)x^i;\] 
    we wish to compute this for $x=-1$. By general theory, this is
    \[\sum_{i=0}^n\binom{n}{i}P(i)x^i=P(x\pd)(1+x)^n.\] 
    The claim is that
    \[(x\pd)^k(1+x)^n\rv_{x=-1}=0\qquad\foralls k<n,\]
    so that by linearity we are done. But this is true by induction: it is evident that each $(x\pd)^{\te{some power}}(1+x)^n$ will look like linear combinations of terms of form $x^i(1+x)^j$. Note
    \[(x\pd)(x^i(1+x)^j)=ix^i(1+x)^j+jx^{i+1}(1+x)^{j-1}.\] 
    In particular, if $x\pd$ is applied less than $n$ times, the exponent on $(1+x)$ will never be zero, so when $x$ is set to $-1$ all these terms vanish. This concludes the proof.

    Another method, much cleaner, is due to Matthew Hase-Liu, using the philosophy of finite differences. Using that the binomial coefficients $\binom{x}{i}$ for $0\le i\le n$ form a basis for the space of polynomials of degree at most $n$, it suffices to show the claim for $P(i)=\binom{i}{k}$. However, 
    \[\sum_{i}(-1)^i\binom{n}{i}\binom{i}{k}=\sum_{i} (-1)^i\binom{n-k}{n-i}\binom{n}{k}=\binom{n}{k}\sum_{i} (-1)^i\binom{n-k}{n-i}=0,\] 
    so we are done.
\end{PRF}

Recall that the standard Ext groups $\Ext_\NB^n(\NB^{\ge\theta}e^\theta,L_0)^*$ act as multiplicity spaces for the BGG resolution in a sense. Unfortunately, these standard Ext groups are not semisimple as $\NB^\theta$-modules; the next best thing we can do to study the structure of these Ext groups, and by extension the structure of the terms appearing in the BGG resolution, is to write down a Jordan-Holder filtration for it, with descriptions of its quotients. We give a Loewy filtration (in the sense that quotients are semisimple) below. 
\begin{THM}
    At $t=0$, letting $\theta=2n$, the $\NB^\theta$-module $\Ext_\NB^n(\NB^{\ge\theta}e^\theta,L_0)^*$ is isomorphic to a quotient of the polynomial ring
    \[\Ext_\NB^n(\NB^{\ge\theta}e^\theta,L_0)^*\cong \mfrac{\BC[X_1,\cdotsc,X_\theta]}{\wan{\tsl p_1,\tsl p_3,\cdotsc,\tsl p_{2n-1}}},\]
    where $\tsl p_i$ is the $i$-th power sum in the dots $X$, and has a filtration $\Ext_\NB^n(\NB^{\ge\theta}e^\theta,L_0)^*=\F^0_{\Ext}\supset \F^1_{\Ext}\supset\cdots $ such that the quotients are
    \[\gr^k\Ext_\NB^n(\NB^{\ge\theta}e^\theta,L_0)^*=L_\theta(\theta)\otimes_\BC q^{-n}\BC[\tsl p_2,\tsl p_4,\cdotsc,\tsl p_{2n}]_{\deg_{\sym}=k}=L_\theta(\theta)\otimes_\BC q^{-n}E_{\theta,k},\]
    where $E_{\theta,k}\coloneqq \BC[\tsl p_2,\cdotsc,\tsl p_{2n}]_{\deg_{\sym}=k}$ denotes the subspace of $\BC[\tsl p_2,\cdotsc,\tsl p_{2n}]$ consisting of symmetric-degree-$k$ elements, where the symmetric degree is defined by $\deg_\sym \tsl p_i=1$ (contrast this with the quantum degree, which has $\deg \tsl p_i=2i$). 
\end{THM}
\vspace{-0.5em}
\begin{PRF}
    To prove this, we utilize the nilBrauer action on the Ext group as realized by $H^n_\flat(\NB^-:L_0)^*$. Note that as $L_0=\BC e^0$, we can ignore the $\Hom_\BK(\sq,L_0)^*$ and calculate the Ext groups directly as the cohomology of $e^0(I^\flat)^{\otimes_{A^\nat}\blt}e^\theta$. Under this identification, the Cartan nilHeckes act on the right, by attaching diagrams below. 

    We have already computed the dimension of $\Ext^{n,*}=\Ext_\NB^n(\NB^{\ge\theta}e^\theta,L_0)^*$ above. Consider the element of lowest quantum degree, corresponding to the term $q^{4+8+\cdots+4n-4}=q^{4\binom{n}{2}}$ in the quantum dimension. This element is constructed as (we remind the reader the orange lines in the middle stand for tensor symbols)
    \[v^+_\theta=v^+=
    \hackcenter
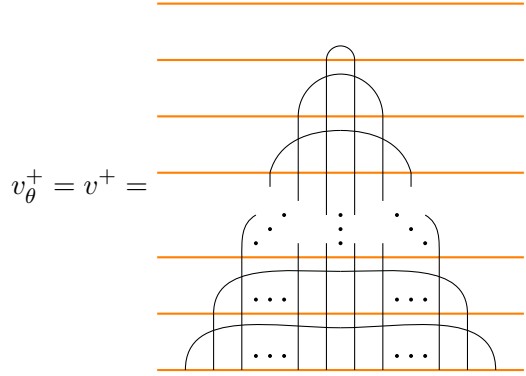
{\begin{tikzpicture}[scale=0.375]
        \draw[thick,orange] (0,0)--(13,0);
        \draw[thick,orange] (0,2)--(13,2);
        \draw[thick,orange] (0,4)--(13,4);
        \draw[thick,orange] (0,7)--(13,7);
        \draw[thick,orange] (0,9)--(13,9);
        \draw[thick,orange] (0,11)--(13,11);
        \draw[thick,orange] (0,13)--(13,13);
        \draw (6,0)--(6,4.5);
        \draw (6,5.5)--(6,11);
        \draw (7,0)--(7,4.5);
        \draw (7,5.5)--(7,11);
        \draw (6,11) arc (180:0:0.5);
        \draw (5,0)--(5,4.5);
        \draw (8,0)--(8,4.5);
        \draw (5,6)--(5,9);
        \draw (8,6)--(8,9);
        \draw (5,9) arc (180:0:1.5);
        \draw (4,7)--(4,6.5);
        \draw (9,7)--(9,6.5);
        \draw (4,7) to[out=75,in=180] (6.5,8.5);
        \draw (6.5,8.5) to[out=0,in=105] (9,7);
        \draw (3,0)--(3,4);
        \draw (10,0)--(10,4);
        \draw (3,4) to[out=90,in=210] (3.5,5.5);
        \draw (10,4) to[out=90,in=-30] (9.5,5.5);
        \draw (2,0)--(2,2);
        \draw (11,0)--(11,2);
        \draw (2,2) to[out=90,in=180](6.5,3.5);
        \draw (6.5,3.5) to[out=0,in=90] (11,2);
        \draw (1,0) to[out=90,in=180] (6.5,1.5);
        \draw (6.5,1.5) to[out=0,in=90] (12,0);
        \fill (3.5,0.5) circle (2pt);
        \fill (4,0.5) circle (2pt);
        \fill (4.5,0.5) circle (2pt);
        \fill (8.5,0.5) circle (2pt);
        \fill (9,0.5) circle (2pt);
        \fill (9.5,0.5) circle (2pt);
        \fill (3.5,2.5) circle (2pt);
        \fill (4,2.5) circle (2pt);
        \fill (4.5,2.5) circle (2pt);
        \fill (8.5,2.5) circle (2pt);
        \fill (9,2.5) circle (2pt);
        \fill (9.5,2.5) circle (2pt);
        \fill (6.5,4.5) circle (2pt);
        \fill (6.5,5) circle (2pt);
        \fill (6.5,5.5) circle (2pt);
        \fill (3.5,4.5) circle (2pt);
        \fill (4,5) circle (2pt);
        \fill (4.5,5.5) circle (2pt);
        \fill (8.5,5.5) circle (2pt);
        \fill (9,5) circle (2pt);
        \fill (9.5,4.5) circle (2pt);
    \end{tikzpicture}}
    \]
    As an example, in the $n=2,3$ cases $v^+$ is
    \[v^+_4=\hackcenter
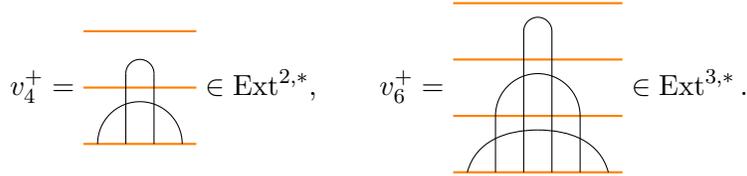
{\begin{tikzpicture}[scale=0.375]
        \draw[thick,orange] (0,0)--(4,0);
        \draw[thick,orange] (0,4)--(4,4);
        \draw[thick,orange] (0,2)--(4,2);
        \draw (0.5,0) arc (180:0:1.5);
        \draw (1.5,0)--(1.5,2.5);
        \draw (2.5,0)--(2.5,2.5);
        \draw (1.5,2.5) arc (180:0:0.5);
    \end{tikzpicture}}\in\Ext^{2,*},\qquad
    v^+_6=\hackcenter{\begin{tikzpicture}[scale=0.375]
        \draw[thick,orange] (0,0)--(6,0);
        \draw[thick,orange] (0,6)--(6,6);
        \draw[thick,orange] (0,2)--(6,2);
        \draw[thick,orange] (0,4)--(6,4);
        \draw (2.5,0)--(2.5,5);
        \draw (3.5,0)--(3.5,5);
        \draw (2.5,5) arc(180:0:0.5);
        \draw (1.5,0)--(1.5,2);
        \draw (4.5,0)--(4.5,2);
        \draw (1.5,2) arc (180:0:1.5);
        \draw (0.5,0) to[out=75,in=180] (3,1.5);
        \draw (3,1.5) to[out=0,in=105] (5.5,0);
    \end{tikzpicture}}\in \Ext^{3,*}.
    \]
    Note that, considered as an element in the complex computing $H^\blt_\flat$ (so that the tensor is over $\NB^\nat$ and therefore allows crossings to pass through), this is the same as the non-intersecting nest of loops on top of two copies of the longest permutation on $n$ strands, denoted $w_0 * w_0$:
    \[
        v^+=\hackcenter
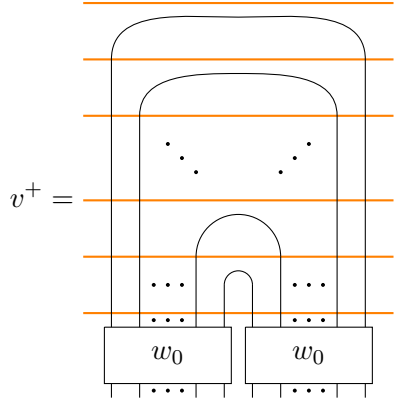
{\begin{tikzpicture}[scale=0.375]
            \draw[thick,orange] (0,0)--(11,0);
            \draw[thick,orange] (0,2)--(11,2);
            \draw[thick,orange] (0,4)--(11,4);
            \draw[thick,orange] (0,7)--(11,7);
            \draw[thick,orange] (0,9)--(11,9);
            \draw[thick,orange] (0,11)--(11,11);
            \draw (5,0)--(5,1);
            \draw (6,0)--(6,1);
            \draw (5,1) arc (180:0:0.5);
            \draw (4,0)--(4,2);
            \draw (7,0)--(7,2);
            \draw (4,2) arc (180:0:1.5);
            \draw (2,0)--(2,7) to[out=90,in=180] (5.5,8.5) to[out=0,in=90] (9,7)--(9,0);
            \draw (1,0)--(1,9) to[out=90,in=180] (5.5,10.5) to[out=0,in=90] (10,9)--(10,0);
            \fill (2.5,1) circle (2pt);
            \fill (3,1) circle (2pt);
            \fill (3.5,1) circle (2pt);
            \fill (7.5,1) circle (2pt);
            \fill (8,1) circle (2pt);
            \fill (8.5,1) circle (2pt);
            \fill (3,6) circle (2pt);
            \fill (3.5,5.5) circle (2pt);
            \fill (4,5) circle (2pt);
            \fill (7,5) circle (2pt);
            \fill (7.5,5.5) circle (2pt);
            \fill (8,6) circle (2pt);
            \fill (2.5,-0.25) circle (2pt);
            \fill (3,-0.25) circle (2pt);
            \fill (3.5,-0.25) circle (2pt);
            \fill (2.5,-2.75) circle (2pt);
            \fill (3,-2.75) circle (2pt);
            \fill (3.5,-2.75) circle (2pt);
            \fill (7.5,-0.25) circle (2pt);
            \fill (8,-0.25) circle (2pt);
            \fill (8.5,-0.25) circle (2pt);
            \fill (7.5,-2.75) circle (2pt);
            \fill (8,-2.75) circle (2pt);
            \fill (8.5,-2.75) circle (2pt);
            \draw (1,0)--(1,-0.5);
            \draw (2,0)--(2,-0.5);
            \draw (4,0)--(4,-0.5);
            \draw (5,0)--(5,-0.5);
            \draw (1,-2.5)--(1,-3);
            \draw (2,-2.5)--(2,-3);
            \draw (4,-2.5)--(4,-3);
            \draw (5,-2.5)--(5,-3);
            \draw (0.75,-0.5)--(5.25,-0.5)--(5.25,-2.5)--(0.75,-2.5)--(0.75,-0.5);
            \node at (3,-1.5) {$w_0$};
            \draw (6,0)--(6,-0.5);
            \draw (7,0)--(7,-0.5);
            \draw (9,0)--(9,-0.5);
            \draw (10,0)--(10,-0.5);
            \draw (6,-2.5)--(6,-3);
            \draw (7,-2.5)--(7,-3);
            \draw (9,-2.5)--(9,-3);
            \draw (10,-2.5)--(10,-3);
            \draw (5.75,-0.5)--(10.25,-0.5)--(10.25,-2.5)--(5.75,-2.5)--(5.75,-0.5);
            \node at (8,-1.5) {$w_0$};
        \end{tikzpicture}}
    \]
    Again as an example
    \[v^+_4=\hackcenter
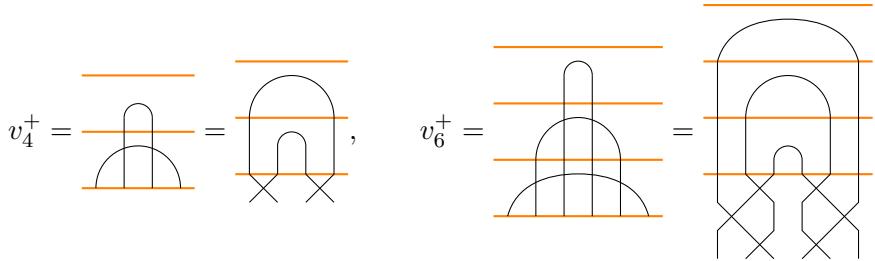
{\begin{tikzpicture}[scale=0.375]
        \draw[thick,orange] (0,0)--(4,0);
        \draw[thick,orange] (0,4)--(4,4);
        \draw[thick,orange] (0,2)--(4,2);
        \draw (0.5,0) arc (180:0:1.5);
        \draw (1.5,0)--(1.5,2.5);
        \draw (2.5,0)--(2.5,2.5);
        \draw (1.5,2.5) arc (180:0:0.5);
    \end{tikzpicture}}
    =
    \hackcenter{\begin{tikzpicture}[scale=0.375]
        \draw[thick,orange] (0,0)--(4,0);
        \draw[thick,orange] (0,4)--(4,4);
        \draw[thick,orange] (0,2)--(4,2);
        \draw (0.5,0)--(0.5,2);
        \draw (3.5,0)--(3.5,2);
        \draw (0.5,2) arc (180:0:1.5);
        \draw (1.5,0)--(1.5,1);
        \draw (2.5,0)--(2.5,1);
        \draw (1.5,1) arc (180:0:0.5);
        \draw (0.5,0)--(1.5,-1);
        \draw (0.5,-1)--(1.5,0);
        \draw (2.5,0)--(3.5,-1);
        \draw (3.5,0)--(2.5,-1);
    \end{tikzpicture}},\qquad
    v^+_6=\hackcenter{\begin{tikzpicture}[scale=0.375]
        \draw[thick,orange] (0,0)--(6,0);
        \draw[thick,orange] (0,6)--(6,6);
        \draw[thick,orange] (0,2)--(6,2);
        \draw[thick,orange] (0,4)--(6,4);
        \draw (2.5,0)--(2.5,5);
        \draw (3.5,0)--(3.5,5);
        \draw (2.5,5) arc(180:0:0.5);
        \draw (1.5,0)--(1.5,2);
        \draw (4.5,0)--(4.5,2);
        \draw (1.5,2) arc (180:0:1.5);
        \draw (0.5,0) to[out=75,in=180] (3,1.5);
        \draw (3,1.5) to[out=0,in=105] (5.5,0);
    \end{tikzpicture}}
    =
    \hackcenter{\begin{tikzpicture}[scale=0.375]
        \draw[thick,orange] (0,0)--(6,0);
        \draw[thick,orange] (0,6)--(6,6);
        \draw[thick,orange] (0,2)--(6,2);
        \draw[thick,orange] (0,4)--(6,4);
        \draw (0.5,0)--(0.5,4);
        \draw (5.5,0)--(5.5,4);
        \draw (0.5,4) to[out=75,in=180] (3,5.5);
        \draw (3,5.5) to[out=0,in=105] (5.5,4);
        \draw (1.5,0)--(1.5,2);
        \draw (4.5,0)--(4.5,2);
        \draw (1.5,2) arc (180:0:1.5);
        \draw (2.5,0)--(2.5,0.5);
        \draw (3.5,0)--(3.5,0.5);
        \draw (2.5,0.5) arc (180:0:0.5);
        \draw (0.5,0)--(0.5,-1)--(2.5,-3);
        \draw (1.5,0)--(2.5,-1)--(2.5,-2)--(1.5,-3);
        \draw (2.5,0)--(0.5,-2)--(0.5,-3);
        \draw (3.5,0)--(5.5,-2)--(5.5,-3);
        \draw (4.5,0)--(3.5,-1)--(3.5,-2)--(4.5,-3);
        \draw (5.5,0)--(5.5,-1)--(3.5,-3);
    \end{tikzpicture}}
    \]
    Indeed, since the multiplication of any two adjacent tensor factors is zero, this term is clearly in the kernel of the complex, and by using the $I^{\otimes\blt}$ incarnation of the nilcohomology one sees this vector is nonzero in the Ext group (indeed there is no $n+1$ term and so there is nothing to mod out by). By using the nested loops on longest permutations description, it is clear that any crossing will kill this vector, i.e. $(v^+)\cdot \sigma=0$ (note the right action notation), which is why we called it $v^+$. Recall from \cite[Theorem 3.8/Corollary 3.9]{brundan2023nilbrauer2} that the Schur $\tsl q$-functions $\tsl q_{\te{odd}}$ and the odd power sums in dots $\tsl p_\te{odd}$ are in the center of nilBrauer. Since the odd power sum in zero variables is zero, this immediately shows that $(v^+)\cdot \tsl p_\te{odd}=0$. While the odd $\tsl q$-functions do not become zero at the zeroth idempotent, by commuting one of them to the top of a diagram, one can witness the entire thing as an image under the differential; hence any diagram with at least one $\tsl q$-function is zero in the cohomology. Hence $\tsl q_\te{odd}$ also act by zero on $v^+$, $(v^+)\cdot \tsl q_\te{odd}=0$. 

    Remark: as, under the usual grading conventions, the highest weight vector of $L_\theta(\theta)$ should have degree $-\binom{\theta}{2}$, the vector $v^+$ which is of degree $-4\binom{n}{2}=-\binom{\theta}{2}+n$ would generate a simple (after passing to the head) shifted upwards by $n$. 

    Now consider the vectors $(v^+)\cdot \tsl p_{2i}$ (note the right action notation). By writing $v^+$ as the nest of loops on $w_0*w_0$, we can commute $\tsl p_{2i}$ (being in the center of nilHecke) past the crossings, and it is visibly evident that $\tsl p_\te{even}$ acting on the nest of loops is nonzero, e.g.
    \[
        \hackcenter{\begin{tikzpicture}[scale=0.375]
        \draw[thick,orange] (0,0)--(4,0);
        \draw[thick,orange] (0,4)--(4,4);
        \draw[thick,orange] (0,2)--(4,2);
        \draw (0.5,0)--(0.5,2);
        \draw (3.5,0)--(3.5,2);
        \draw (0.5,2) arc (180:0:1.5);
        \draw (1.5,0)--(1.5,1);
        \draw (2.5,0)--(2.5,1);
        \draw (1.5,1) arc (180:0:0.5);
        \draw (0.5,0)--(0.5,-1);
        \fill (0.5,-0.3) circle (5pt);
        \fill (0.5,-0.7) circle (5pt);
    \end{tikzpicture}}
    +
    \hackcenter{\begin{tikzpicture}[scale=0.375]
        \draw[thick,orange] (0,0)--(4,0);
        \draw[thick,orange] (0,4)--(4,4);
        \draw[thick,orange] (0,2)--(4,2);
        \draw (0.5,0)--(0.5,2);
        \draw (3.5,0)--(3.5,2);
        \draw (0.5,2) arc (180:0:1.5);
        \draw (1.5,0)--(1.5,1);
        \draw (2.5,0)--(2.5,1);
        \draw (1.5,1) arc (180:0:0.5);
        \draw (1.5,0)--(1.5,-1);
        \fill (1.5,-0.3) circle (5pt);
        \fill (1.5,-0.7) circle (5pt);
    \end{tikzpicture}}
    +
    \hackcenter{\begin{tikzpicture}[scale=0.375]
        \draw[thick,orange] (0,0)--(4,0);
        \draw[thick,orange] (0,4)--(4,4);
        \draw[thick,orange] (0,2)--(4,2);
        \draw (0.5,0)--(0.5,2);
        \draw (3.5,0)--(3.5,2);
        \draw (0.5,2) arc (180:0:1.5);
        \draw (1.5,0)--(1.5,1);
        \draw (2.5,0)--(2.5,1);
        \draw (1.5,1) arc (180:0:0.5);
        \draw (2.5,0)--(2.5,-1);
        \fill (2.5,-0.3) circle (5pt);
        \fill (2.5,-0.7) circle (5pt);
    \end{tikzpicture}}
    +
    \hackcenter{\begin{tikzpicture}[scale=0.375]
        \draw[thick,orange] (0,0)--(4,0);
        \draw[thick,orange] (0,4)--(4,4);
        \draw[thick,orange] (0,2)--(4,2);
        \draw (0.5,0)--(0.5,2);
        \draw (3.5,0)--(3.5,2);
        \draw (0.5,2) arc (180:0:1.5);
        \draw (1.5,0)--(1.5,1);
        \draw (2.5,0)--(2.5,1);
        \draw (1.5,1) arc (180:0:0.5);
        \draw (3.5,0)--(3.5,-1);
        \fill (3.5,-0.3) circle (5pt);
        \fill (3.5,-0.7) circle (5pt);
    \end{tikzpicture}}
    =
    2\hackcenter{\begin{tikzpicture}[scale=0.375]
        \draw[thick,orange] (0,0)--(4,0);
        \draw[thick,orange] (0,4)--(4,4);
        \draw[thick,orange] (0,2)--(4,2);
        \draw (0.5,0)--(0.5,2);
        \draw (3.5,0)--(3.5,2);
        \draw (0.5,2) arc (180:0:1.5);
        \fill (0.5,0.7) circle (5pt);
        \fill (0.5,1.1) circle (5pt);
        \draw (1.5,0)--(1.5,1);
        \draw (2.5,0)--(2.5,1);
        \draw (1.5,1) arc (180:0:0.5);
    \end{tikzpicture}}
    +
    2\hackcenter{\begin{tikzpicture}[scale=0.375]
        \draw[thick,orange] (0,0)--(4,0);
        \draw[thick,orange] (0,4)--(4,4);
        \draw[thick,orange] (0,2)--(4,2);
        \draw (0.5,0)--(0.5,2);
        \draw (3.5,0)--(3.5,2);
        \draw (0.5,2) arc (180:0:1.5);
        \draw (1.5,0)--(1.5,1.25);
        \draw (2.5,0)--(2.5,1.25);
        \draw (1.5,1.25) arc (180:0:0.5);
        \fill (1.5,0.6) circle (5pt);
        \fill (1.5,1) circle (5pt);
    \end{tikzpicture}}
    \neq 0.
    \]
    Hence the $(v^+)\cdot \tsl p_\te{even}$ are nonzero. It is also clear that they are mutually distinct (for example by looking at quantum degree). Moreover, since any further crossing also commutes with $\tsl p_\te{even}$, we have $(v^+\tsl p_\te{even})\cdot \sigma=0$. In this way we obtain many mutually distinct `highest weight' (in the sense of being killed by all crossings) vectors. More generally 
    \[\{(v^+)\cdot P:P\in\BC[\tsl p_2,\tsl p_4,\cdotsc,\tsl p_\theta]\}\] 
    are all mutually distinct nonzero highest weight vectors. Note that a basis for this lattice is given by
    \[\{(v^+)\cdot \tsl p_\te{even partition}\},\] 
    where an even partition means a partition whose parts are all even.

    Now consider
    \[\mfrac{(\jidl^\theta)^k}{(\jidl^\theta)^{k+1}}\otimes_{\NB^\theta}\Ext^{n,*}.\] 
    This is a semisimple $\NB^\theta$-representation. Under the identification above, inside of it are distinct highest weight vectors $(v^+)\cdot \tsl p_\lbd$ for even partitions $\lbd$ of length $k$. By the linear independence of these vectors and the semisimplicity of the action over the nilHecke $\NB^\theta$, we know that, when $\NB^\theta$ is allowed to act on them, they each generate a copy of $L_\theta(\theta)$ (appropriately shifted). This gives $\dim \BC[\tsl p_2,\cdotsc, \tsl p_{2n}]_{\deg_\sym=k}$ many copies of $L_\theta(\theta)$, each shifted by $q^{-n}$, in the $k$-th layer of $\Ext^{n,*}$. Together, across all $k$, these would have dimension equal to $q^{-n}l_\theta(\theta)\dim \BC[\tsl p_2,\cdotsc, \tsl p_{2n}]$, which is already the full dimension of $\Ext^{n,*}$. Hence this must be everything. We immediately get the desired filtration with quotients of the claimed dimension.
    
    And in particular, $v^+$ generates all of $\Ext^{n,*}$ via the action of dots, and since the dimension of $\BC[X_1,\cdotsc,X_\theta]/\wan{\tsl p_\te{odd}}$ is already the correct dimension of $\Ext^{n,*}$, we conclude it must be the case that 
    \[\Ext^{n,*}\cong \BC[X_1,\cdotsc,X_\theta]/\wan{\tsl p_\te{odd}},\] 
    where the isomorphism is via the action on $v^+$.

    Let us quickly check that $\BC[X_1,\cdotsc,X_\theta]/\wan{\tsl p_\te{odd}}$ is indeed the correct dimension. By PIE, this dimension is
    \begin{align*}
        \dim \BC[X_1,\cdotsc,X_\theta]/\wan{\tsl p_\te{odd}}&= \frac{1}{(1-q^{-2})^\theta}\Big(1-q^{-2}-q^{-6}-\cdots -q^{-2(\theta-1)}\\
        &\hspace{10em} +q^{-(2+6)}+q^{-(2+10)}+\cdots+q^{-(2+2(\theta-1))}\\
        &\hspace{10em} +q^{-(6+10)}+q^{-(6+14)}+\cdots+q^{-(6+2(\theta-1))}\\
        &\hspace{11em} \vdots\\
        &\hspace{10em}+\cdots+q^{-(2(\theta-3)+2(\theta-1))}\\
        &\hspace{12em}-\cdots\Big),
    \end{align*}
    which can be collapsed to
    \begin{align*}
        \dim \BC[X_1,\cdotsc,X_\theta]/\wan{\tsl p_\te{odd}}&=\frac{(1-q^{-2})(1-q^{-6})\cdots (1-q^{-2(\theta-1)})}{(1-q^{-2})^\theta}\\
        &= \frac{(1-q^{-2})(1-q^{-6})\cdots (1-q^{-4n+2})}{(1-q^{-2})^{2n}}\\
        &=\frac{(1+q^{-2}+q^{-4})(1+q^{-2}+\cdots+q^{-8})\cdots(1+q^{-2}+\cdots+q^{-4n+4})}{(1-q^{-2})^n}\\
        &=q^{(-4)+(-8)+\cdots+(-4n+4)}\frac{(1+q^2+q^4)\cdots (1+\cdots+q^{4n-4})}{(1-q^{-2})^n}\\
        &=q^{-4\binom{n}{2}}\frac{(1+q^2+q^4)\cdots (1+\cdots+q^{4n-4})}{(1-q^{-2})^n}.
    \end{align*}
    Note that the $q^{-4\binom{n}{2}}$ cancels out with the degree of the vector $v^+$, and the remaining part is recognizable from the earlier theorem statement on the dimension of $\Ext^{n,*}$. This concludes the proof.
\end{PRF}



\subsection{BGG}
The results of the previous subsection, together with Theorem \ref{thm:nbspecseq}, show that there is a resolution of $L_0$ at $t=0$ by extensions of proper standard modules. This is the `BGG resolution' of the one-dimensional simple. 

Indeed, apply \ref{thm:nbspecseq} to $M=L_0$, which lies in the even block, namely $\vtheta=0$. Then the spectral sequence in Theorem \ref{thm:nbspecseq} has first page terms $E_1^{p,q}$ nonzero only when $-(p+q)=\frac{1}{2}(\vtheta-2p)=-p$, thanks to our result on the concentration of standard Ext groups. In other words $q=0$ on $E_1$ and $p\le 0$, i.e. the first page is concentrated on the negative $x$-axis:
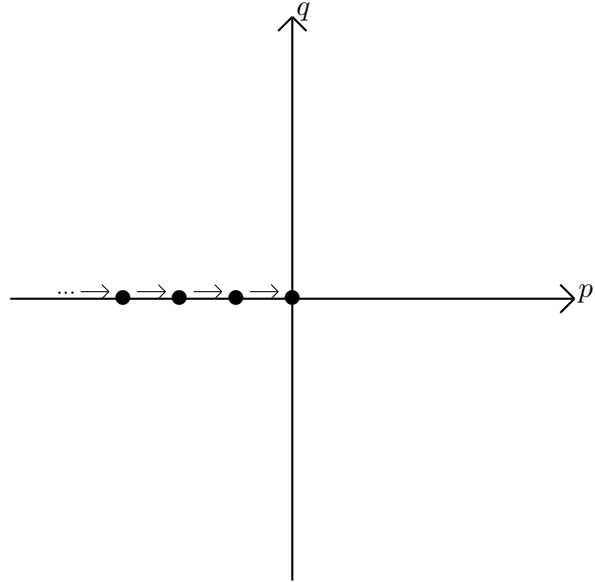
\begin{center}
    \begin{tikzpicture}[scale=0.375]
        \draw[thick] (-10,0)--(10,0);
        \draw[thick] (10,0)--(9.5,0.5);
        \draw[thick] (10,0)--(9.5,-0.5);
        \node at (10.4,0.2) {$p$};
        \draw[thick] (0,-10)--(0,10);
        \draw[thick] (0,10)--(0.5,9.5);
        \draw[thick] (0,10)--(-0.5,9.5);
        \node at (0.4,10.2) {$q$};
        \node (A) at (0,0) {\Large $\blt$};
        \node (B) at (-2,0) {\Large $\blt$};
        \node (C) at (-4,0) {\Large $\blt$};
        \node (D) at (-6,0) {\Large $\blt$};
        \node (E) at (-8,0) {${}^{\cdots}$};
        \draw (-7.5,0.25)--(-6.5,0.25);
        \draw (-6.5,0.25)--(-6.75,0);
        \draw (-6.5,0.25)--(-6.75,0.5);
        \draw (-5.5,0.25)--(-4.5,0.25);
        \draw (-4.5,0.25)--(-4.75,0);
        \draw (-4.5,0.25)--(-4.75,0.5);
        \draw (-3.5,0.25)--(-2.5,0.25);
        \draw (-2.5,0.25)--(-2.75,0);
        \draw (-2.5,0.25)--(-2.75,0.5);
        \draw (-1.5,0.25)--(-0.5,0.25);
        \draw (-0.5,0.25)--(-0.75,0);
        \draw (-0.5,0.25)--(-0.75,0.5);
    \end{tikzpicture}
\end{center}
Then, as $L_0$ is concentrated in homological degree 0 (so that $E_\infty^{p,q}$ must be $L_0$ sitting at the origin), it is clear that the spectral sequence must converge after one step, i.e. $E_1^{\blt,0}$ is a resolution of $L_0$. Moreover, recalling that $\Delta(2n)=\Delta_{2n}^{\ol{l_{2n}(2n)}}=\NB^{\ge 2n}e^{2n}$ has an obvious right $\NB^{2n}$-action, we have
\[E_1^{-n,0}=\Delta_{2n}^{\ol{l_{2n}(2n)}}\otimes_{\NB^{2n}}\Ext^n_\NB(\Delta_{2n}^{\ol{l_{2n}(2n)}},L_0)^*,\] 
where the structure of the standard Ext group has already been determined. Hence each term of this resolution is an extension of proper standard modules whose structure is completely understood. So in conclusion we have shown:
\begin{THM}[BGG Resolution]
    At parameter $t=0$, the 1-dimensional simple $L_0$ has a BGG resolution
    \[\cdots\lto C^{-n}_\te{BGG}(L_0)\lto C^{-(n-1)}_\te{BGG}(L_0)\lto\cdots\lto C^0_\te{BGG}(L_0)\lto L_0\lto 0\] 
    where the terms have character
    \[\chi(C^{-n}_\te{BGG}(L_0))=\frac{q^{-n}}{(1-q^{-4})(1-q^{-8})\cdots(1-q^{-4n})}\chi(\ol\Delta_{2n})\] 
    and admit filtrations $C^{-n}_\te{BGG}(L_0)=\F^0_\te{BGG}\supset \F^1_\te{BGG}\supset\cdots$ such that
    \[\gr^kC^{-n}_\te{BGG}(L_0)=\ol\Delta_{2n}\otimes_\BC q^{-n}\BC[\tsl p_2,\tsl p_4,\cdotsc, \tsl p_{2n}]_{\deg_\sym=k},\] 
    where $\deg_\sym \tsl p_i=1$. 
\end{THM}
This provides a categorification of the character formula Equation (\ref{eqn:bww}) for $L_0$, thorough in the sense that a description of each $\NB$-module in the resolution is given. 

Note that Theorem \ref{thm:nbspecseq} provides a categorification via spectral sequences of the character formula (\ref{eqn:bww}) for all simples, even though in general we do not have a description of the module structure of each term. Indeed, plugging in $M=L_{2n+\vtheta}\in\D^-(\Mod^\vtheta\,\NB)$ for $\vtheta=0,1$, one obtains
\begin{THM}[BGG Spectral Sequence]
    There is a spectral sequence
    \[E^{p,q}_1=\Delta(\vtheta-2p)\otimes_{\NB^{\vtheta-2p}}\Ext_\NB^{-(p+q)}(\Delta(\vtheta-2p),L_{2n+\vtheta})^*\specseqimplies E^{p,q}_\infty=\begin{cases} L_{2n+\vtheta} & p=q=0\\ 0 &\te{else}\end{cases};\]
    the terms $E^{p,q}_1$ of the first page are extensions of proper standard modules $\ol\Delta_{\vtheta-2p}$, and by the fact that simples and proper standards form bases of the Grothendieck group, this must categorify the character formula (\ref{eqn:bww}) in all cases.
\end{THM}

\section{Proper Standard Ext Groups}

In the previous section we have managed to compute the nilHecke structure of the standard Ext groups (or equivalently nilcohomology) because in this specific case of the nilBrauer we had a very good guess for the `highest weight vector'. In general this may be difficult to do, and in order to say more about the finer structure of the Ext groups $\Ext_\NB^\blt(\NB^{\ge\theta}e^\theta,L_0)^*=\Ext_\NB^\blt(\Delta_\theta^{\ol{l_\theta(\theta)}},L_0)^*=H^\blt(\NB^-:L_0)^*$ one would first have to discuss a different Ext group, namely $\Ext_\NB^\blt(\ol\Delta_\theta^{\ol{l_\theta(\theta)}},L_0)^*$. (More generally one would consider $\Ext_A^\blt(\bigoplus_{\lbd\in\theta}\ol\Delta_\lbd^{\ol{l_\lbd(\theta)}},M)^*$.) As these are Ext groups out of proper standard modules, we shall call them ``proper standard Ext groups'', or ``small Ext groups''. Of course, independently of the question of determining the nilHecke structure on the standard Ext groups, these proper standard Ext groups are interesting for homological reasons. In some sense they are harder to handle. If one knows a reasonably informative filtration of projectives over the Cartan (e.g. a Jacobson/Loewy filtration), one can then get a spectral sequence whose first page are proper standard Ext groups and which converges to (the graded pieces of) standard Ext groups. The amount of information you can get out of this depends, of course, on how informative the starting filtration was. In our case we have such a filtration, but instead of computing the Cartan structure of the standard Ext groups by using this spectral sequence, we will work backwards using what we have proved earlier to obtain results about proper standard Ext groups.

Let us say a few words about the differences between standard and proper standard Ext groups. Big and small Vermas inherit different properties of actual standard modules from a highest weight context: while standard modules inherit the property of being freely generated from a single vector, proper standard modules $\ol\Delta_\lbd$ inherit the property of containing exactly one copy of the simple $L_\lbd$ in its JH series, with every other simple having $\mu\in\phi>\theta$ (using lowest weight notation here). Other properties, such as BGG reciprocity and Ext vanishing, require both standard and proper standard objects to state. In particular proper standards lose the property of being freely generated on a single vector -- that vector, say of $L_\lbd(\theta)$, has the additional stipulation of being killed by the Jacobson radical of the Cartan algebra $A^\theta$. 

\subsection{From the nilHecke over $\Gamma$ to nilBrauer}
Recall that the Cartan algebras $A^\theta$ in our case are $\NB^\theta=\Gamma[X_1,\cdotsc,X_\theta,\sigma_1,\cdotsc,\sigma_{\theta-1}]/\sim$, the nilHecke over the ground ring $\Gamma\cong\BC[\tsl q_1,\tsl q_3,\cdots]$. As usual, this algebra is isomorphic to the matrix algebra over the symmetric polynomials in the $X$'s over $\Gamma$, which in turn is Morita-equivalent to just the ring of symmetric polynomials over $\Gamma$, which we denote 
\[\Sym_\Gamma(\theta)\coloneqq \Gamma[X_1,\cdotsc,X_\theta]^{S_\theta}=\Gamma[\tsl e_1,\cdotsc,\tsl e_\theta].\]
Note the usage of the slanted text font to denote symmetric polynomials\footnote{This is done in an effort to distinguish e.g. the elementary symmetric polynomials $\tsl e_i$ from the idempotents $e^\theta$, or the Schur $\tsl q$-functions $\tsl q_i$ from the quantum parameter $q$.}. This algebra has two natural gradings: the quantum grading, or the $q$-grading, in which\footnote{We have elected not to decorate the degree with a subscript here because this was the first degree one considered.}
\[\deg \tsl e_i=2i,\qquad \deg \tsl q_i=2i,\] 
and the symmetric grading, in which
\[\deg_\sym \tsl e_i=1,\qquad\deg_\sym \tsl q_i=1.\]
Note well that the ($q$-graded) algebra $\Sym_\Gamma(\theta)$ has (up to grading shift) a unique simple module, namely the one-dimensional $\BC$, appearing as the head. The Loewy filtration of this algebra as a module over itself, namely the projective indecomposable, is also very simple: it is the symmetric-degree filtration, given by
\[\Sym_{\Gamma}(\theta)=\Sym_{\Gamma}(\theta)_{\deg_\sym\ge 0}\supset \Sym_{\Gamma}(\theta)_{\deg_\sym\ge 1}\supset \Sym_{\Gamma}(\theta)_{\deg_\sym\ge 2}\supset\cdots,\]
so that one visibly has that $\Hed \Sym_\Gamma(\theta)\cong\BC$, as well as other things such as $\Sym_\Gamma(\theta)_{\deg_\sym=1}\coloneqq \Sym_{\Gamma}(\theta)_{\deg_\sym\ge 1}/\Sym_{\Gamma}(\theta)_{\deg_\sym\ge 2}\cong\BC\set{\tsl e_1,\cdotsc,\tsl e_\theta,\tsl q_1,\tsl q_3,\cdots}$, which has a (quantum-)dimension one can easily read off. For convenience, let
\[M_\theta(k)\coloneqq \Sym_\Gamma(\theta)_{\deg_\sym=k}.\]
Note well that the $q$-dimension of this is
\[\dim M_\theta(k)=[t^k] \frac{1}{(1-tq^{-2})(1-tq^{-4})\cdots(1-tq^{-2\theta})}\prod_{0< i\equiv2\mod 4}\frac{1}{1-tq^{-i}},\] 
where the $[t^k]$ in front indicates that we are taking the coefficient in front of $t^k$, and note $\prod_{0< i\equiv2\mod 4}\frac{1}{1-tq^{-i}}=\dim \Gamma$.

As $\Sym_\Gamma(\theta)$ is Morita-equivalent to $\NB^\theta$, we then have that the projective indecomposable of $\NB^\theta$, namely $P_\theta(\theta)=q^{\binom{\theta}{2}}\Gamma[X_1,\cdotsc,X_\theta]$, has a filtration $\F=\F_\te{sym}$
\[P_\theta(\theta)=\F^0P_\theta(\theta)\supset \F^1P_\theta(\theta)\supset\F^2P_\theta(\theta)\supset\cdots\] 
in which the quotients, which are semisimple, are
\[\gr_\F^kP_\theta(\theta)=L_\theta(\theta)\otimes_\BC M_\theta(k).\]

Now we use the fact that $\jota^\theta_!$ is exact to obtain from this a filtration for $\Delta_\theta$ in which the quotients are $\ol\Delta_\theta\otimes_\BC M_\theta(k)$:
\[\Delta_\theta=\F^0 \Delta_\theta\supset \F^1\Delta_\theta\supset \F^2\Delta_\theta\supset\cdots,\] 
with
\[\gr_\F^k\Delta_\theta=\ol\Delta_\theta\otimes_\BC M_\theta(k).\] 

\subsection{The spectral sequence}
The idea now is to compute $\Ext(\Delta_\theta^{\ol{l_\theta(\theta)}},L_0)^*$ in terms of $\Ext(\ol\Delta_\theta^{\ol{l_\theta(\theta)}},L_0)^*$ by using the filtration above to obtain a spectral sequence. 

From the filtration of $\Delta$'s by $\ol\Delta$'s above, one obtains a spectral sequence going from proper standard Ext groups to standard Ext groups. To be more precise, consider the complex $\hom_{\NB}(\sq,Q^\blt)^*$ for an injective resolution $L_0\lto Q^\blt$ of $L_0$. This is a covariant functor which is exact, and so each $$0\lto \F^{k+1}\Delta_\theta\lto \F^k \Delta_\theta \lto \ol\Delta_\theta\otimes_\BC M_\theta(k)\lto 0,$$ or rather each $$0\lto \F^{k+1}(\Delta_\theta^{\ol{l_\theta(\theta)}})\lto \F^k (\Delta_\theta^{\ol{l_\theta(\theta)}}) \lto \ol\Delta_\theta^{\ol{l_\theta(\theta)}}\otimes_\BC M_\theta(k)\lto 0,$$ under this functor gives rise to 
\[0\lto \F^{k+1}\hom_\NB(\Delta_\theta^{\ol{l_\theta(\theta)}},Q^\blt)^*\lto \F^k\hom_\NB(\Delta_\theta^{\ol{l_\theta(\theta)}},Q^\blt)^*\lto \hom_\NB(\ol\Delta_\theta^{\ol{l_\theta(\theta)}}\otimes_\BC M_\theta(k),Q^\blt)^*\lto 0.\]
So there is a filtration on the complex $\hom_\NB(\Delta_\theta^{\ol{l_\theta(\theta)}},Q^\blt)^*$ whose associated graded is $\hom_\NB(\ol\Delta_\theta^{\ol{l_\theta(\theta)}},Q^\blt)^*\otimes_\BC M_\theta(k)$, and by the theory of spectral sequences we have a spectral sequence
\[E_1^{p,q}=\Ext^{p+q}_\NB(\ol\Delta_\theta^{\ol{l_\theta(\theta)}},L_0)^*\otimes_\BC M_\theta(p)\implies E_\infty^{p,q}=\gr^p \Ext^{p+q}_\NB(\ol\Delta_\theta^{\ol{l_\theta(\theta)}},L_0)^*.\]
Pictorially, this spectral sequence roughly looks like:
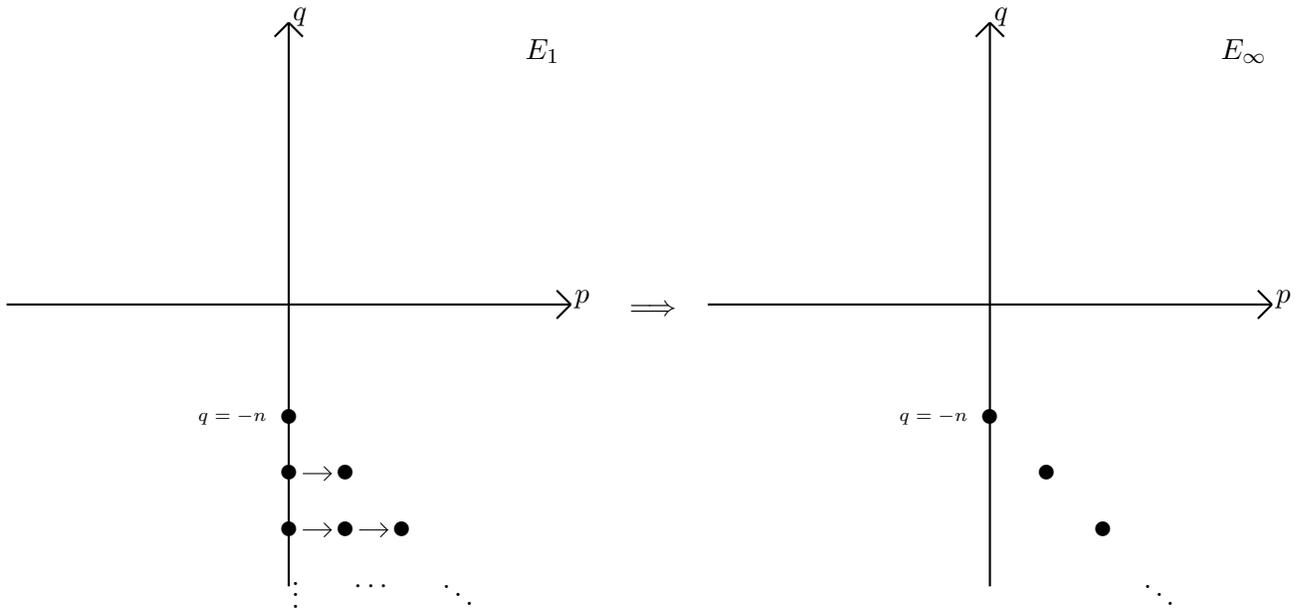
\[\hackcenter{\begin{tikzpicture}[scale=0.375]
    \draw[thick] (-10,0)--(10,0);
        \draw[thick] (10,0)--(9.5,0.5);
        \draw[thick] (10,0)--(9.5,-0.5);
        \node at (10.4,0.2) {$p$};
        \draw[thick] (0,-10)--(0,10);
        \draw[thick] (0,10)--(0.5,9.5);
        \draw[thick] (0,10)--(-0.5,9.5);
        \node at (0.4,10.2) {$q$};
        \node at (9,9) {$E_1$};
        \node (A) at (0,-4) {\Large $\blt$};
        \node at (-2,-4) {\tiny $q=-n$};
        \node (B) at (2,-6) {\Large $\blt$};
        \node (C) at (0,-6) {\Large $\blt$};
        \node (D) at (4,-8) {\Large $\blt$};
        \node (E) at (2,-8) {\Large $\blt$};
        \node (F) at (0,-8) {\Large $\blt$};
        \node at (0.25,-10) {$\vdots$};
        \node at (6,-10) {$\ddots$};
        \node at (3,-10) {$\cdots$};
        \draw (0.5,-6)--(1.5,-6)--(1.25,-6.25);
        \draw (1.5,-6)--(1.25,-5.75);
        \draw (0.5,-8)--(1.5,-8)--(1.25,-8.25);
        \draw (1.5,-8)--(1.25,-7.75);
        \draw (2.5,-8)--(3.5,-8)--(3.25,-8.25);
        \draw (3.5,-8)--(3.25,-7.75);
\end{tikzpicture}}
\specseqimplies
\hackcenter{\begin{tikzpicture}[scale=0.375]
    \draw[thick] (-10,0)--(10,0);
        \draw[thick] (10,0)--(9.5,0.5);
        \draw[thick] (10,0)--(9.5,-0.5);
        \node at (10.4,0.2) {$p$};
        \draw[thick] (0,-10)--(0,10);
        \draw[thick] (0,10)--(0.5,9.5);
        \draw[thick] (0,10)--(-0.5,9.5);
        \node at (0.4,10.2) {$q$};
        \node at (9,9) {$E_\infty$};
        \node (A) at (0,-4) {\Large $\blt$};
        \node at (-2,-4) {\tiny $q=-n$};
        \node (B) at (2,-6) {\Large $\blt$};
        \node (D) at (4,-8) {\Large $\blt$};
        \node at (6,-10) {$\ddots$};
\end{tikzpicture}}
\]

Let us check the row $q=n$. By inspection, it must be the case that
\[\Ext^n_\NB(\ol\Delta_\theta^{\ol{l_\theta(\theta)}},L_0)^*=\gr^0\Ext^n(\Delta_\theta^{\ol{l_\theta(\theta)}},L_0)^*.\] 
Note that by tensor-hom one has
\[\rhom_{\NB^\theta}(L_\theta(\theta)^{ \ol{l_\theta(\theta)}},\rhom_\NB(\NB^{\ge\theta}e^\theta,L_0))=\rhom_{\NB}(\NB^{\ge\theta}e^\theta\otimes_{\NB^\theta}L_\theta(\theta)^{ \ol{l_\theta(\theta)}},L_0)=\rhom_\NB(\ol\Delta_{\theta}^{\ol{l_\theta(\theta)}},L_0),\]
and by the Grothendieck spectral sequence, since we have already since $\rhom_\NB(\NB^{\ge\theta }e^\theta,L_0)$ is appropriately concentrated, we know that the last equality in the following is true:
\begin{align*}
    \Hed\Ext_\NB^n(\NB^{\ge\theta}e^\theta,L_0)^*&= \hom_{\NB^\theta}(\Ext_\NB^n(\NB^{\ge\theta}e^\theta,L_0)^*,L_\theta(\theta)^{\ol{l_\theta(\theta)}})\\
    &=\hom_{\NB^\theta}(L_\theta(\theta)^{\ol{l_\theta(\theta)}},\Ext_\NB^n(\NB^{\ge\theta}e^\theta,L_0))\\
    &=\Ext^n_\NB(\ol\Delta_\theta^{\ol{l_\theta(\theta)}},L_0),
\end{align*}
and the head of $\Ext^{n,*}$ as we have seen above already is simply $L_\theta(\theta)$. Hence
\begin{PROP}
    At $t=0$, letting $\theta=2n$, the nilHecke action on the proper standard Ext groups of the trivial representation is
    \[\Ext^n_{\NB}(\ol\Delta_\theta^{\ol{l_\theta(\theta)}},L_0)^*\cong L_\theta(\theta).\] 
    Equivalently the claim is that as vector spaces
    \[\Ext^n_\NB(\ol\Delta_\theta,L_0)\cong\BC.\] 
\end{PROP}

\bibliographystyle{alpha}
\bibliography{references}

\

\noindent \small\textsc{Department of Mathematics, Columbia University}

\noindent\small\textsc{Email:} \href{mailto:fz2326@columbia.edu}{{\tt{\textcolor{black}{fz2326@columbia.edu}}}}


\end{document}